\documentclass[10pt,reqno,twoside]{amsart}

\usepackage{amsfonts}
\usepackage{paralist}
  \usepackage{graphics}
  \usepackage{epsfig}
 \usepackage[colorlinks=true]{hyperref}
\usepackage{hyperref}

\numberwithin{equation}{section}
\usepackage{latexsym}
\usepackage{amsmath}
\usepackage{amssymb}
\usepackage{mathrsfs}
\usepackage{graphicx}
\usepackage{ifthen}
\usepackage[T1]{fontenc}
\usepackage[latin1]{inputenc}
\usepackage{color}
\newtheorem{Satz}{Theorem}[section]

\newtheorem{Def}[Satz]{Definition}

\newcommand{\oo}{\overline \Omega}

\newcommand{\po}{\partial\Omega}

\newcommand{\dist}{\text {dist}}

\newcommand{\diag}{\text {diag}}

\newcommand{\diam}{\text {diam}}

\newcommand{\supp}{\text {supp}}

\newcommand{\sign}{\text {sign}}

\begin{document}

\title[Boundary separated and clustered layer positive solutions]
{
Boundary separated  and   clustered layer positive solutions  for an elliptic  Neumann problem with large exponent
}

\author[Yibin Zhang]{Yibin Zhang}
\address{College of Sciences, Nanjing Agricultural University, Nanjing 210095, China}
\email{yibin10201029@njau.edu.cn}

\subjclass[2010]{Primary 35B25; Secondary 35B38,  35J25.}

\keywords{
Boundary separated and clustered layer positive solutions;
Large exponent;
Anisotropic Green's function.}

\begin{abstract}
Given a  smooth bounded domain $\mathcal{D}$ in $\mathbb{R}^N$ with $N\geq3$,
we study the existence and the profile of positive solutions for
the following elliptic Nenumann problem
$$
\left\{\aligned
&-\Delta
\upsilon+\upsilon=\upsilon^p,\,\,\,\ \,
\upsilon>0
\,\,\,\ \,
\textrm{in}\,\,\,\,\,\,
\mathcal{D},\\[1mm]
&\frac{\partial \upsilon}{\partial\nu}=0\,\,\,\,\,
\ \,\ \,\qquad\,\ \,\,\,
\,\ \,\,\quad\quad\,\,\
\textrm{on}\,\,\,\,
\partial\mathcal{D},
\endaligned\right.
$$
where  $p>1$ is a large exponent and $\nu$
denotes the outer unit normal vector to the boundary $\partial\mathcal{D}$.
For suitable domains $\mathcal{D}$,  by a constructive way
we prove that, for any non-negative integers $k$, $l$ with  $k+l\geq1$,
if $p$ is large enough,
such a problem  has a family of positive solutions with  $k$
boundary layers and $l$ interior layers which concentrate along  $k+l$ distinct $(N-2)$-dimensional minimal
submanifolds  of $\partial\mathcal{D}$, or collapse to
the same $(N-2)$-dimensional minimal submanifold of $\partial\mathcal{D}$ as $p\rightarrow+\infty$.\\
\end{abstract}

\maketitle

\section{Introduction}
In this paper, we are interested in
the following classical elliptic Neumann problem
\begin{equation}\label{a1}
\left\{\aligned
&-d^2\Delta
\upsilon+\upsilon=\upsilon^p,\,\,\,\ \,
\upsilon>0
\,\,\,\ \,
\textrm{in}\,\,\,\,\,\,
\mathcal{D},\\[1mm]
&\frac{\partial \upsilon}{\partial\nu}=0\,\,\,\,\,
\ \,\ \,\quad\qquad\,\,\,\,\,
\,\ \,\,\quad\quad\,\,\
\textrm{on}\,\,\,\,
\partial\mathcal{D},
\endaligned\right.
\end{equation}
where $\mathcal{D}$ is a smooth bounded domain in
$\mathbb{R}^N$,   $N\geq2$,
$d>0$ is a  parameter, $p>1$ is  a exponent
and $\nu$
denotes the outer unit normal vector to
the boundary $\partial\mathcal{D}$.

Problem (\ref{a1}) has received considerable
attention in the last three decades, because it
appears in many different mathematical
models: for example,  it arises from the study of
steady states for the logarithmic Keller-Segel system in
chemotaxis \cite{KS} and the shadow system of
the Gierer-Meinhardt model
in biological pattern formation \cite{GM}.
In particular, it has been shown that the
solutions of (\ref{a1}) exhibit a variety of
interesting
concentration phenomena as either the exponent $p$ tends to
some  critical values or the parameter $d$ tends
to zero.

Let us first define $p_{h+1}^*:=(N-h+2)/(N-h-2)$,
$0\leq h\leq N-3$
as the {\it $(h+1)$-th  critical exponent} (recall that $h$ is an integer)
and set $p_{N-1}^*:=+\infty$.
In the subcritical case, i.e. $N=2$, or $N\geq3$ and $p<p_1^*=(N+2)/(N-2)$,
compactness of Sobolev's  embedding ensures the existence of a positive least
energy solution  of (\ref{a1}).
For $d>0$ sufficiently small,  Lin, Ni and Takagi in
 \cite{LNT,NT1,NT2}  proved that this  least energy solution
has exactly one bounded, very sharp spike located on the boundary and
 near the {\it most curved} part of the boundary,
i.e., the region where the mean curvature of the boundary
attains  its maximum.
Higher energy solutions of (\ref{a1}) with multiple boundary peaks
as well as multiple interior
peaks have been established  in \cite{CW,DFW,DFW1,DY,GPW,GW,GW1,GWW,K,L,W1,WW}.
It turns out that multiple boundary spikes tend
to cluster around  critical points of the mean curvature
of the boundary, while the
location of multiple interior spikes is determined by the distance between the peaks and
the boundary.
In particular,  Gui and Wei  in \cite{GW1} proved that
 for any non-negative integers $k$, $l$ with $k+l\geq1$,
problem (\ref{a1}) has a solution with exactly
$k$ different boundary spikes and $l$ different interior spikes provided that $d$ is small and $p$
is subcritical.
Generally, such spiky solutions are called
{\it solutions with $0$-dimensional concentration sets}.

In the critical case, i.e. $p=p_1^*$, such type of concentration phenomena
occurs, but the situation is quite different.
The lack of compactness of Sobolev's embedding makes it non-obvious to apply
variational techniques to obtain a nonconstant least energy solution of (\ref{a1})
when $d$ is sufficiently small. The first existence result for a
nonconstant least energy solution of (\ref{a1}) in general domains and $d$
small were obtained by Adimurthi-Mancini \cite{AM}
and Wang \cite{W}. The profile and asymptotic behavior of this least energy
solution has been  clarified in the subsequent works \cite{APY,NPT,R}.
As in the subcritical case, the least energy solution has a unique maximum point
or peak that lies on the boundary and goes, as $d$ tends to zero,
to a maximum point of the  mean curvature of the boundary.
Unlike the subcritical case, the least energy
solution blows up, as $d$ tends to zero,
at a maximum point of the mean curvature of  the boundary.
Higher energy solutions with  one or more separated or clustered
boundary blow-up points have been exhibited
for instance in \cite{AM1,AMY,G,GG,LWW,R,R1,WY}, where these blow-up points
are nothing but degenerate or non-degenerate
critical points of the mean curvature of the boundary.
A major difference with the subcritical case is that the condition of
positivity for the mean curvature at these critical points
turns out to be necessary for the boundary bubbling phenomenon to take place
\cite{GL,R}. Note that in striking contrast with the subcritical case, there
are no solutions blowing up at only interior points when $d$ tends
to zero, namely at least one blow-up point has to lie on the boundary,
even  all of the blow-up points for solutions with uniformly bounded energy
have to lie only on the boundary
as established in \cite{CNY,E,R2}. Moreover, for the subcritical case,
spiky solutions with only interior peaks always exist, but these interior peaks
must stay in the domain  as  $d$ tends
to zero \cite{W2}; while for the critical case, interior peaks of solutions
could not remain in the domain as $d$ tends
to zero \cite{G,GGZ}.

The {\it almost first critical case}, i.e. $p=p_1^*\pm\epsilon$
with $\epsilon$  positive and sufficiently small,
has been widely studied. In the slightly subcritical case, i.e.
$p=p_1^*-\epsilon$,
the existence of solutions with simple or non-simple
blow-up points located on the boundary near
critical points
of the mean curvature of the boundary  with {\it negative value}
as $\epsilon\rightarrow0$ for fixed $d$ and $N\geq4$ was established
in \cite{CN,RW1,W3};
while  in the slightly supercritical case, i.e.
$p=p_1^*+\epsilon$, it was proved in \cite{DMP,RW1,W3} that  there also exists a solution
with simple or non-simple blow-up points located on the boundary near
critical points of the mean curvature of the boundary  with {\it positive value}
as $\epsilon\rightarrow0$ for fixed $d$ and $N\geq4$.
In \cite{RW}, it was proved that if $N=3$ and $\epsilon$ tends to zero from below
or above, then a solution with one interior blow-up point may exist for finite $d$.
Moreover, if $N=2$,
$d$ is fixed and the exponent $p$ goes to $+\infty$,
Musso and Wei in \cite{MW} proved that  for any non-negative integers $k$, $l$ with
$k+l\geq1$, problem (\ref{a1}) admits  a solution with exactly
$k$ different boundary spikes and $l$ different interior spikes, whose location can be characterized by
critical points of a certain combination of Green's function and its regular part.

It seems natural to ask if problem (\ref{a1}) has {\it solutions that
exhibit concentration phenomena  on
$h$-dimensional  subsets of
$\overline{\mathcal{D}}$  for every $1\leq h\leq N-1$},
as conjectured by Ni in \cite{N}.
In particular, given an $h$-dimensional submanifold $\Gamma$ of
$\partial\mathcal{D}$ and assuming that either $h\geq N-2$ or
$p\leq p_{h+1}^{*}+\epsilon$ with $\epsilon$  positive and sufficiently small,
one question is to ask whether problem (\ref{a1}) admits a solution
that concentrates along $\Gamma$ as either $d$ or $\epsilon$
tends to zero.
For results in this direction, we first mention
the {\it $(h+1)$-th subcritical case}, i.e. $p<p_{h+1}^*$, where
Malchiodi and Montenegro \cite{MM2,MM3}
proved that, given any $N\in\mathbb{N}$ and any $p>1$, there
exist solutions  concentrating on the
{\it whole boundary} if the sequence $d$ satisfies
some gap condition, corresponding to $h=N-1$.
Furthermore, the result was  extended in \cite{MM1,M1} for
general $h$, and the concentration set $\Gamma$ is
{\it an  embedded closed
minimal submanifold of $\partial\mathcal{D}$}
which is in addition {\it nondegenerate } in the sense that
its Jacobi operator is invertible.  Indeed, this phenomenon is
rather subtle compared with pointwise concentration:
existence can only be achieved along a sequence
 of parameters $d=d_i\rightarrow0$.
The sequence of  parameters $d_i$ must be suitable away from certain values where
resonance occurs, and the topological type of the solution changes:
unlike the pointwise concentration, the Morse index of these solutions
is very large and grows as $d_i\rightarrow0$.
Del Pino, Mahmoudi and Musso \cite{DMM} extended this type of results to the
{\it $(h+1)$-th critical case
 $p=p_{h+1}^*$} and  proved that if  $\Gamma$ is an  embedded closed
minimal submanifold of $\partial\mathcal{D}$ with dimension $h\leq N-7$
(in particular, $N\geq8$)
which is nondegenerate, and a certain weighted average of sectional curvature of
$\partial\mathcal{D}$ is positive along $\Gamma$, then problem (\ref{a1})
admits a solution for a suitable sequence of parameters $d_i\rightarrow0$
which blows up along $\Gamma$. Further for the {\it
almost $(h+1)$-th critical case} $p=p_{h+1}^*\pm\epsilon$
with $0<\epsilon\rightarrow0$  but $d$ fixed,
the same conclusion was established by Deng, Mahmoudi and Musso
\cite{DMM1} under
analogous assumptions. Meanwhile,
for any integer $h=1,\ldots,N-3$,
Manna and Pistoia \cite{MP} proved that in some suitable domains $\mathcal{D}$,
problem (\ref{a1}) has a solution which blows up along an $h$-dimensional minimal
submanifold $\Gamma$ of $\partial\mathcal{D}$ as $p$ approaches  from either below or
above the $(h+1)$th critical exponent $p_{h+1}^*$.

In the present paper we consider the {\it almost $(N-1)$-th critical case},
i.e. $h=N-2$, and give a positive answer when $d$ is fixed and
$p$  goes to $+\infty$. More precisely, we find some domains $\mathcal{D}$
such that if $d=1$ and $p$ is large enough, then for any positive
integer $m$, problem (\ref{a1}) has a positive solution with $m$ distinct mixed interior and
boundary layers which  concentrate along  $m$ distinct $(N-2)$-dimensional minimal
submanifolds  of $\partial\mathcal{D}$, or collapse to
the same $(N-2)$-dimensional minimal submanifold of $\partial\mathcal{D}$ as $p$ goes to $+\infty$.

Let  $\Omega$  be  a smooth bounded domain in $\mathbb{R}^2$ such that
\begin{equation*}\label{1.4}
\aligned
\overline{\Omega}\subset
\{
(x_1,x_2)\in\mathbb{R}^2|\,
\,x_1>0
\}
\,\,\qquad\,\,
\textrm{or}
\,\,\,\qquad\,
\overline{\Omega}\subset
\{
(x_1,x_2)\in\mathbb{R}^2|\,
\,x_1,\,x_2>0
\}.
\endaligned
\end{equation*}
Let  $n=1$ or $n=2$ be fixed.
Fix $k_1,k_n\in\mathbb{N}$ with $h:=k_1+k_n=N-2$ and set
\begin{equation}\label{1.5}
\aligned
\mathcal{D}:=
\big\{
(y_1,y_n,x')\in\mathbb{R}^{k_1+1}\times\mathbb{R}^{k_n+1}\times\mathbb{R}^{2-n}|\,
\,(|y_1|,|y_n|,x')\in\Omega
\big\}.
\endaligned
\end{equation}
Then $\mathcal{D}$ is a smooth bounded domain in $\mathbb{R}^N$ which is  $\Upsilon$-invariant
for the action of the group $\Upsilon:=\mathcal{O}(k_1+1)\times \mathcal{O}(k_n+1)$
on $R^N$ given by
$$
\aligned
(g_1,g_n)(y_1,y_n,x'):=
(g_1y_1,g_ny_n,x').
\endaligned
$$
Here $\mathcal{O}(k_i+1)$ denotes the group of linear isometries of $\mathbb{R}^{k_i+1}$.
For $p>1$ large enough we shall look for $\Upsilon$-invariant solutions of problem (\ref{a1})
with $d=1$, i.e. solutions $\upsilon$ of the form
\begin{equation}\label{1.6}
\aligned
\upsilon(y_1,y_n,x')=u(|y_1|,|y_n|,x').
\endaligned
\end{equation}
Then a simple calculation shows that $\upsilon$ solves problem (\ref{a1})
with $d=1$ if and only if $u$ satisfies
\begin{equation}\label{1.7}
\left\{\aligned
&-\Delta u-\sum_{i=1}^n\frac{k_i}{\,x_i\,}\frac{\partial u}{\partial x_i}+u=u^p,
\,\,\,\,\,
u>0\,\,\,\,\,\,
\textrm{in}\,\,\,\,\,
\Omega,\\
&\frac{\partial u}{\partial\nu}=0\,\,\,\,\,
\ \,\ \,\qquad\qquad\qquad\qquad\,\,\,\,
\,\ \,\,\quad\quad\,\,\
\textrm{on}\,\,\,\,
\partial\Omega.
\endaligned\right.
\end{equation}
Thus, we are led to study the more general anisotropic problem
\begin{equation}\label{1.1}
\left\{\aligned
&-\nabla(a(x)\nabla u)+a(x)u=a(x)u^p,\,\,\,\,
u>0\,\,\,\,\,
\textrm{in}\,\,\,\,\,
\Omega,\\[1mm]
&\frac{\partial u}{\partial\nu}=0\,\,
\qquad\quad\qquad\qquad\qquad\qquad\qquad
\ \ \ \ \,\,
\textrm{on}\,\,\,
\partial\Omega,
\endaligned\right.
\end{equation}
where  $\Omega$ is a smooth  bounded   domain in $\mathbb{R}^2$,
$a(x)$ is a positive smooth function over $\overline{\Omega}$,
$p>1$ is a large exponent and $\nu$ denotes the outer unit  normal vector to  $\partial\Omega$.
Note that if
\begin{equation*}\label{1.8}
\aligned
a(x)=a(x_1,x_n,x'):=x_1^{k_1}x_n^{k_n},
\endaligned
\end{equation*}
then problem (\ref{1.1}) can be rewritten as equation (\ref{1.7}).

Our goal is to construct solutions $u_p$ to problem (\ref{1.1}) with $m$
distinct mixed interior and boundary spikes which concentrate at points
$\xi^*_1,\ldots,\xi^*_m$ of $\partial\Omega$, or accumulate to the same
point $\xi_*$ of $\po$ as $p$ goes to $+\infty$. They
correspond, via (\ref{1.6}), to $\Upsilon$-invariant solutions $\upsilon_p$
of problem (\ref{a1})$\large|_{d=1}$ with $m$ distinct mixed
interior and boundary layers which concentrate along the $\Upsilon$-orbits
$\Xi(\xi^*_1),\ldots,\Xi(\xi^*_m)$ of $\partial\mathcal{D}$, or collapse to
the same $\Upsilon$-orbit $\Xi(\xi_*)$ of $\partial\mathcal{D}$
as $p$ goes to $+\infty$. Here
$$
\aligned
\Xi(\zeta):=\big\{
(y_1,y_n,x')\in\partial\mathcal{D}|
\,\,
(|y_1|,|y_n|,x')=\zeta\in\partial\Omega
\big\}
\endaligned
$$
is a $(N-2)$-dimensional minimal submanifold of $\partial\mathcal{D}$
diffeomorphic to $\mathbb{S}^{k_1}\times\mathbb{S}^{k_n}$ (note that $k_1+k_n=h=N-2$), where
$\mathbb{S}^{k_i}$ is the unit sphere in $\mathbb{R}^{k_i+1}$.

Let us define the linear differential operator
$$
\aligned
\Delta_au=\frac1{a(x)}\nabla(a(x)\nabla u)=\Delta u+\nabla\log a(x)\nabla u
\endaligned
$$
and the   Green's function associated with the Neumann problem
\begin{equation}\label{1.2}
\left\{\aligned
&-\Delta_aG(x,y)+G(x,y)=\delta_y(x),\,\,\,\,\,\,\,
x\in\Omega,\\
&\frac{\partial G}{\partial\nu_x}(x,y)=0,
\qquad\,\,
\qquad\qquad\qquad\,\,
x\in\partial\Omega,
\endaligned\right.
\end{equation}
for every $y\in\oo$.
The regular part of $G(x,y)$ is defined depending on whether $y$ lies in
the domain or on its boundary as
\begin{equation}\label{1.3}
\aligned
H(x,y)=\left\{\aligned
&G(x,y)+\frac{1}{2\pi}\log|x-y|,\,\quad\,y\in\Omega,\\[1mm]
&G(x,y)+\frac1{\pi}\log|x-y|,\,\ \quad\,y\in\po.
\endaligned\right.
\endaligned
\end{equation}

Our  first result
concerns the existence of
solutions of problem (\ref{1.1})
whose interior and boundary spikes are uniformly far away from each other
and interior spikes
lie in the domain with  distance
to the boundary uniformly approaching zero.

\vspace {1mm}
\vspace {1mm}
\vspace {1mm}

\noindent{\bf Theorem 1.1.}\,\,\,{\it
Let $k$, $l$  be non-negative integers with $m=k+l\geq1$ and assume that there exist
$m$ different points  $\xi^*_1,\ldots,\xi^*_m\in\partial\Omega$
such that
each $\xi_i^*$ is either a strict local
maximum or a strict local minimum point of $a(x)$ on $\po$  and
satisfies  for all $i=1,\ldots,l$,
$\partial_{\nu}a(\xi_i^*):=\langle\nabla a(\xi_i^*),\,\nu(\xi_i^*)\rangle>0$.
Then, there exists $p_0$ such that for any $p>p_0$, there is a
family of  positive solutions $u_p$ for problem {\upshape (\ref{1.1})}
with $k$ different boundary spikes and $l$ different interior spikes
located at distance $O\left(1/p\right)$ from $\po$. More precisely,
$$
\aligned
u_{p}(x)=\sum\limits_{i=1}^{m}\frac{1}{\gamma\mu_i^{2/(p-1)}}\left[\,\log
\frac1{(\delta_i^2+|x-\xi^p_i|^2)^2}
+c_i H(x,\xi^p_i)\right]+o\left(1\right),
\endaligned
$$
where $o(1)\rightarrow0$, as $p\rightarrow+\infty$, on each compact subset of
$\oo\setminus\{\xi^p_1,\ldots,\xi^p_m\}$,
the parameters $\gamma$, $\delta_i$ and $\mu_i$ satisfy
$$
\aligned
\gamma=p^{p/(p-1)}\varepsilon^{2/(p-1)},\,\,
\quad\quad\quad\,\,
\delta_i=\mu_i\varepsilon,
\quad\quad\quad\,\,
\varepsilon=e^{-p/4},
\quad\quad\quad\,\,
\frac1{C}<\mu_i<Cp^{2(m^2+1)},
\endaligned
$$
for some $C>0$,
 $(\xi^p_1,\ldots,\xi^p_m)\in\Omega^l\times(\po)^{k}$ satisfies
$$
\aligned
\xi^p_i\rightarrow\xi^*_i
\,\quad\,\textrm{for all}\,\,\,i,
\,\qquad\quad\,
\textrm{and}
\,\qquad\quad\,
\dist(\xi^p_i,\po)=O\left(1/p\right)
\quad\,\,\forall\,\,i=1,\ldots,l,
\endaligned
$$
and $c_i=8\pi$ for $i=1,\ldots,l$, but $c_i=4\pi$
for $i=l+1,\ldots,m$. In particular, for any  $d>0$,
as $p\rightarrow+\infty$,
$$
\aligned
pu_{p}^{p+1}\rightharpoonup8\pi e\sum_{i=1}^l\delta_{\xi^*_i}+4\pi e\sum_{i=l+1}^m\delta_{\xi^*_i}
\,\quad\textrm{weakly in the sense of measures in}\,\,\,\,
\overline{\Omega},
\endaligned
$$
$$
\aligned
u_{p}\rightarrow0
\,\quad\textrm{uniformly in}\,\,\,\,
\overline{\Omega}\setminus\bigcup_{i=1}^mB_{d}(\xi^*_i),
\endaligned
$$
and
$$
\aligned
\sup\limits_{x\in\oo\cap B_{d}(\xi^*_i)}u_p(x)\rightarrow\sqrt{e}.
\endaligned
$$
}

The corresponding result for problem  (\ref{a1})$\large|_{d=1}$ can be stated as follows.

\vspace{1mm}
\vspace{1mm}
\vspace{1mm}
\vspace{1mm}

\noindent{\bf Theorem 1.2.}\,\,\,{\it
Let $k$, $l$ be non-negative integers with $k+l\geq1$
 and $\mathcal{D}$ be as in {\upshape (\ref{1.5})}.
If the assumption of Theorem 1.1 holds, then
there exists $p_0$ such that for any $p>p_0$,
problem {\upshape (\ref{a1})}$\large|_{d=1}$
has a positive solution $\upsilon_p$
with $k$ boundary layers and $l$ interior layers which concentrate along $k+l$ different
$(N-2)$-dimensional minimal submanifolds of $\partial\mathcal{D}$,
namely the $\Upsilon$-orbit $\Xi(\xi^*_i)$ of $\xi^*_i$ for every $i=1,\ldots,k+l$, as $p\rightarrow+\infty$.
}

\vspace{1mm}
\vspace{1mm}
\vspace{1mm}
\vspace{1mm}

Our  next  result
concerns the existence of solutions of problem (\ref{1.1})
with mixed interior and boundary spikes which  accumulate to the same point
of the boundary.

\vspace{1mm}
\vspace{1mm}
\vspace{1mm}
\vspace{1mm}

\noindent{\bf Theorem 1.3.}\,\,\,{\it
Let $k$, $l$ be non-negative integers with $m=k+l\geq1$
and assume that $\xi_*\in\po$ is a strict local maximum point of $a(x)$
and satisfies $\partial_{\nu}a(\xi_*):=\langle\nabla a(\xi_*),\,\nu(\xi_*)\rangle=0$.
Then, there exists $p_0$ such that for any $p>p_0$, there is a
family of positive solutions $u_p$ for problem {\upshape (\ref{1.1})}
with $k$ different boundary spikes and $l$ different interior spikes
which accumulate to $\xi_*$ as $p\rightarrow+\infty$. More precisely,
$$
\aligned
u_{p}(x)=\sum\limits_{i=1}^{m}\frac{1}{\gamma\mu_i^{2/(p-1)}}\left[\,\log
\frac1{(\delta_i^2+|x-\xi^p_i|^2)^2}
+c_i H(x,\xi^p_i)\right]+o\left(1\right),
\endaligned
$$
where $o(1)\rightarrow0$, as $p\rightarrow+\infty$, on each compact subset of
$\oo\setminus\{\xi^p_1,\ldots,\xi^p_m\}$,
the parameters $\gamma$, $\delta_i$ and $\mu_i$ satisfy
$$
\aligned
\gamma=p^{p/(p-1)}\varepsilon^{2/(p-1)},\,\,
\quad\quad\quad\,\,
\delta_i=\mu_i\varepsilon,
\quad\quad\quad\,\,
\varepsilon=e^{-p/4},
\quad\quad\quad\,\,
\frac1{C}<\mu_i<Cp^{2(m^2+1)},
\endaligned
$$
for some $C>0$,
$(\xi^p_1,\ldots,\xi^p_m)\in\Omega^l\times(\po)^{k}$  satisfies
$$
\aligned
\xi^p_i\rightarrow\xi_*
\,\quad\,\textrm{for all}\,\,\,i,
\qquad
|\xi^p_i-\xi^p_k|>p^{-2(m^2+1)}
\quad\forall\,\,i\neq k,
\qquad
\textrm{and}
\qquad
\dist(\xi^p_i,\po)>p^{-2(m^2+1)}
\quad\forall\,\,i=1,\ldots,l,
\endaligned
$$
 and $c_i=8\pi$ for $i=1,\ldots,l$, but $c_i=4\pi$
for $i=l+1,\ldots,m$.
In particular, for any  $d>0$,
as $p\rightarrow+\infty$,
$$
\aligned
pu_{p}^{p+1}\rightharpoonup4\pi e(k+2l)\delta_{\xi_*}
\,\quad\textrm{weakly in the sense of measures in}\,\,\,\,
\overline{\Omega},
\endaligned
$$
$$
\aligned
u_{p}\rightarrow0
\,\quad\textrm{uniformly in}\,\,\,\,
\overline{\Omega}\setminus B_{d}(\xi_*),
\endaligned
$$
and
$$
\aligned
\sup\limits_{x\in\oo\cap B_{d/p^{4(m^2+1)}}(\xi_i^p)}u_p(x)\rightarrow\sqrt{e}.
\endaligned
$$
}

The corresponding result for problem {\upshape (\ref{a1})}$\large|_{d=1}$ can be stated as follows.

\vspace{1mm}
\vspace{1mm}
\vspace{1mm}
\vspace{1mm}

\noindent{\bf Theorem 1.4.}\,\,\,{\it
Let $k$, $l$ be non-negative integers with $k+l\geq1$
 and $\mathcal{D}$ be as in {\upshape (\ref{1.5})}.
If the assumption of Theorem 1.3 holds, then
there exists $p_0$ such that for any $p>p_0$,
problem {\upshape (\ref{a1})}$\large|_{d=1}$ has a positive solution $\upsilon_p$
with $k$ boundary layers and $l$ interior layers which collapse to the same
$(N-2)$-dimensional minimal submanifold of $\partial\mathcal{D}$, namely
the $\Upsilon$-orbit $\Xi(\xi_*)$ of $\xi_*$, as $p\rightarrow+\infty$.
}

\vspace{1mm}
\vspace{1mm}
\vspace{1mm}
\vspace{1mm}


Let us remark that the assumptions in
Theorems $1.3$-$1.4$  contain the following two cases:\\
(A1) \,
$\xi_*\in\po$ is a strict
local maximum point of $a(x)$ restricted on $\po$;\\
(A2) \,
$\xi_*\in\po$ is a strict local maximum
point of $a(x)$ restricted in $\Omega$
and satisfies $\partial_{\nu}a(\xi_*)=\langle\nabla a(\xi_*),\,\nu(\xi_*)\rangle=0$.\\
In fact, arguing as for the proof of Theorem $1.3$, we can easily find that if
(A1) holds,  then problem (\ref{1.1}) has positive solutions with
arbitrarily many
boundary spikes which  accumulate to $\xi_*$ along $\po$; while
if (A2) holds,  then problem (\ref{1.1}) has positive solutions with
arbitrarily many
interior spikes which  accumulate to $\xi_*$ along the
inner normal direction of $\po$.
For the latter  case, it seems that this paper is the {\it first} one in the literature
obtaining this type of concentration phenomenon for  {\it positive} solutions of
some two-dimensional anisotropic nonlinear elliptic Neumann problems, see \cite{AP} as an instance.

The general strategy for proving  our main results  relies on a very
well known  Lyapunov-Schmidt reduction.
In Section $2$ we provide an appropriate
approximation for a solution of problem (\ref{1.1}) and
give a basic estimate for the scaling error term created by the choice of our approximation.
Then we rewrite  problem (\ref{a1}) in terms of a linearized operator for which
a solvability theory, subject to suitable orthogonality conditions,
is performed  through solving a linearized  problem in Section $3$.
In Section $4$ we
solve  an auxiliary nonlinear problem.
In Section $5$ we reduce the problem of finding  spike solutions
of (\ref{1.1}) to that of finding a critical point of a finite-dimensional function.
Section $6$ concerns with  an asymptotic expansion for  the finite-dimensional function
appeared in Section $5$. Finally, in Section $7$ we provide the detailed proof
of Theorems $1.1$ and $1.3$.

In this paper, the letter $C$ will always denote
a generic positive  constant independent of $p$,
which could be changed from one line to another.
The symbol $o(t)$ (respectively $O(t)$) will denote a quantity for which
$\frac{o(t)}{|t|}$ tends to zero
(respectively, $\frac{O(t)}{|t|}$ stays bounded )
as parameter $t$ goes to zero.
Moreover, we will use the notation
$o(1)$ (respectively $O(1)$)
to stand for a quantity which tends to zero
 (respectively, which remains uniformly bounded) as $p\rightarrow+\infty$.

\section{An approximation for the solution}
In this section we provide an appropriate
approximation for a solution of problem (\ref{1.1}) and
give a basic estimate for the scaling error term created by the choice of our approximation.
Since the function $H(x,y)$ defined in (\ref{1.3}) plays an essential role in
our construction,  we shall first state its asymptotic behavior without proof, see \cite{AP} for details.

Consider the vector function $T(x)=(T_1(x),T_2(x))$ as the solution of
\begin{equation}\label{2.1a}
\aligned
\Delta_xT-T=\frac{x}{|x|^2}\,\,\quad\,\,\textrm{in}
\,\,\,\,\mathbb{R}^2,\,\,\quad\,\,T(x)\in L^{\infty}_{loc}(\mathbb{R}^2).
\endaligned
\end{equation}
Then standard elliptic regularity theory implies that for any $1<q<2$,
$T(x)\in W^{2,q}_{loc}(\mathbb{R}^2)\cap C^{\infty}(\mathbb{R}^2\setminus\{0\})$,
and the Sobolev embeddings  yield that
$T(x)\in W^{1,1/\alpha}(B_r(0))\cap C^{\alpha}(\overline{B_r(0)})$
for any $r>0$ and $0<\alpha<1$.

\vspace{1mm}
\vspace{1mm}
\vspace{1mm}
\vspace{1mm}

\noindent{\bf Lemma 2.1 (\cite{AP}).}\,\,{\it
Let $T(x)$ be the function described in {\upshape (\ref{2.1a})}.
There exists a function $H_1(x,y)$ such that\\
\indent {\upshape (i)}  for every $x,y\in\oo$,
\begin{equation*}\label{2.1b}
\aligned
H(x,y)=H_1(x,y)+\left\{
\aligned
&\frac1{2\pi}\nabla\log a(y)\cdot T(x-y),
\,\quad\,\,y\in\Omega,\\[1mm]
&\,\frac1{\pi}\,\nabla\log a(y)\cdot T(x-y),
\,\,\quad\,y\in\po,
\endaligned\right.
\endaligned
\end{equation*}
\indent {\upshape (ii)} the mapping
$y\in\oo\mapsto H_1(\cdot,y)$ belongs to
$ C^1\big(\Omega, C^{1}(\overline{\Omega})\big)
\cap C^1\big(\partial\Omega, C^{1}(\overline{\Omega})\big)$.\\
In this way, $y\in\oo\mapsto H(\cdot,y)\in C\big(\Omega, C^{\alpha}(\overline{\Omega})\big)
\cap C\big(\partial\Omega, C^{\alpha}(\overline{\Omega})\big)$ and
$H(x,y)\in C^\alpha\big(\oo\times\Omega\big)\cap C^\alpha\big(\oo\times\po\big)
\cap C^1\big(\oo\times\Omega\setminus\{x=y\}\big)
\cap C^1\big(\oo\times\po\setminus\{x=y\}\big)$
for any $\alpha\in(0,1)$,
and the corresponding Robin's function $y\in\oo\mapsto H(y,y)$ belongs to $C^1(\Omega)\cap C^1(\partial\Omega)$.
}

\vspace{1mm}
\vspace{1mm}
\vspace{1mm}
\vspace{1mm}

Let  $d>0$ be a sufficiently small but fixed number such that for any
$y\in\Omega$ with $\dist(y,\po)<d$, we can define a reflection
of $y$ across $\po$ along the outer normal direction, $y^*\in\Omega^c$,
and  get that  $|y-y^*|=2\dist(y,\partial\Omega)$.
Set
$$
\aligned
\Omega_d:=\left\{
y\in\Omega\big|\,\,
\dist(y,\partial\Omega)<d\,
\right\}.
\endaligned
$$

\vspace{1mm}
\vspace{1mm}

\noindent{\bf Lemma 2.2 (\cite{AP}).}\,\,{\it
There exists a mapping
$y\in\Omega_d\mapsto \mathrm{z}(\cdot,y)\in C\big(\Omega_d, C^{\alpha}(\overline{\Omega})\big)
\cap L^\infty\big(\Omega_d, C^{\alpha}(\overline{\Omega})\big)$
for any $\alpha\in(0,1)$ such that for any $x\in\oo$ and  $y\in\Omega_d$,
\begin{equation}\label{2.1c}
\aligned
H(x,y)=\frac1{2\pi}\log\frac{1}{|x-y^*|}+\mathrm{z}(x,y).
\endaligned
\end{equation}
Even more, for any $x\in\overline{\Omega}$ and  $y\in\Omega_d$,
\begin{equation}\label{2.1d}
\aligned
\mathrm{z}(x,y)=\frac1{2\pi}\nabla\log a(y)\cdot T(x-y)-\frac1{2\pi}\nabla\log a(y^*)\cdot T(x-y^*)+\tilde{\mathrm{z}}(x,y),
\endaligned
\end{equation}
where the mapping
$y\in\overline{\Omega_d}\mapsto \tilde{\mathrm{z}}(\cdot,y)$ belongs to
$ C^1\big(\overline{\Omega_d},\,C^{1}(\overline{\Omega})\big)$.
}

\vspace{1mm}
\vspace{1mm}
\vspace{1mm}
\vspace{1mm}

\noindent{\bf Corollary 2.3.}\,\,{\it
Under the assumptions in Lemma 2.2, the Robin's function
$y\in\Omega\mapsto H(y,y)$ satisfies
\begin{equation}\label{2.1e}
\aligned
H(y,y)=\frac1{2\pi}\log\frac{1}{|y-y^*|}+\mathrm{z}(y),
\,\,\quad\forall\,\,y\in\Omega_d,
\endaligned
\end{equation}
where $\mathrm{z}\in C^1\big(\overline{\Omega_d}\big)$ and
\begin{equation}\label{2.1f}
\aligned
\mathrm{z}(y)=\frac1{2\pi}\nabla\log a(y)\cdot T(0)-\frac1{2\pi}\nabla\log a(y^*)\cdot T(y-y^*)+\tilde{\mathrm{z}}(y,y),
\,\,\quad\forall\,\,y\in\Omega_d.
\endaligned
\end{equation}
}

\vspace{1mm}

The key ingredient to describe the shape of  the approximate solution of (\ref{1.1})
is based on the standard bubble
\begin{equation}\label{2.1}
\aligned
U_{\delta,\xi}(x)=\log\frac{8\delta^2}{(\delta^2+|x-\xi|^2)^2},\quad
\,\,\delta>0,\,\,\,\,\xi\in\mathbb{R}^2.
\endaligned
\end{equation}
It is well known from \cite{CL} that these are   all the solutions of the
problem
\begin{equation*}\label{2.2}
\left\{\aligned
&-\Delta
u=e^u\,\ \,\,\,\,
\textrm{in}\,\,\,\,
\mathbb{R}^2,\\
&\int_{\mathbb{R}^2}e^u<+\infty.
\endaligned\right.
\end{equation*}
The configuration space for $m$ concentration points $\xi=(\xi_1,\ldots,\xi_m)$ we
try to seek is the following
\begin{eqnarray}\label{2.3}
\mathcal{O}_p:=\left\{\,\xi=(\xi_1,\ldots,\xi_m)\in\Omega^l\times(\po)^{m-l}
\left|\,
\min_{i,j=1,\ldots,m,\,i\neq j}|\xi_i-\xi_j|>\frac{1}{p^\kappa},
\,\,\,\,\,\,\,
\min_{1\leq i\leq l}\dist(\xi_i,\po)>\frac{1}{p^\kappa}
\right.\right\},
\end{eqnarray}
where $l\in\{0,\ldots,m\}$ and $\kappa$ is given by
\begin{eqnarray}\label{2.3a}
\kappa=2(m^2+1).
\end{eqnarray}
Let  $m\in\mathbb{N}^*$ and
$\xi=(\xi_1,\ldots,\xi_m)\in\mathcal{O}_p$ be fixed.
Given number $\mu_i>0$,
$i=1,\ldots,m$, yet to be chosen, we define
\begin{equation}\label{2.4}
\aligned
\gamma=p^{\frac{p}{p-1}}e^{-\frac{p}{2(p-1)}},
\,\,\qquad\qquad\,\,\varepsilon=e^{-\frac{1}{4}p},
\,\,\qquad\qquad\ \,\,\delta_i=\varepsilon\mu_i,
\endaligned
\end{equation}
and
\begin{equation}\label{2.5}
\aligned
U_i(x)=\frac{1}{\gamma\mu_i^{2/(p-1)}}\left[U_{\delta_i,\xi_i}(x)
+\frac1p\omega_1\left(\frac{x-\xi_i}{\delta_i}\right)
+\frac1{p^2}\omega_2\left(\frac{x-\xi_i}{\delta_i}\right)
\right].
\endaligned
\end{equation}
Here, $\omega_j$, $j=1, 2$, are radial solutions of
\begin{equation}\label{2.6}
\aligned
\Delta\omega_j+\frac{8}{(1+|z|^2)^2}\omega_j=\frac{8}{(1+|z|^2)^2}f_j(|z|)
\quad\,\textrm{in}\,\,\,\,\mathbb{R}^2,
\endaligned
\end{equation}
with
\begin{equation}\label{2.7}
\aligned
f_1=\frac12U_{1,0}^2,\,\quad\quad\,
f_2=\omega_1U_{1,0}-\frac13U_{1,0}^3-\frac12\omega_1^2-\frac18U_{1,0}^4+\frac{1}{2}\omega_1U_{1,0}^2,
\endaligned
\end{equation}
having asymptotic
\begin{equation}\label{2.8}
\aligned
\omega_j(r)=\frac{C_j}{2}\log(1+r^2)+O\left(\frac1{1+r}\right),
\,\,\quad\,\,
\partial_{r}\omega_j(r)=\frac{C_jr}{1+r^2}+O\left(\frac{1}{1+r^2}\right)
\quad\,\,\,\textrm{as}\,\,\,r\rightarrow+\infty,\,\,\,r=|z|,
\endaligned
\end{equation}
where
\begin{equation*}\label{2.9}
\aligned
C_j=8\int_{0}^{\infty}t\frac{t^2-1}{(t^2+1)^3}f_j(t)dt,
\endaligned
\end{equation*}
in particular,
\begin{eqnarray}\label{2.10}
\omega_1(z)=\frac12U_{1,0}^2(z)+6\log(|z|^2+1)
+\frac{2\log8-10}{|z|^2+1}+\frac{|z|^2-1}{|z|^2+1}
\times\left\{-\frac12\log^28+2\log^2(|z|^2+1)\right.&&\nonumber\\
\left.
+4\int_{|z|^2}^{\infty}\frac{ds}{s+1}
\log\frac{s+1}{s}
-8\log|z|\log(|z|^2+1)\right\},
\,\,\qquad\qquad\qquad\qquad\qquad\qquad\qquad\,\,&&
\end{eqnarray}
and
\begin{equation}\label{2.11}
\aligned
C_1=12-4\log8
\endaligned
\end{equation}
(see \cite{CI,EMP}).
We now approximate the solution of problem (\ref{1.1}) by
\begin{equation}\label{2.12}
\aligned
U_\xi(x):=\sum_{i=1}^mPU_i(x)=\sum_{i=1}^m\big[U_i(x)+H_i(x)\big],
\endaligned
\end{equation}
where $H_i$  is a correction term defined as the solution of
\begin{equation}\label{2.13}
\aligned
\left\{\aligned
&-\Delta_a
H_i+H_i=\nabla\log a(x)\nabla U_i-U_i\,\,\,\,\,\,
\textrm {in}\,\,\,\,\,\,\Omega,\\[1mm]
&\frac{\partial H_i}{\partial \nu}=-\frac{\partial U_i}{\partial\nu}\,
\,\qquad\qquad\qquad\qquad\qquad\,\,\,
\textrm{on}\,\,\,\,\po.
\endaligned\right.
\endaligned
\end{equation}
In order to understand the asymptotic behavior of the correction term $H_i$, let us  first use the convention
\begin{equation}\label{2.14}
\aligned
c_i=\left\{
\aligned
&8\pi,
\qquad\textrm{if}\quad\,\xi_i\in\Omega,
\\[1mm]
&4\pi,
\qquad\textrm{if}\quad \xi_i\in\po.
\endaligned
\right.
\endaligned
\end{equation}
Furthermore, we have the following result, whose proof is postponed to the Appendix.

\vspace{1mm}
\vspace{1mm}
\vspace{1mm}
\vspace{1mm}

\noindent{\bf Lemma 2.4.}\,\,{\it
For any $0<\alpha<1$  and
$\xi=(\xi_1,\ldots,\xi_m)\in\mathcal{O}_p$,   then we have
\begin{eqnarray}\label{2.15}
H_i(x)=\frac{1}{\gamma\mu_i^{2/(p-1)}}\left[
\left(
1-\frac{C_1}{4p}-\frac{C_2}{4p^2}
\right)c_i H(x,\xi_i)
-\log(8\delta_i^2)+\left(\frac{C_1}{p}+\frac{C_2}{p^2}\right)\log\delta_i
+O\left(\delta_i^{\alpha/2}\right)
\right]
\end{eqnarray}
uniformly in $\oo$,
where
$H$ is the regular part of Green's function defined in
{\upshape(\ref{1.3})}.
}

\vspace{1mm}
\vspace{1mm}
\vspace{1mm}
\vspace{1mm}

From Lemma 2.4 we have that away from each point $\xi_i$,
namely  $|x-\xi_i|\geq 1/p^{2\kappa}$ for any $i=1,\ldots,m$,
\begin{equation}\label{2.20}
\aligned
U_{\xi}(x)=\sum\limits_{i=1}^{m}
\frac{1}{\gamma\mu_i^{2/(p-1)}}
\left[\left(1
-\frac{C_1}{4p}-\frac{C_2}{4p^2}\right)c_i G(x,\xi_i)
+\frac1pO\left(\delta_i^{\alpha/2}\right)\right].
\endaligned
\end{equation}
While for
$|x-\xi_i|< 1/p^{2\kappa}$  with some $i$,  by  (\ref{2.1}),
(\ref{2.5}), (\ref{2.8}), (\ref{2.15}) and the
fact that $H(\cdot,\xi_k)\in C^{\alpha}(\oo)$ for any $\xi_k\in\oo$ and any
$\alpha\in(0,1)$
we obtain
$$
\aligned
PU_i(x)=&\frac{1}{\gamma\mu_i^{2/(p-1)}}
\left[U_{1,0}\left(\frac{x-\xi_i}{\delta_i}\right)
+\frac1p\omega_1\left(\frac{x-\xi_i}{\delta_i}\right)
+\frac1{p^2}\omega_2\left(\frac{x-\xi_i}{\delta_i}\right)+\left(
1-\frac{C_1}{4p}-\frac{C_2}{4p^2}
\right)c_i H(\xi_i,\xi_i)\right.\\[1mm]
&\left.
-\log(8\delta_i^4)+\left(\frac{C_1}{p}+\frac{C_2}{p^2}\right)\log\delta_i
+O\left(
|x-\xi_i|^\alpha
+\delta_i^{\alpha/2}\right)
\right],
\endaligned
$$
and for any $k\neq i$,
$$
\aligned
PU_k(x)=&\frac{1}{\gamma\mu_k^{2/(p-1)}}
\left[
\log\frac{8\delta_k^2}{(\delta_k^2+|x-\xi_k|^2)^2}
+\frac1p\omega_1\left(\frac{x-\xi_k}{\delta_k}\right)
+\frac1{p^2}\omega_2\left(\frac{x-\xi_k}{\delta_k}\right)
\right.\\[1mm]
&\left.
+\left(
1-\frac{C_1}{4p}-\frac{C_2}{4p^2}
\right)c_k H(x,\xi_k)
-\log(8\delta_k^2)+\left(\frac{C_1}{p}+\frac{C_2}{p^2}\right)\log\delta_k
+O\left(\delta_k^{\alpha/2}\right)
\right]\\[1mm]
=&\frac{1}{\gamma\mu_k^{2/(p-1)}}
\left[
\left(
1-\frac{C_1}{4p}-\frac{C_2}{4p^2}
\right)c_k G(\xi_i,\xi_k)
+O\left(
|x-\xi_i|^\alpha
+\delta_k^{\alpha/2}\right)
\right].
\endaligned
$$
Hence for $|x-\xi_i|< 1/p^{2\kappa}$,
\begin{eqnarray}\label{2.21}
U_{\xi}(x)=\frac{1}{\gamma\mu_i^{2/(p-1)}}
\left[\,p+U_{1,0}\left(\frac{x-\xi_i}{\delta_i}\right)
+\frac1p\omega_1\left(\frac{x-\xi_i}{\delta_i}\right)
+\frac1{p^2}\omega_2\left(\frac{x-\xi_i}{\delta_i}\right)
+O\left(
|x-\xi_i|^\alpha
+\sum_{k=1}^m\delta_k^{\alpha/2}
\right)
\right]
\end{eqnarray}
is a good approximation for a solution of problem (\ref{1.1})
provided that the concentration
parameters $\mu_i$, $i=1,\ldots,m$, are the solution of the nonlinear system
\begin{eqnarray}\label{2.22}
\log\big(8\mu_i^4\big)=\left(
1-\frac{C_1}{4p}-\frac{C_2}{4p^2}
\right)c_i H(\xi_i,\xi_i)
+\left(\frac{C_1}{p}+\frac{C_2}{p^2}\right)\log\delta_i
+\left(
1-\frac{C_1}{4p}-\frac{C_2}{4p^2}
\right)\sum_{k\neq i}^m
\left(\frac{\mu_i}{\mu_k}\right)^{2/(p-1)}
c_k G(\xi_i,\xi_k).
\end{eqnarray}
Indeed, the parameters  $\mu=(\mu_1,\ldots,\mu_m)$
are well defined in system (\ref{2.22}), which is stated
as follows and proved in the Appendix.

\vspace{1mm}
\vspace{1mm}
\vspace{1mm}

\noindent{\bf Lemma 2.5.}\,\,{\it For
any points $\xi=(\xi_1,\ldots,\xi_m)\in\mathcal{O}_p$
and any  $p>1$ large enough, system {\upshape(\ref{2.22})}
has a unique solution
 $\mu=(\mu_1,\ldots,\mu_m)$  satisfying
\begin{equation}\label{2.23}
\aligned
1/C\leq\mu_i\leq Cp^\kappa
\qquad\forall\,\,\,i=1,\ldots,m,
\endaligned
\end{equation}
for some  $C>0$. Moreover,  for any $i=1,\ldots,m$, one has
\begin{equation}\label{2.60}
\aligned
\big|D_{\xi}\log\mu_i\big|\leq Cp^{\kappa},
\endaligned
\end{equation}
and
\begin{equation}\label{2.28}
\aligned
\mu_i=e\large^{-\frac{3}{4}+\frac14c_iH(\xi_i,\xi_i)
+\frac14\sum\large_{k=1,\,k\neq i}^m
c_k G(\xi_i,\xi_k)}\left[\,1+O\left(\frac{\log^2p}p\right)\right].
\endaligned
\end{equation}
}
\vspace{1mm}
\vspace{1mm}

\noindent{\bf Remark 2.6.}\,\,
Observe that
for $|x-\xi_i|=\delta_i|z|< 1/p^{2\kappa}$,
by (\ref{2.1}), (\ref{2.4}), (\ref{2.8}) and (\ref{2.23}),
$$
\aligned
p+U_{1,0}(z)
+\frac1p\omega_1(z)
+\frac1{p^2}\omega_2(z)
\geq p-2\log\left(1+|z|^2\right)+O\left(1\right)
\geq8\kappa\log p+
4\log\mu_i
+O\left(1\right)>7\kappa\log p.
\endaligned
$$
Hence by (\ref{2.21}), we can easily get that
$0<U_\xi\leq2\sqrt{e}$ in $B_{1/p^{2\kappa}}(\xi_i)$,
and
$\sup_{B_{1/p^{2\kappa}}(\xi_i)}U_\xi\rightarrow\sqrt{e}$
as $p\rightarrow+\infty$.
Moreover, by the maximum principle, we see that $G(x,\xi_i)>0$
over $\overline{\Omega}$ and thus by (\ref{2.20}),
$U_\xi$ is a  positive, uniformly bounded function over $\overline{\Omega}$.
In conclusion, $0<U_\xi\leq2\sqrt{e}$ over $\overline{\Omega}$.

\vspace{1mm}
\vspace{1mm}
\vspace{1mm}
\vspace{1mm}

Let us  perform  the change of variables
\begin{equation*}\label{2.29}
\aligned
\upsilon(y)=\varepsilon^{2/(p-1)}u(\varepsilon y),\,\,\,
\,\quad\forall\,\,\,y\in\Omega_{p}:=e^{p/4}\Omega.
\endaligned
\end{equation*}
Then by the definition of $\varepsilon$ in (\ref{2.4}),
$u(x)$ solves equation (\ref{1.1}) if and only if the function
$\upsilon(y)$ satisfies
\begin{equation}\label{2.30}
\begin{array}{ll}
\left\{\aligned
&-\Delta_{a(\varepsilon y)}\upsilon+\varepsilon^2\upsilon=\upsilon^p,
\,\,\,\,\,\,
\upsilon>0\,\,\,\,\,\,
\textrm{in}\,\,\,\,\,\Omega_p,\\[1mm]
&\frac{\partial \upsilon}{\partial\nu}=0
\qquad\qquad\qquad\qquad\qquad\quad
\,\textrm{on}\,\,\,\partial\Omega_p.
\endaligned\right.
\end{array}
\end{equation}
We
write $\xi_i'=\xi_i/\varepsilon$, $i=1,\ldots,m$
and define the initial approximate solution of (\ref{2.30}) as
\begin{equation}\label{2.31}
\aligned
V_{\xi'}(y)=\varepsilon^{2/(p-1)}U_\xi(\varepsilon y),
\endaligned
\end{equation}
with
$\xi'=(\xi_1',\ldots,\xi_m')$
and
$U_\xi$  defined in (\ref{2.12}).
Let us set
$$
\aligned
S_p(\upsilon)=-\Delta_{a(\varepsilon y)}\upsilon+\varepsilon^2\upsilon-\upsilon^p_{+},
\,\qquad\,\textrm{where}\,\,\,\,\upsilon_{+}=\max\{\upsilon,0\},
\endaligned
$$
and introduce the  functional
\begin{equation*}\label{2.32}
\aligned
I_p(\upsilon)=\frac12\int_{\Omega_p}a(\varepsilon y)\left(
|\nabla \upsilon|^2+\varepsilon^2\upsilon^2
\right)dy-\frac1{p+1}\int_{\Omega_p}a(\varepsilon y)\upsilon^{p+1}_{+}dy,
\,\quad\,\upsilon\in H^1(\Omega_p),
\endaligned
\end{equation*}
whose nontrivial critical points are solutions of problem (\ref{2.30}).
In fact, by the maximum principle, problem (\ref{2.30}) is equivalent to
$$
\aligned
S_p(\upsilon)=0,\,\quad\,\upsilon_{+}\not\equiv0
\,\quad\,\textrm{in}\,\,\,\,\Omega_p,\,\qquad\qquad\,
\frac{\partial\upsilon}{\partial\nu}=0
\,\quad\,\textrm{on}\,\,\,\partial\Omega_p.
\endaligned
$$
We will look for solutions of problem (\ref{2.30}) in the form
$\upsilon=V_{\xi'}+\phi$, where $\phi$ will represent a higher-order correction
in the expansion of $\upsilon$. Observe that
\begin{equation*}\label{2.33}
\aligned
S_p(V_{\xi'}+\phi)=\mathcal{L}(\phi)+R_{\xi'}+N(\phi)=0,
\endaligned
\end{equation*}
where
\begin{equation*}\label{2.34}
\aligned
\mathcal{L}(\phi)=-\Delta_{a(\varepsilon y)}\phi+\varepsilon^2\phi-W_{\xi'}\phi\,\quad\,\textrm{with}\,\quad\,
W_{\xi'}=pV_{\xi'}^{p-1},
\endaligned
\end{equation*}
and
\begin{equation}\label{2.35}
\aligned
R_{\xi'}=-\Delta_{a(\varepsilon y)}V_{\xi'}+\varepsilon^2V_{\xi'}-V_{\xi'}^p,
\quad\quad\quad
N(\phi)=-\big[(V_{\xi'}+\phi)_+^p-V_{\xi'}^p-pV_{\xi'}^{p-1}\phi\big].
\endaligned
\end{equation}
In terms of $\phi$,
problem (\ref{2.30}) becomes
\begin{equation}\label{2.36}
\aligned
\left\{\aligned
&\mathcal{L}(\phi)=-\big[
R_{\xi'}+N(\phi)
\big]
\quad\textrm{in}\,\,\,\,\,\,\Omega_p,\\
&\frac{\partial \phi}{\partial\nu}=0
\quad\qquad\qquad\qquad\,\,\,\,
\textrm{on}\,\,\,\,
\partial\Omega_p.
\endaligned\right.\endaligned
\end{equation}

For any  $\xi=(\xi_1,\ldots,\xi_m)\in\mathcal{O}_{p}$ and $h\in L^\infty(\Omega_p)$, let us  introduce a weighted $L^\infty$-norm
defined as
\begin{equation}\label{2.37}
\aligned
\|h\|_{*}=
\sup_{y\in\Omega_p}\left|\left(
\sum\limits_{i=1}^m\frac{\mu_i^\sigma}{(\mu_i+|y-\xi'_i|)^{2+\sigma}}
+\varepsilon^2
\right)^{-1}h(y)
\right|,
\endaligned
\end{equation}
where $\sigma>0$ is  small but fixed  independent of $p$.
With respect to the $\|\cdot\|_{*}$-norm,
the error term $R_{\xi'}$ defined in (\ref{2.35}) can be estimated as follows.

\vspace{1mm}
\vspace{1mm}
\vspace{1mm}
\vspace{1mm}

\noindent{\bf Proposition 2.7.}\,\,{\it
Let $m$ be a positive integer.
There exist constants $C>0$ and $p_m>1$ such that
for any $\xi=(\xi_1,\ldots,\xi_m)\in\mathcal{O}_{p}$ and any $p>p_m$,
\begin{equation}\label{2.38}
\aligned
\|R_{\xi'}\|_{*}\leq\frac{C}{p^4}.
\endaligned
\end{equation}
}

\begin{proof}
Observe that, by (\ref{2.12}), (\ref{2.13}) and (\ref{2.31}),
$$
\aligned
-\Delta_{a(\varepsilon y)}V_{\xi'}+\varepsilon^2V_{\xi'}
=\varepsilon^{2}\sum_{i=1}^m\varepsilon^{2/(p-1)}\left[-\Delta_{a}\big(U_i+H_i\big)+\big(U_i+H_i\big)\right]
=-\varepsilon^{2p/(p-1)}\sum_{i=1}^m\Delta U_i.
\endaligned
$$
Then by (\ref{2.1}), (\ref{2.4}), (\ref{2.5}) and (\ref{2.6}),
\begin{equation}\label{2.39}
\aligned
-\Delta_{a(\varepsilon y)}V_{\xi'}+\varepsilon^2V_{\xi'}
=\sum_{i=1}^m\frac{1}{p^{p/(p-1)}\mu_i^{2p/(p-1)}}
e^{U_{1,0}(z)}\left[1
-\frac{1}{p}f_1(z)
-\frac{1}{p^2}f_2(z)+\frac{1}{p}\omega_1(z)
+\frac{1}{p^2}\omega_2(z)
\right]
\endaligned
\end{equation}
with  $z=(\varepsilon y-\xi_i)/\delta_i$.
By (\ref{2.1}),  (\ref{2.8}) and (\ref{2.23})
we get, if $|\varepsilon y-\xi_i|=\delta_i|z|\geq1/p^{2\kappa}$ for any $i=1,\ldots,m$,
then
$$
\aligned
U_{1,0}(z)=-p+O\left(\log p\right),\,\qquad\quad\,\omega_{j}(z)=\frac{1}{4}C_{j}p+O\left(\log p\right),
\,\,\quad\,j=1,\,2,
\endaligned
$$
and hence, by  (\ref{2.7})  and (\ref{2.39}),
\begin{equation}\label{2.40}
\aligned
-\Delta_{a(\varepsilon y)}V_{\xi'}+\varepsilon^2V_{\xi'}
=\sum_{i=1}^m\frac{e^{U_{1,0}(z)}}{p^{p/(p-1)}\mu_i^{2p/(p-1)}}
O\left(p^2
\right).
\endaligned
\end{equation}
On the other hand, in the same region, by (\ref{2.20}) and (\ref{2.31}) we get
\begin{equation*}\label{2.41}
\aligned
V_{\xi'}^p
=O\left(\frac1p\left(\frac{\log p}{p}\right)^p
\right),
\endaligned
\end{equation*}
which, together with (\ref{2.4}), (\ref{2.23}) and (\ref{2.40}), implies
\begin{eqnarray}\label{2.42}
\left|\left(
\sum\limits_{i=1}^m\frac{\mu_i^\sigma}{(\mu_i+|y-\xi'_i|)^{2+\sigma}}
+\varepsilon^2
\right)^{-1}
R_{\xi'}(y)
\right|
\leq C
\left[p\sum_{i=1}^m
\left|\frac{y-\xi'_i}{\mu_i}\right|^{\sigma-2}+\frac1p\left(\frac{\sqrt{e}\log p}{p}\right)^p\right]
=o\left(e^{-p/4}
\right).
\end{eqnarray}
Let us  fix an index $i\in\{1,\ldots,m\}$
and the region $|\varepsilon y-\xi_i|=\delta_i|z|<1/p^{2\kappa}$.
By (\ref{2.21}), (\ref{2.31}) and the relation
\begin{equation}\label{2.43}
\aligned
\left(
\frac{p\varepsilon^{2/(p-1)}}{\gamma\mu_i^{2/(p-1)}}
\right)^p=\frac{1}{p^{p/(p-1)}\mu_i^{2p/(p-1)}},
\endaligned
\end{equation}
we get, for  $|\varepsilon y-\xi_i|=\delta_i|z|< 1/p^{2\kappa}$,
\begin{equation}\label{2.44}
\aligned
V^p_{\xi'}(y)=\frac{1}{p^{p/(p-1)}\mu_i^{2p/(p-1)}}
\left\{
1+\frac{U_{1,0}(z)}p+\frac{\omega_1(z)}{p^2}
+\frac1{p^3}\left[\omega_2(z)+O\left(
p^2\delta_i^\alpha|z|^\alpha
+p^2\sum_{k=1}^m\delta_k^{\alpha/2}\right)
\right]\right\}^p.
\endaligned
\end{equation}
From a Taylor expansion of the exponential and logarithmic functions
\begin{eqnarray}\label{2.45}
\left(1+\frac{a}p+\frac{b}{p^2}+\frac{c}{p^3}\right)^p
=e^a\left[1+\frac1p\left(b-\frac{a^2}2\right)+\frac1{p^2}\left(c-a b+\frac{a^3}3
+\frac{b^2}2-\frac{a^2b}2+\frac{a^4}8\right)
\right.
&&\nonumber\\
\left.+O\left(\frac{\log^6(|z|+2)}{p^3}\right)\right],
\,\quad\qquad\qquad\qquad\qquad\qquad\qquad\qquad\,\,\,
&&
\end{eqnarray}
which holds for $|z|\leq Ce^{p/8}$
provided
$-4\log(|z|+2)\leq a(z)\leq C$ and $|b(z)|+|c(z)|\leq C\log(|z|+2)$, so we have
that for $|\varepsilon y-\xi_i|=\delta_i|z|\leq\sqrt{\delta_i}/p^{2\kappa}$,
\begin{eqnarray*}\label{2.46}
V^p_{\xi'}(y)=\frac{1}{p^{p/(p-1)}\mu_i^{2p/(p-1)}}e^{U_{1,0}(z)}
\left[\,
1+\frac1p\left(\omega_1
-\frac12U_{1,0}^2
\right)(z)+\frac1{p^2}\left(\omega_2
-\omega_1U_{1,0}+\frac13U_{1,0}^3+\frac12\omega_1^2
\right.\right.&&\nonumber\\
\left.\left.
-\frac12\omega_1U_{1,0}^2
+\frac18U_{1,0}^4
\right)(z)
+O\left(\frac{\log^6(|z|+2)}{p^3}
+\delta_i^\alpha|z|^\alpha
+\sum_{k=1}^m\delta_k^{\alpha/2}\right)
\right],
\,\qquad\qquad\qquad\qquad
&&
\end{eqnarray*}
which combined with (\ref{2.7}) and (\ref{2.39})   gives
$$
\aligned
R_{\xi'}(y)=-\Delta_{a(\varepsilon y)}V_{\xi'}+\varepsilon^2V_{\xi'}-V_{\xi'}^p=
\frac{1}{p^{p/(p-1)}\mu_i^{2p/(p-1)}}e^{U_{1,0}(z)}
O\left(\frac{\log^6(|z|+2)}{p^3}
+\delta_i^\alpha|z|^\alpha
+\sum_{k=1}^m\delta_k^{\alpha/4}\right).
\endaligned
$$
Hence, in this region we get
\begin{eqnarray}\label{2.47}
\left|\left(
\sum\limits_{i=1}^m\frac{\mu_i^\sigma}{(\mu_i+|y-\xi'_i|)^{2+\sigma}}
+\varepsilon^2
\right)^{-1}
R_{\xi'}(y)
\right|
\leq\frac{C}{p^4\big(\big|\frac{y-\xi'_i}{\mu_i}\big|+1\big)^{2-\sigma}}
\log^6\left(\left|\frac{y-\xi'_i}{\mu_i}\right|+2\right)
=O\left(\frac1{p^4}
\right).
\end{eqnarray}
Finally, in the remaining region
$\sqrt{\delta_i}/p^{2\kappa}<|\varepsilon y-\xi_i|=\delta_i|z|<1/p^{2\kappa}$, we have that,
by (\ref{2.7}) and (\ref{2.39}),
\begin{equation*}\label{2.48}
\aligned
-\Delta_{a(\varepsilon y)}V_{\xi'}+\varepsilon^2V_{\xi'}=
\frac{e^{U_{1,0}(z)}}{p^{p/(p-1)}\mu_i^{2p/(p-1)}}O\left(
p^2
\right),
\endaligned
\end{equation*}
and by (\ref{2.44}),
\begin{equation*}\label{2.49}
\aligned
V^p_{\xi'}(y)=
\frac{e^{U_{1,0}(z)}}{p^{p/(p-1)}\mu_i^{2p/(p-1)}}O\left(
1
\right),
\endaligned
\end{equation*}
since $(1+\frac{s}p)^p\leq e^s$.
Thus, in this region,
\begin{eqnarray}\label{2.50}
\left|\left(
\sum\limits_{i=1}^m\frac{\mu_i^\sigma}{(\mu_i+|y-\xi'_i|)^{2+\sigma}}
+\varepsilon^2
\right)^{-1}
R_{\xi'}(y)
\right|
\leq Cp
\left|\frac{y-\xi'_i}{\mu_i}\right|^{\sigma-2}
=o\left(e^{-p/8}
\right).
\end{eqnarray}
Combining (\ref{2.37}), (\ref{2.42}), (\ref{2.47}) with  (\ref{2.50}),
we conclude that estimate (\ref{2.38}) holds.
\end{proof}

\vspace{1mm}

\section{Analysis of the linearized operator}
In this section, we prove bounded invertibility of the operator $\mathcal{L}$,
uniformly on $\xi\in\mathcal{O}_p$, by using the weighted
$L^{\infty}$-norm defined in (\ref{2.37}). Let us recall that
$\mathcal{L}(\phi)=-\Delta_{a(\varepsilon y)}\phi+\varepsilon^2\phi-W_{\xi'}\phi$,
where $W_{\xi'}=pV_{\xi'}^{p-1}$.
As in Proposition 2.7, we have  the following expansions with respect to the potential $W_{\xi'}$.

\vspace{1mm}
\vspace{1mm}
\vspace{1mm}
\vspace{1mm}

\noindent{\bf Lemma 3.1.}\,\,{\it
Let $m$ be a positive integer. There exist constants
$D_0>0$  and $p_m>1$  such that
\begin{equation}\label{2.51}
\aligned
W_{\xi'}(y)\leq D_0\sum_{i=1}^{m}\frac{1}{\mu_i^2}
e^{U_{1,0}\big(\frac{y-\xi'_i}{\mu_i}\big)}
\endaligned
\end{equation}
for any points $\xi=(\xi_1,\ldots,\xi_m)\in\mathcal{O}_{p}$ and any $p\geq p_m$.
Furthermore,
\begin{equation}\label{2.52}
\aligned
W_{\xi'}(y)=\frac{8}{\mu_i^2(1+|z|^2)^2}\left[1+
\frac1p\left(\omega_1-U_{1,0}-\frac12U_{1,0}^2\right)(z)
+O\left(\frac{\log^4\big(|z|+2\big)}{p^2}\right)\right]
\endaligned
\end{equation}
for any $|\varepsilon y-\xi_i|\leq\sqrt{\delta_i}/p^{2\kappa}$, where $z=(y-\xi'_i)/\mu_i$.
}

\vspace{1mm}
\vspace{1mm}
\vspace{1mm}
\vspace{1mm}

\noindent{\it Proof.}\,\,\,If  $|\varepsilon y-\xi_i|=\delta_i|z|< 1/p^{2\kappa}$ for some $i=1,\ldots,m$,
by (\ref{2.21}), (\ref{2.31}) and (\ref{2.43}),
$$
\aligned
W_{\xi'}(y)=&p\left(
\frac{\varepsilon^{2/(p-1)}}{\gamma\mu_i^{2/(p-1)}}
\right)^{p-1}
\left[\,p+U_{1,0}(z)
+\frac1p\omega_1(z)
+\frac1{p^2}\omega_2(z)
+O\left(
\delta_i^\alpha|z|^\alpha
+\sum_{k=1}^m\delta_k^{\alpha/2}\right)
\right]^{p-1}\\[1mm]
=&\frac{1}{\mu_i^2}
\left[
1+\frac{1}pU_{1,0}(z)+\frac{1}{p^2}\omega_1(z)
+\frac{1}{p^3}\omega_2(z)+\frac{1}{p}O\left(
\delta_i^\alpha|z|^\alpha
+\sum_{k=1}^m\delta_k^{\alpha/2}\right)
\right]^{p-1}.
\endaligned
$$
In this region, using the fact that $(1+a/p)^{p-1}\leq e^{(p-1)a/p}$
and $U_{1,0}(z)\geq-p+O\big(\log p\big)$,
we get
$$
\aligned
W_{\xi'}(y)\leq \frac{C}{\mu_i^2}
e^{U_{1,0}(z)}
e^{-U_{1,0}(z)/p}
=O\left(
\frac{1}{\mu_i^2}
e^{U_{1,0}(z)}
\right).
\endaligned
$$
In particular, from a slight modification of formula (\ref{2.45}), namely
$$
\aligned
\left(1+\frac{a}p+\frac{b}{p^2}+\frac{c}{p^3}\right)^{p-1}=e^a\left[1+\frac1p\left(b-a-\frac{a^2}2\right)
+O\left(\frac{\log^4(|z|+2)}{p^2}\right)\right],
\endaligned
$$
we conclude that if $|\varepsilon y-\xi_i|=\delta_i|z|\leq\sqrt{\delta_i}/p^{2\kappa}$, then
$$
\aligned
W_{\xi'}(y)
=\frac{1}{\mu_i^2}e^{U_{1,0}(z)}\left[1+\frac1p\left(\omega_1-U_{1,0}-\frac12U_{1,0}^2
\right)(z)
+O\left(\frac{\log^4(|z|+2)}{p^2}\right)\right].
\endaligned
$$
Additionally, if $|\varepsilon y-\xi_i|=\delta_i|z|\geq1/p^{2\kappa}$  for
all $i$, then by (\ref{2.4}), (\ref{2.20}), (\ref{2.23}) and (\ref{2.31}),
$$
\aligned
\,\,\quad\qquad\qquad\,\,
W_{\xi'}(y)=p\varepsilon^2U_\xi^{p-1}(\varepsilon y)
\leq p\varepsilon^2\left(
\frac{C\log p}{\gamma}
\right)^{p-1}=O\left(
\left(\frac{\log p}{p}\right)^{p-1}
\right)=o\left(\frac{1}{\mu_i^2}
e^{U_{1,0}(z)}
\right).
\qquad\qquad\quad\,\,\,\,
\square
\endaligned
$$

\vspace{1mm}
\vspace{1mm}
\vspace{1mm}

\noindent{\bf Remark 3.2.}\,\,\,As for $W_{\xi'}$, we mention that
if $|\varepsilon y-\xi_i|< 1/p^{2\kappa}$ for some $i=1,\ldots,m$, then
$$
\aligned
p\left[V_{\xi'}(y)+O\left(\frac1{p^3}\right)\right]^{p-2}\leq
Cp\left(\frac{p\varepsilon^{2/(p-1)}}{\gamma\mu_i^{2/(p-1)}}
\right)^{p-2}e^{\frac{p-2}{p}U_{1,0}\big(\frac{y-\xi'_i}{\mu_i}\big)}
=O\left(\frac{1}{\mu_i^2}e^{U_{1,0}\big(\frac{y-\xi'_i}{\mu_i}\big)}\right).
\endaligned
$$
Since this estimate is true if  $|\varepsilon y-\xi_i|\geq1/p^{2\kappa}$
for all $i$, we get
\begin{equation*}\label{2.53}
\aligned
p\left[V_{\xi'}(y)+O\left(\frac1{p^3}\right)\right]^{p-2}\leq
 C\sum\limits_{i=1}^{m}\frac{1}{\mu_i^2}e^{U_{1,0}\big(\frac{y-\xi'_i}{\mu_i}\big)}.
\endaligned
\end{equation*}

\vspace{1mm}
\vspace{1mm}
\vspace{1mm}

Let
\begin{equation}\label{3.2}
\aligned
Z_{0}(z)=\frac{|z|^2-1}{|z|^2+1},
\,\,\quad\qquad\quad\,\,
Z_{j}(z)=\frac{z_j}{|z|^2+1},\,\,\,\,j=1,\,2.
\endaligned
\end{equation}
It is well known (see \cite{BP,CL}) that
\begin{itemize}
  \item any bounded
solution to
\begin{equation}\label{3.3}
\aligned
\Delta
\phi+\frac{8}{(1+|z|^2)^2}\phi=0
\,\quad\textrm{in}\,\,\,\mathbb{R}^2,
\endaligned
\end{equation}
is a linear combination of $Z_j$, $j=0,1,2$;
  \item any bounded
solution to
\begin{equation}\label{3.4}
\aligned
\Delta
\phi+\frac{8}{(1+|z|^2)^2}\phi=0
\,\quad\textrm{in}\,\,\,\mathbb{R}^2_{+},
\,\qquad\qquad\,\frac{\partial\phi}{\partial\nu}=0
\,\quad\textrm{on}\,\,\,\partial\mathbb{R}^2_{+},
\endaligned
\end{equation}
where $\mathbb{R}^2_{+}=\{(z_1,z_2):\,z_2>0\}$,
is a linear combination of $Z_j$, $j=0,1$.
\end{itemize}
Now we consider
the following linear problem: given $h\in C(\overline{\Omega}_p)$ and points
$\xi=(\xi_1,\ldots,\xi_m)\in\mathcal{O}_{p}$, we find a function $\phi\in H^2(\Omega_p)$ and scalars
$c_{ij}\in\mathbb{R}$, $i=1,\ldots,m$, $j=1,J_i$,
such that
\begin{equation}\label{3.1}
\left\{\aligned
&\mathcal{L}(\phi)=-\Delta_{a(\varepsilon y)}\phi+\varepsilon^2\phi-W_{\xi'}\phi=h
+\frac1{a(\varepsilon y)}\sum\limits_{i=1}^m\sum\limits_{j=1}^{J_i}c_{ij}\chi_i\,Z_{ij}\,\,\ \,
\,\textrm{in}\,\,\,\,\,\,\Omega_p,\\
&\frac{\partial\phi}{\partial\nu}=0\,\,\,\,\,\,\,\,
\ \ \ \ \ \ \ \ \ \ \ \ \ \,\,
\qquad\qquad\qquad\quad\qquad\qquad\qquad\qquad
\ \,\ \ \ \,\,\,\ \,
\ \,\textrm{on}\,\,\,\,\partial\Omega_{p},\\[1mm]
&\int_{\Omega_p}\chi_i\,Z_{ij}\phi=0
\,\qquad\qquad\qquad\quad
\qquad\qquad\qquad\forall\,\,i=1,\ldots,m,\,\,\,j=1, J_i,
\endaligned\right.
\end{equation}
where
$J_i=2$ if $i=1,\ldots,l$ while $J_i=1$ if $i=l+1,\ldots,m$,
and $Z_{ij}$, $\chi_i$, are defined as follows.

Let $\chi:\mathbb{R}\rightarrow[0,1]$ be
a smooth, non-increasing cut-off function  such
that for a large but fixed number $R_0>0$,  $\chi(r)=1$ if $r\leq
R_0$,  and $\chi(r)=0$ if $r\geq R_0+1$.

For $i=1,\ldots,l$ (corresponding to interior spike case), we define
\begin{equation}\label{3.5}
\aligned
\chi_i(y)=\chi\left(
\frac{|y-\xi'_i|}{\mu_i}
\right),
\,\qquad\,
Z_{ij}(y)=\frac{1}{\mu_i}Z_j\left(\frac{y-\xi_i'}{\mu_i}\right),
\,\quad\,j=0,1,2.
\endaligned
\end{equation}

For $i=l+1,\ldots,m$ (corresponding to boundary spike case), we have to straighten
the boundary first. More precisely, at the boundary point $\xi_i\in\po$, we
define a rotation map  $A_i: \mathbb{R}^2\mapsto\mathbb{R}^2$ such that
$A_i\nu_{\Omega}(\xi_i)=\nu_{\mathbb{R}_+^2}(0)$.
Let $\mathcal{G}(x_1)$ be the defining function
for the boundary $A_i(\po-\{\xi_i\})$ in a small
neighborhood $B_\delta(0,0)$ of the origin,
that is, there exist $R_1>0$, $\delta>0$ small and
a smooth function $\mathcal{G}:(-R_1,R_1)\mapsto\mathbb{R}$
satisfying $\mathcal{G}(0)=0$, $\mathcal{G}'(0)=0$ and such that
$A_i(\Omega-\{\xi_i\})\cap B_\delta(0,0)=\{(x_1,x_2):\,-R_1<x_1<R_1,\,x_2>\mathcal{G}(x_1)\}\cap B_\delta(0,0)$.
Furthermore, we consider the flattening change of variables
$F_i: B_\delta(0,0)\cap\overline{A_i(\Omega-\{\xi_i\})}\mapsto\mathbb{R}^2$
be defined by $F_i=(F_{i1}, F_{i2})$, where
\begin{equation}\label{3.6}
\aligned
F_{i1}=x_1+\frac{x_2-\mathcal{G}(x_1)}{\,1+|\mathcal{G}'(x_1)|^2\,}\mathcal{G}'(x_1)\ \qquad\textrm{and}\ \qquad\,F_{i2}=x_2-\mathcal{G}(x_1).
\endaligned
\end{equation}
Then for any $i=l+1,\ldots,m$, we set
\begin{equation}\label{3.7}
\aligned
F_i^p(y)=e^{p/4}F_i\big(A_i(e^{-p/4}y-\xi_i)\big)=\frac{1}{\varepsilon}F_i\big(A_i(\varepsilon y-\xi_i)\big),
\endaligned
\end{equation}
and define
\begin{equation}\label{3.8}
\aligned
\chi_i(y)=\chi\left(
\frac{1}{\mu_i}|F_i^p(y)|
\right),
\,\qquad\,
Z_{ij}(y)=\frac{1}{\mu_i}Z_{j}\left(\frac{1}{\mu_i}F_i^p(y)\right),
\,\quad\,j=0,1.
\endaligned
\end{equation}
It is important to note that $F^p_i$, $i=l+1,\ldots,m$,
preserves the  Neumann boundary condition and
\begin{equation}\label{3.9}
\aligned
\Delta_{a(\varepsilon y)}Z_{i0}+\frac{8\mu_i^2}{(\mu_i^2+|y-\xi'_i|^2)^2}Z_{i0}
=O\left(\frac{\varepsilon\mu_i}{
(\mu_i+|y-\xi'_i|)^3}\right),
\,\,\,
\,\ \,\forall\,\,i=1,\ldots,m.
\endaligned
\end{equation}

\vspace{1mm}
\vspace{1mm}
\vspace{1mm}
\vspace{1mm}

\noindent {\bf Proposition 3.3.}\,\,{\it
Let $m$ be a positive integer.
Then there exist constants $C>0$  and $p_m>1$ such
that for any  $p>p_m$,    any points
$\xi=(\xi_1,\ldots,\xi_m)\in\mathcal{O}_{p}$ and any $h\in C(\overline{\Omega}_p)$,
there is a unique solution $\phi\in
H^2(\Omega_p)$ of problem {\upshape(\ref{3.1})}
for some coefficients  $c_{ij}\in\mathbb{R}$,  $i=1,\ldots,m$, $j=1,J_i$,  which satisfies
\begin{equation}\label{3.10}
\aligned
\|\phi\|_{L^{\infty}(\Omega_p)}\leq Cp\|h\|_{*}.
\endaligned
\end{equation}

}\indent The proof of this result will be split into  four steps which we state and prove next.

{\bf Step 1:} Constructing a suitable barrier.

\vspace{1mm}
\vspace{1mm}
\vspace{1mm}

\noindent{\bf Lemma 3.4.}\,\,{\it There exist constants $R_1>0$ and $C>0$,
independent of $p$, such that
for any sufficiently large
$p$,  any points $\xi=(\xi_1,\ldots,\xi_m)\in\mathcal{O}_p$ and any $\sigma\in(0,1)$, there is
a function
$$
\aligned
\psi:\,\,\Omega_p\setminus\bigcup_{i=1}^mB_{R_1\mu_i}(\xi'_i)\,\,\mapsto\mathbb{R}
\endaligned
$$
smooth  and
positive so that
$$
\aligned
\mathcal{L}(\psi)=-\Delta_{a(\varepsilon y)}\psi+\varepsilon^2\psi-W_{\xi'}\psi&\geq\sum_{i=1}^m\frac{\mu_i^\sigma}{|y-\xi'_i|^{2+\sigma}}+\varepsilon^2
\,\,\,\quad\,\,
\textrm{in}\,\,\,\,\,\,\,\Omega_p\setminus\bigcup_{i=1}^mB_{R_1\mu_i}(\xi'_i),\\
\frac{\partial\psi}{\partial\nu}&\geq0
\,\,\,\qquad\qquad\qquad\qquad\quad\,\,
\textrm{on}\,\,\,\,\partial\Omega_p\setminus\bigcup_{i=1}^mB_{R_1\mu_i}(\xi'_i),\\
\psi&>0
\,\,\,\qquad\qquad\qquad\qquad\quad\,\,
\textrm{in}\,\,\,\,\,\,\,
\Omega_p\setminus\bigcup_{i=1}^mB_{R_1\mu_i}(\xi'_i),\\
\psi&\geq1
\,\,\,\qquad\qquad\qquad\qquad\quad\,\,
\textrm{on}\,\,\,\,
\Omega_p\cap\left(\bigcup_{i=1}^m\partial B_{R_1\mu_i}(\xi'_i)\right).
\endaligned
$$
Moreover, $\psi$ is  uniformly bounded, i.e.
$$
\aligned
1<\psi\leq C\,\,\,\quad\textrm{in}\,\,\,
\,\Omega_p\setminus\bigcup_{i=1}^mB_{R_1\mu_i}(\xi'_i).
\endaligned
$$
}

\begin{proof}
Let us take
$$
\aligned
\psi=\sum_{i=1}^m
\left(1-
\frac{\mu_i^\sigma}{|y-\xi'_i|^\sigma}
\right)+
C_1
\Psi_0(y),
\endaligned
$$
where $\Psi_0$ is the unique solution of
$$
\aligned
\left\{
\aligned
&-\Delta_{a(\varepsilon y)}\Psi_0+\varepsilon^2\Psi_0=\varepsilon^2
\,\quad\,\,\textrm{in}\,\,\,\,\Omega_p,\\
&\frac{\partial\Psi_0}{\partial\nu}=\varepsilon
\,\qquad\qquad\qquad\qquad\,
\textrm{on}\,\,\,\,\partial\Omega_p.
\endaligned
\right.
\endaligned
$$
Observing that $\Psi_0$ is uniformly bounded in $\Omega_p$, it is directly checked that,
choosing the positive constant $C_1$ larger if necessary, $\psi$ satisfies all the properties
of the lemma for large enough numbers $R_1$ and $p$.
\end{proof}

\vspace{1mm}
\vspace{1mm}
\vspace{1mm}

{\bf Step 2:} An auxiliary linear equation. Given $h\in C^{0,\alpha}(\overline{\Omega}_p)$
and $\xi=(\xi_1,\ldots,\xi_m)\in\mathcal{O}_p$, we first study the linear equation
\begin{equation}\label{3.11}
\left\{\aligned
&\mathcal{L}(\phi)=-\Delta_{a(\varepsilon y)}\phi+\varepsilon^2\phi-W_{\xi'}\phi=h
\,\,\ \,
\,\textrm{in}\,\,\,\,\,\,\Omega_p,\\
&\frac{\partial\phi}{\partial\nu}=0\,\,\,\,\,
\qquad\qquad\qquad\qquad\qquad\qquad
\ \,\textrm{on}\,\,\,\,\partial\Omega_{p}.
\endaligned\right.
\end{equation}
For the solution of (\ref{3.11}) satisfying the orthogonality conditions with respect to
$Z_{ij}$, $i=1,\ldots,m$, $j=0,1,J_i$, we prove the following a priori estimate.

\vspace{1mm}
\vspace{1mm}
\vspace{1mm}

\noindent{\bf Lemma 3.5.}\,\,{\it There exist $R_0>0$ and  $p_m>1$ such that for any $p>p_m$
and any solution $\phi$ of {\upshape (\ref{3.11})} with the orthogonality conditions
\begin{equation}\label{3.12}
\aligned
\int_{\Omega_p}\chi_iZ_{ij}\phi=0\,\,\,\,\,
\,\,\,\,\forall\,\,i=1,\ldots,m,\,\,j=0,1,J_i,
\endaligned
\end{equation}
we have
\begin{equation*}\label{3.13}
\aligned
\|\phi\|_{L^{\infty}(\Omega_p)}\leq C
\|h\|_{*},
\endaligned
\end{equation*}
where $C>0$ is independent of $p$.
}

\vspace{1mm}
\vspace{1mm}

\begin{proof}
Take $R_0=2R_1$, with
$R_1$ as the constant of Lemma 3.4.
Since  $\xi=(\xi_1,\ldots,\xi_m)\in\mathcal{O}_p$ and $\varepsilon\mu_i=o(1/p^\kappa)$
for $p$ large enough, we find $B_{R_1\mu_i}(\xi'_i)$ disjointed.
Let $h$ be bounded and $\phi$ be a bounded solution to (\ref{3.11}) satisfying (\ref{3.12}).
We first consider the following inner norm of $\phi$:
$$
\aligned
\|\phi\|_i=\sup_{y\in\overline{\Omega}_p\cap\left(\bigcup_{i=1}^mB_{R_1\mu_i}(\xi'_i)\right)}
|\phi(y)|,
\endaligned
$$
and claim that there is a constant $C>0$ independent of  $p$ such that
\begin{equation}\label{3.14}
\aligned
\|\phi\|_{L^{\infty}(\Omega_p)}\leq C\left(\|\phi\|_i+
\|h\|_{*}\right).
\endaligned
\end{equation}
Indeed, set
$$
\aligned
\widetilde{\phi}(y)=C_1\left(\|\phi\|_i+
\|h\|_{*}
\right)\psi(y)
\,\qquad\forall\,\,\,y\in\overline{\Omega}_p\setminus\bigcup_{i=1}^mB_{R_1\mu_i}(\xi'_i),
\endaligned
$$
where $\psi$ is the positive, uniformly bounded barrier constructed by the previous lemma and
the constant $C_1>0$ is chosen larger if necessary,  independent of $p$.
Then for
$y\in\Omega_p\setminus\bigcup_{i=1}^mB_{R_1\mu_i}(\xi'_i)$,
$$
\aligned
\mathcal{L}(\widetilde{\phi}\pm\phi)(y)\geq C_{1}\,\|h\|_{*}\left\{
\sum_{i=1}^m\frac{\mu_i^\sigma}{|y-\xi'_i|^{2+\sigma}}+\varepsilon^2
\right\}\pm
h(y)\geq|h(y)|\pm h(y)\geq0,
\endaligned
$$
for
$y\in\partial\Omega_p\setminus\bigcup_{i=1}^mB_{R_1\mu_i}(\xi'_i)$,
$$
\aligned
\frac{\partial}{\partial\nu}(\widetilde{\phi}\pm\phi)(y)\geq 0,
\endaligned
$$
and for
$y\in\Omega_p\cap\left(\bigcup_{i=1}^m\partial B_{R_1\mu_i}(\xi'_i)\right)$,
$$
\aligned
(\widetilde{\phi}\pm\phi)(y)>\|\phi\|_{i}\pm\phi(y)\geq
|\phi(y)|\pm\phi(y)\geq 0.
\endaligned
$$
From  the maximum principle (see
\cite{PW}), it follows that
$-\widetilde{\phi}\leq\phi\leq\widetilde{\phi}$ on
$\overline{\Omega}_p\setminus\bigcup_{i=1}^mB_{R_1\mu_i}(\xi'_i)$,
which implies estimate (\ref{3.14}).

We prove the lemma by contradiction. Assume that there are  sequences of
parameters $p_n\rightarrow+\infty$,
points $\xi^n=(\xi_1^n,\ldots,\xi_m^n)\in\mathcal{O}_{p_n}$,
functions $h_n$,  and associated solutions $\phi_n$ of
equation (\ref{3.11}) with orthogonality conditions (\ref{3.12})
such that
\begin{equation}\label{3.15}
\aligned
\|\phi_n\|_{L^{\infty}(\Omega_{p_n})}=1
\,\,\quad\,\,
\textrm{and}
\,\,\quad\,\,\|h_n\|_{*}\rightarrow0,
\,\,\quad\,\,\textrm{as}\,\,\,\,n\rightarrow+\infty.
\endaligned
\end{equation}
For each $k\in\{1,\ldots,l\}$, we have $\xi_k^n\in\Omega$ and we consider
$\widehat{\phi}^n_k(z)=\phi_n\big(\mu_k^nz+(\xi^n_k)'\big)$,
where $\mu^n=(\mu^n_1,\ldots,\mu_m^n)$,  $(\xi^n_k)'=\xi^n_k/\varepsilon_n$ and
 $\varepsilon_n=\exp\left\{-p_n/4\right\}$.
Note that
$$
\aligned
h_n(y)=
\big(-\Delta_{a(\varepsilon_n y)}\phi_n+\varepsilon_n^2\phi-W_{(\xi^n)'}\phi_n\big)\big|_{y=\mu_{k}^n z+(\xi^n_k)'}
=(\mu_k^n)^{-2}\left[
-\Delta_{\widehat{a}_n}\widehat{\phi}_k^n
+\varepsilon_n^2(\mu_k^n)^{2}\widehat{\phi}_k^n
-(\mu_k^n)^{2}\widehat{W}^n\widehat{\phi}_k^n
\right](z),
\endaligned
$$
where
$$
\aligned
\widehat{a}_n(z)=a(\varepsilon_n\mu_{k}^nz+\xi^n_k),
\qquad\qquad
\widehat{W}^n(z)=W_{(\xi^n)'}(\mu_{k}^nz+(\xi^n_k)').
\endaligned
$$
By the expansion of $W^n$ in (\ref{2.52})
and elliptic regularity,
$\widehat{\phi}^n_k$
converges uniformly over
compact sets to a bounded solution $\widehat{\phi}^{\infty}_k$ of equation
$(\ref{3.3})$, which satisfies
\begin{equation}\label{3.16}
\aligned
\int_{\mathbb{R}^2}\chi Z_j\widehat{\phi}_k^{\infty}=0
\quad\,\,\,\textrm{for}\,\,\,\,j=0,\,1,\,2.
\endaligned
\end{equation}
Thus $\widehat{\phi}^{\infty}_k$ is
a linear combination of $Z_j$, $j=0,1,2$.
Notice that $\int_{\mathbb{R}^2}\chi Z_jZ_{t}=0$ for $j\neq t$
and $\int_{\mathbb{R}^2}\chi Z_j^2>0$.
Hence (\ref{3.16}) implies  $\widehat{\phi}_k^{\infty}\equiv0$.

As for each $k\in\{l+1,\ldots,m\}$, we have $\xi_k^n\in\partial\Omega$ and we consider
$\widehat{\phi}^n_k(z)=\phi_n\big((A^n_k)^{-1}\mu_k^nz+(\xi^n_k)'\big)$,
where $A_k^n: \mathbb{R}^2\rightarrow\mathbb{R}^2$ is a rotation map
such that $A_k^n\nu_{\Omega_{p_n}}\big((\xi_k^n)'\big)=\nu_{\mathbb{R}_+^2}\big(0\big)$.
Similarly to the above argument, we can get that
$\widehat{\phi}^n_k$
converges uniformly over compact sets
to a bounded solution $\widehat{\phi}^{\infty}_k$ of equation
$(\ref{3.4})$, which satisfies
\begin{equation}\label{3.17}
\aligned
\int_{\mathbb{R}_{+}^2}\chi Z_j\widehat{\phi}_k^{\infty}=0
\quad\,\,\,\textrm{for}\,\,\,\,j=0,\,1.
\endaligned
\end{equation}
Thus $\widehat{\phi}^{\infty}_k$ is
a linear combination of $Z_j$, $j=0,1$.
Notice that $\int_{\mathbb{R}_{+}^2}\chi Z_jZ_t=0$ for $j\neq t$
and $\int_{\mathbb{R}_{+}^2}\chi Z_j^2>0$.
Hence (\ref{3.17}) implies  $\widehat{\phi}_k^{\infty}=0$.
Furthermore, we find that
$\lim_{n\rightarrow+\infty}\|\phi_n\|_i=0$.
But (\ref{3.14}) and (\ref{3.15}) tell us
$\liminf_{n\rightarrow+\infty}\|\phi_n\|_i>0$,
which is a contradiction.
\end{proof}

\vspace{1mm}
\vspace{1mm}
\vspace{1mm}

{\bf Step 3:} Proving  an a priori estimate for solutions to
(\ref{3.11}) that satisfy orthogonality conditions with respect to
$Z_{ij}$, $j=1,J_i$ only.

\vspace{1mm}
\vspace{1mm}
\vspace{1mm}

\noindent{\bf Lemma 3.6.}\,\,{\it For $p$  large enough, if
$\phi$ solves {\upshape (\ref{3.11})} and satisfies
\begin{equation}\label{3.18}
\aligned
\int_{\Omega_p}\chi_iZ_{ij}\phi=0\,\,\,\,\,
\,\,\,\,\forall\,\,i=1,\ldots,m,\,\,j=1,J_i,
\endaligned
\end{equation}
then
\begin{equation}\label{3.19}
\aligned
\|\phi\|_{L^{\infty}(\Omega_p)}\leq Cp\, \|h\|_{*},
\endaligned
\end{equation}
where $C>0$ is independent of $p$.}

\vspace{1mm}

\begin{proof}
According to the results in Lemma 3.4 of \cite{DW} and Lemma 3.2 of \cite{MW},
for simplicity  we only consider the validity of estimate (\ref{3.19})
when the $m$ concentration points $\xi=(\xi_1,\ldots,\xi_m)\in\mathcal{O}_p$ satisfy
the relation $|\xi_i-\xi_k|\leq2d$ for any $i,k=1,\ldots,m$, $i\neq k$ and for
any $d>0$ sufficiently small, fixed and independent of $p$.
Let $R>R_0+1$ be a large but fixed number. Denote for $i=1,\ldots,m$,
\begin{equation}\label{3.20}
\aligned
\widehat{Z}_{i0}(y)=Z_{i0}(y)-\frac1{\mu_i}
+a_{i0}G(\varepsilon y,\xi_i),
\endaligned
\end{equation}
where
\begin{equation}\label{3.21}
\aligned
a_{i0}=\frac1{\mu_i\big[H(\xi_i,\xi_i)-\frac{4}{c_i}\log(\varepsilon \mu_i R)\big]}.
\endaligned
\end{equation}
Note that by estimate (\ref{2.23}),  expansions (\ref{2.1c}) and (\ref{2.1e}), and definitions (\ref{2.4}), (\ref{3.2}), (\ref{3.5}) and (\ref{3.8}),
\begin{equation}\label{3.22}
\aligned
\frac18p\leq-\log(\varepsilon \mu_i R)
\leq \frac12p,
\endaligned
\end{equation}
and
\begin{equation}\label{3.23}
\aligned
\widehat{Z}_{i0}(y)=O\left(
\frac{\,G(\varepsilon y,\xi_i)\,}{p\mu_i}
\right).
\endaligned
\end{equation}
Let   $\eta_1$ and $\eta_2$ be  radial smooth cut-off functions in $\mathbb{R}^2$ such that
$$
\aligned
&0\leq\eta_1\leq1;\,\,\,\ \,\,\,|\nabla\eta_1|\leq C\,\,\ \textrm{in}\,\,\,\mathbb{R}^2;
\,\,\,\ \,\,\,\eta_1\equiv1\,\,\ \textrm{in}\,\,\,B_R(0);\,\,\,\,\ \,\,\,
\,\eta_1\equiv0\,\,\ \textrm{in}\,\,\,\mathbb{R}^2\setminus B_{R+1}(0);\\[1mm]
&0\leq\eta_2\leq1;\,\,\,\ \,\,\,|\nabla\eta_2|\leq C\,\,\ \textrm{in}\,\,\,\mathbb{R}^2;
\,\,\,\ \,\,\,\eta_2\equiv1\,\,\ \textrm{in}\,\,\,B_{3d}(0);\,\,\,
\,\,\,\,\,\,\eta_2\equiv0\,\,\ \textrm{in}\,\,\,\mathbb{R}^2\setminus B_{6d}(0).
\endaligned
$$
We set, for  $i\in\{1,\ldots,l\}$,
\begin{equation}\label{3.24}
\aligned
\eta_{i1}(y)=
\eta_1\left(\frac{1}{\mu_i}\big|y-\xi_i'\big|\right),
\,\,\quad\quad\,\,
\eta_{i2}(y)=
\eta_2\left(\varepsilon\big|y-\xi'_i\big|\right),
\endaligned
\end{equation}
and for  $i\in\{l+1,\ldots,m\}$,
\begin{equation}\label{3.25}
\aligned
\eta_{i1}(y)=
\eta_1\left(\frac{1}{\mu_i}\big|F_i^p(y)\big|\right),
\,\,\quad\quad\,\,
\eta_{i2}(y)=
\eta_2\left(\varepsilon\big|F_i^p(y)\big|\right).
\endaligned
\end{equation}

Now we define the test function
\begin{equation}\label{3.26}
\aligned
\widetilde{Z}_{i0}(y)=\eta_{i1}Z_{i0}+(1-\eta_{i1})\eta_{i2}\widehat{Z}_{i0}.
\endaligned
\end{equation}
Given $\phi$ satisfying (\ref{3.11}) and (\ref{3.18}), let
\begin{equation}\label{3.27}
\aligned
\widetilde{\phi}=\phi+\sum\limits_{i=1}^{m}d_i\widetilde{Z}_{i0}+\sum_{i=1}^m\sum\limits_{j=1}^{J_i}e_{ij}\chi_iZ_{ij}.
\endaligned
\end{equation}
We will first prove the existence of
$d_i$ and $e_{ij}$  such
that $\widetilde{\phi}$ satisfies the orthogonality condition
\begin{equation}\label{3.28}
\aligned
\int_{\Omega_p}\chi_iZ_{ij}\widetilde{\phi}=0
\,\,\quad\,\,\forall\,\,i=1,\ldots,m,\,\,j=0,1,J_i.
\endaligned
\end{equation}
Multiplying  (\ref{3.27}) by $\chi_iZ_{ij}$, $i=1,\ldots,m$,
$j=0,1,J_i$  and using
orthogonality conditions (\ref{3.18}) and (\ref{3.28})  together with the fact that
$\chi_i\chi_k\equiv0$ if $i\neq k$, we get
\begin{equation}\label{3.29}
\aligned
d_i\int_{\Omega_p}\chi_iZ_{i0}\widetilde{Z}_{i0}
+\sum_{k\neq i}^md_k\int_{\Omega_p}\chi_iZ_{i0}\widetilde{Z}_{k0}
+\sum\limits_{t=1}^{J_i}e_{it}\int_{\Omega_p}\chi_i^2Z_{i0}Z_{it}
=-\int_{\Omega_p}\chi_iZ_{i0}\phi,
\endaligned
\end{equation}
\begin{equation}\label{3.30}
\aligned
d_i\int_{\Omega_p}\chi_iZ_{ij}\widetilde{Z}_{i0}
+\sum_{k\neq i}^md_k\int_{\Omega_p}\chi_iZ_{ij}\widetilde{Z}_{k0}
+\sum_{t=1}^{J_i}e_{it}\int_{\Omega_p}\chi^2_iZ_{ij}Z_{it}=0,
\,\quad\ \ \,\ \,j=1,\,J_i.
\endaligned
\end{equation}
Remark that for any $i=1,\ldots,l$,
$\widetilde{Z}_{i0}$
coincides with $Z_{i0}$ in $B_{R\mu_i}(\xi'_i)$,
while for any $i=l+1,\ldots,m$, $\widetilde{Z}_{i0}$
coincides with $Z_{i0}$ in the region
$\{y\in\Omega_p:|F_i^p(y)|\leq R\mu_i\}$.
Moreover, from definitions (\ref{3.6})-(\ref{3.7}) we can write
$z=\frac{1}{\mu_i}F^p_i(y)$ and its inverse
$y=\big(\frac{1}{\mu_i}F^p_i\big)^{-1}(z)=\xi'_i+\frac1{\varepsilon}A_i^{-1}F^{-1}_{i}(\varepsilon\mu_i z)$
such that
$\det\big(\nabla\big(\frac{1}{\mu_i}F_i^p\big)^{-1}(z)\big)=\mu_i^2+O(\varepsilon\mu_i^3|z|)$
holds in the upper half-ball $B^+_{R\mu_i}(0):=B_{R\mu_i}(0)\cap\mathbb{R}^2_{+}$.
Then for any $i=1,\ldots,l,\,$   $j=1,2\,$ and $\,t=1,2$,
$$
\aligned
\int_{\Omega_p}\chi_iZ_{i0}\widetilde{Z}_{i0}=
\int_{\mathbb{R}^2}\chi Z^2_{0}=C_0>0,
\qquad\qquad\qquad\qquad
\int_{\Omega_p}\chi_i^2Z_{i0}Z_{it}=
\int_{\mathbb{R}^2}\chi^2 Z_{0}Z_t=0,
\endaligned
$$
$$
\aligned
\int_{\Omega_p}\chi_iZ_{ij}\widetilde{Z}_{i0}=
\int_{\mathbb{R}^2}\chi Z_{j}Z_{0}=0,
\qquad\qquad\qquad\qquad
\int_{\Omega_p}\chi^2_iZ_{ij}Z_{it}=
\int_{\mathbb{R}^2}\chi^2 Z_{j}Z_{t}
=C_j\delta_{jt},
\endaligned
$$
where $\delta_{jt}$ denotes the Kronecker's symbol, but for any $i=l+1,\ldots,m$ and $j=t=J_i=1$,
$$
\aligned
\int_{\Omega_p}\chi_iZ_{i0}\widetilde{Z}_{i0}=
\int_{\mathbb{R}_{+}^2}\chi Z^2_{0}[1+O\big(\varepsilon\mu_i|z|\big)]=\frac{C_0}2+O\left(\varepsilon\mu_i\right),
\qquad
\int_{\Omega_p}\chi_i^2Z_{i0}Z_{i1}=
\int_{\mathbb{R}_{+}^2}\chi^2 Z_{0}Z_1[1+O\big(\varepsilon\mu_i|z|\big)]=O\left(\varepsilon\mu_i\right),
\endaligned
$$
$$
\aligned
\int_{\Omega_p}\chi_iZ_{i1}\widetilde{Z}_{i0}=
\int_{\mathbb{R}_{+}^2}\chi  Z_{1} Z_{0}  [1+O\big(\varepsilon\mu_i|z|\big) ] =O\left(\varepsilon\mu_i\right),
\,\,\qquad\,\,
\int_{\Omega_p}\chi^2_iZ_{i1}^2=
\int_{\mathbb{R}_{+}^2}\chi^2  Z_{1}^2  [1+O\big(\varepsilon\mu_i|z|\big) ]
=\frac{C_1}{2} +O\left(\varepsilon\mu_i\right).
\endaligned
$$
Moreover, from  (\ref{3.23}) and (\ref{3.26})
it follows  that  for any $i=1,\ldots,m$ and
$j=0,1,J_i$,
$$
\aligned
\int_{\Omega_p}\chi_iZ_{ij}\widetilde{Z}_{k0}=
O\left(\frac{\mu_i\log p}{p\mu_k}
\right),\,
\quad\,\forall\,\,k\neq i.
\endaligned
$$
Thus by (\ref{3.30}),
\begin{equation}\label{3.31}
\aligned
e_{ij}=\left(-d_i\int_{\Omega_p}\chi_iZ_{ij}\widetilde{Z}_{i0}
-\sum_{k\neq i}^md_k\int_{\Omega_p}\chi_iZ_{ij}\widetilde{Z}_{k0}
\right)
\left/
\int_{\Omega_p}\chi^2_iZ^2_{ij},
\right.
\qquad
i=1,\ldots,m,\,\,\,
j=1,J_i.
\endaligned
\end{equation}
Furthermore,
\begin{equation}\label{3.32}
\aligned
|e_{ij}|\leq\left\{\aligned
&C\sum_{k\neq i}^m|d_k|\frac{\mu_i\log p}{p\mu_k},
\,\quad\qquad\qquad\quad\,\,\,\,\forall\,\,\,i=1,\ldots,l,\,\ \ \,\,\,j=1,2,
\\
&C\varepsilon\mu_i|d_i|+
C\sum_{k\neq i}^m|d_k|\frac{\mu_i\log p}{p\mu_k},
\quad\quad\forall\,\,\,i=l+1,\ldots,m,\,\,j=1.
\endaligned
\right.
\endaligned
\end{equation}
We need just to show that $d_i$ is well defined. From (\ref{3.29}) we can easily get
that for any $i=1,\ldots,l$,
\begin{equation}\label{3.33}
\aligned
d_iC_0
+\sum_{k\neq i}^md_kO\left(\frac{\mu_i\log p}{p\mu_k}
\right)
=-\int_{\Omega_p}\chi_iZ_{i0}\phi,
\endaligned
\end{equation}
and for any  $i=l+1,\ldots,m$,
\begin{equation}\label{3.34}
\aligned
\frac12d_iC_0\big[1+O\big(\varepsilon\mu_i\big)\big]
+\sum_{k\neq i}^md_kO\left(\frac{\mu_i\log p}{p\mu_k}
\right)+e_{i1}O\big(\varepsilon\mu_i\big)
=-\int_{\Omega_p}\chi_iZ_{i0}\phi,
\endaligned
\end{equation}
where $e_{i1}$ is defined in (\ref{3.31}) and satisfies
$$
\aligned
e_{i1}=d_iO\left(\varepsilon\mu_i\right)+
\sum_{k\neq i}^md_kO\left(\frac{\mu_i\log p}{p\mu_k}\right).
\endaligned
$$
We denote $\mathcal{A}$ the coefficient matrix of equations (\ref{3.33})-(\ref{3.34}). By the
above estimates, it is clear that
$\mathcal{M}^{-1}\mathcal{A}\mathcal{M}$ is diagonally dominant and thus invertible, where
$\mathcal{M}=\diag(\mu_1,\ldots,\mu_m)$. Hence $\mathcal{A}$ is also invertible and
$(d_1,\ldots,d_m)$ is well defined.

Estimate (\ref{3.19}) is a  direct consequence of the following two claims.

\vspace{1mm}
\vspace{1mm}
\vspace{1mm}

\noindent{\bf Claim 1.}\,\,{\it
Let $\mathcal{L}=-\Delta_{a(\varepsilon y)}+\varepsilon^2-W_{\xi'}$, then for any $i=1,\ldots,m$
and $j=1, J_i$,
\begin{equation}\label{3.35}
\aligned
\big\|\mathcal{L}(\chi_iZ_{ij})\big\|_{*}\leq\frac{C}{\mu_i},
\,\quad\quad\quad\,\,\,\quad\,
\big\|\mathcal{L}(\widetilde{Z}_{i0})\big\|_{*}\leq
C\frac{\log p}{p\mu_i}.
\endaligned
\end{equation}
}

\noindent{\bf Claim 2.}\,\,{\it
For any $i=1,\ldots,m$ and $j=1,J_i$,
\begin{equation*}\label{3.36}
\aligned
|d_i|\leq Cp\mu_i\|h\|_{*},
\,\quad\qquad\quad\,\,\quad
|e_{ij}|\leq C\mu_i\log p\,\|h\|_{*}.
\endaligned
\end{equation*}
}

In fact,  the definition of
$\widetilde{\phi}$ in (\ref{3.27}) tells us
\begin{equation}\label{3.37}
\aligned\left\{\aligned
&
\mathcal{L}(\widetilde{\phi})=h+\sum\limits_{i=1}^{m}d_i\mathcal{L}(\widetilde{Z}_{i0})
+\sum_{i=1}^m\sum_{j=1}^{J_i}e_{ij}\mathcal{L}(\chi_iZ_{ij})
\,\quad\,\textrm{in}\,\,\,\ \,\,\Omega_p,\\
&\frac{\partial\widetilde{\phi}}{\partial\nu}=0
\,\qquad\qquad\qquad\qquad\quad\qquad\,
\,\qquad\qquad\qquad\,\,
\,\textrm{on}\,\,\,\,\,\partial\Omega_p.
\endaligned\right.
\endaligned
\end{equation}
Then by Lemma 3.5, we obtain
\begin{equation}\label{3.38}
\aligned
\|\widetilde{\phi}\|_{L^{\infty}(\Omega_p)}\leq
C\left\{\|h\|_{*}
+\sum\limits_{i=1}^{m}|d_i|\big\|\mathcal{L}(\widetilde{Z}_{i0})\big\|_{*}
+\sum_{i=1}^m\sum_{j=1}^{J_i}|e_{ij}|\big\|\mathcal{L}(\chi_i Z_{ij})\big\|_{*}
\right\}\leq C \log p\,\|h\|_{*}.
\endaligned
\end{equation}
Using the definition of $\widetilde{\phi}$
again and the fact that
\begin{equation}\label{3.39}
\aligned
\big\|\widetilde{Z}_{i0}\big\|_{L^{\infty}(\Omega_p)}\leq\frac{C}{\mu_i}
\quad\quad\,\textrm{and}\,\quad\quad
\big\|\chi_iZ_{ij}\big\|_{L^{\infty}(\Omega_p)}\leq\frac{C}{\mu_i},
\,\,\quad\,\forall\,\,\,i=1,\ldots,m,\,\,j=1,J_i,
\endaligned
\end{equation}
estimate (\ref{3.19}) then follows from estimate (\ref{3.38}) and Claim 2.

\vspace{1mm}
\vspace{1mm}
\vspace{1mm}
\vspace{1mm}

\noindent{\bf Proof of Claim 1.}
Observe that
\begin{equation}\label{3.40}
\aligned
\mathcal{L}(\chi_iZ_{ij})=\,
\chi_i\mathcal{L}(Z_{ij})
-2\nabla\chi_i\nabla Z_{ij}
-Z_{ij}\big[\Delta\chi_i
+\varepsilon\nabla\log a(\varepsilon y)\nabla\chi_i\big].
\endaligned
\end{equation}
For any $i=1,\ldots,l$ and $j=1,2$, we write $z_i:=y-\xi'_i$ and
note that in the region
$|z_i|=|y-\xi'_i|\leq\mu_i(R_0+1)$, by  (\ref{2.52}), (\ref{3.2}), (\ref{3.3}) and (\ref{3.5}),
$$
\aligned
\mathcal{L}(Z_{ij})=&
\big(-\Delta-W_{\xi'}\big)\left[\frac1{\mu_i}Z_{j}\left(\frac{y-\xi'_i}{\mu_i}\right)
\right]
-\varepsilon\nabla\log a(\varepsilon y)\nabla\left[\frac1{\mu_i}Z_{j}\left(\frac{y-\xi'_i}{\mu_i}\right)
\right]+\frac{\varepsilon^2}{\mu_i}Z_{j}\left(\frac{y-\xi'_i}{\mu_i}\right)\\[1mm]
=&
O\left(\frac1{p\mu_i}\cdot\frac{8\mu_i^2}{(\mu_i^2+|y-\xi'_i|^2)^2}\right)
+O\left(
\frac{\varepsilon}{\mu_{i}^2+|y-\xi'_i|^2}
\right)+O\left(
\frac{\varepsilon^2}{(\mu_i^2+|y-\xi'_i|^2)^{1/2}}
\right),
\endaligned
$$
and then, by (\ref{3.40}),
\begin{equation}\label{3.41}
\aligned
\mathcal{L}(\chi_iZ_{ij})
=O\left(\frac{\mu_ip^{-1}}{(\mu_i^2+|y-\xi'_i|^2)^2}\right)
+O\left(
\frac{\mu_i^{-1}}{\mu_i^2+|y-\xi'_i|^2}
\right)
+O\left(
\frac{\mu_i^{-2}}{(\mu_{i}^2+|y-\xi'_i|^2)^{1/2}}
\right).
\endaligned
\end{equation}
As for any $i=l+1,\ldots,m$,
owing to $F_i^p(\xi'_i)=(0,0)$ and $\nabla F_i^p(\xi'_i)=A_i$, we know
\begin{equation}\label{3.42}
\aligned
\nabla_y=A_i\nabla_{z_i}+O(\varepsilon|z_i|)
\nabla_{z_i}
\,\quad\,\,\,\quad\,
\textrm{and}
\,\quad\,\,\,\quad\,
-\Delta_y=-\Delta_{z_i}+O(\varepsilon|z_i|)
\nabla_{z_i}^2+O(\varepsilon)\nabla_{z_i},
\endaligned
\end{equation}
where
\begin{equation}\label{3.43}
\aligned
z_i:=F_i^p(y)=\frac{1}{\varepsilon}F_i\big(A_i(\varepsilon y-\xi_i)\big)
=A_i(y-\xi_i')\big\{1+O\big(\varepsilon A_i(y-\xi_i')\big)\big\}.
\endaligned
\end{equation}
In the region
$|z_i|=|F_i^p(y)|\leq\mu_i(R_0+1)$, by (\ref{2.52}), (\ref{3.2}), (\ref{3.4}), (\ref{3.8}),
(\ref{3.42}) and (\ref{3.43}),
$$
\aligned
\mathcal{L}(Z_{i1})=&
\big(-\Delta_{z_i}-W_{\xi'}\big)\left[\frac1{\mu_i}Z_{1}\left(\frac{z_i}{\mu_i}\right)
\right]
+
\big(\varepsilon^2+O(\varepsilon|z_i|)
\nabla_{z_i}^2+O(\varepsilon)\nabla_{z_i}\big)\left[\frac1{\mu_i}Z_{1}\left(\frac{z_i}{\mu_i}\right)
\right]
\\[1mm]
=&
O\left(\frac1{p\mu_i}\cdot\frac{8\mu_i^2}{(\mu_i^2+|z_i|^2)^2}\right)
+O\left(
\frac{\varepsilon}{\mu_{i}^2+|z_i|^2}
\right)+O\left(
\frac{\varepsilon^2}{(\mu_i^2+|z_i|^2)^{1/2}}
\right).
\endaligned
$$
Thus by (\ref{3.40}),
\begin{equation}\label{3.44}
\aligned
\mathcal{L}(\chi_iZ_{i1})
=O\left(\frac{\mu_ip^{-1}}{(\mu_i^2+|F_i^p(y)|^2)^2}\right)
+O\left(
\frac{\mu_i^{-1}}{\mu_i^2+|F_i^p(y)|^2}
\right)
+O\left(
\frac{\mu_i^{-2}}{(\mu_{i}^2+|F_i^p(y)|^2)^{1/2}}
\right).
\endaligned
\end{equation}
Hence by (\ref{3.41}), (\ref{3.43}), (\ref{3.44}) and
the definition of  $\|\cdot\|_*$ in (\ref{2.37}),  we obtain
$\big\|\mathcal{L}(\chi_iZ_{ij})\big\|_{*}=O\left(1/\mu_i\right)$
for all $i=1,\ldots,m$ and $j=1,J_i$.

We now prove the second inequality in (\ref{3.35}).
In fact,
\begin{eqnarray}\label{3.45}
&&\mathcal{L}(\widetilde{Z}_{i0})=\eta_{i1}\mathcal{L}(Z_{i0}-\widehat{Z}_{i0})
+\eta_{i2}\mathcal{L}(\widehat{Z}_{i0})
-(Z_{i0}-\widehat{Z}_{i0})\Delta_{a(\varepsilon y)}\eta_{i1}
-2\nabla\eta_{i1}\nabla(Z_{i0}-\widehat{Z}_{i0})
\nonumber\\[1mm]
&&\qquad\qquad\,-2\nabla\eta_{i2}\nabla\widehat{Z}_{i0}-\widehat{Z}_{i0}\Delta_{a(\varepsilon y)}\eta_{i2}.
\end{eqnarray}
Recalling that $z_i=y-\xi'_i$  for any $i=1,\ldots,l$, but
$z_i=F_i^p(y)$  for any $i=l+1,\ldots,m$,  we
now consider the four regions
$$
\aligned
\Omega_{1}=\left\{y\in\Omega_p\big|\,|z_i|\leq\mu_iR\right\},
\,\quad\,\quad\quad\quad\quad\quad\quad\quad\,
\Omega_{2}=\left\{y\in\Omega_p\big|\,\mu_iR<|z_i|\leq\mu_i(R+1)\right\},\\[1mm]
\Omega_{3}=\left\{y\in\Omega_p\left|\,\mu_i(R+1)<|z_i|\leq3d/\varepsilon\right.\right\},
\,\quad\quad\quad\quad\,\,\,\,\,\,\,
\Omega_{4}=\left\{y\in\Omega_p\left|\,3d/\varepsilon<|z_i|\leq
6d/\varepsilon\right.\right\}.
\quad\quad
\endaligned
$$
Notice first that, by (\ref{3.2}), (\ref{3.5}), (\ref{3.8}) and (\ref{3.43}),
\begin{equation}\label{3.46}
\aligned
\left|
Z_{i0}-\frac{1}{\mu_i}
\right|=\frac{2\mu_i}{\mu_i^2+|z_i|^2}
=O\left(\frac{\mu_i}{(\mu_i+|y-\xi_i'|)^2}\right),
\endaligned
\end{equation}
and for $\mu_iR<|z_i|\leq6d/\varepsilon$, by (\ref{3.20})-(\ref{3.21}),
\begin{equation}\label{3.47}
\aligned
Z_{i0}-\widehat{Z}_{i0}=\frac1{\mu_i}
-a_{i0}G(\varepsilon y,\xi_i)=
\frac1{\mu_i\big[H(\xi_i,\xi_i)-\frac{4}{c_i}\log(\varepsilon \mu_i R)\big]}\left[
\frac{4}{c_i}\log\frac{|y-\xi_i'|}{\mu_iR}+O\big(
\varepsilon^\alpha|y-\xi_i'|^\alpha\big)
\right].
\endaligned
\end{equation}
In $\Omega_1$, by (\ref{3.45}),
$$
\aligned
\mathcal{L}(\widetilde{Z}_{i0})=\mathcal{L}(Z_{i0})
=\left[-\Delta_{a(\varepsilon y)}Z_{i0}-\frac{8\mu_i^2}{(\mu_i^2+|y-\xi'_i|^2)^2}Z_{i0}
\right]+\left[\frac{8\mu_i^2}{(\mu_i^2+|y-\xi'_i|^2)^2}-W_{\xi'}\right]Z_{i0}
+\varepsilon^2Z_{i0}.
\endaligned
$$
Note that, by (\ref{2.52}),
\begin{equation}\label{3.48}
\aligned
\left[\frac{8\mu_i^2}{(\mu_i^2+|y-\xi'_i|^2)^2}-W_{\xi'}\right]Z_{i0}
=O\left(\frac{\mu_ip^{-1}}{(\mu_i^2+|y-\xi'_i|^2)^2}\right),
\qquad\forall\,\,y\in\Omega_1\cup\Omega_2.
\endaligned
\end{equation}
Hence in $\Omega_1$, by  (\ref{3.9}),
\begin{equation}\label{3.49}
\aligned
\mathcal{L}(\widetilde{Z}_{i0})=\mathcal{L}(Z_{i0})=
O\left(\frac{1}{p\mu_i^3}\right).
\endaligned
\end{equation}
In $\Omega_2$, by (\ref{3.20}) and (\ref{3.45}),
\begin{eqnarray*}
\mathcal{L}(\widetilde{Z}_{i0})
=\left[-\Delta_{a(\varepsilon y)}Z_{i0}-\frac{8\mu_i^2}{(\mu_i^2+|y-\xi'_i|^2)^2}Z_{i0}
\right]+\left[\frac{8\mu_i^2}{(\mu_i^2+|y-\xi'_i|^2)^2}-W_{\xi'}\right]Z_{i0}
+\varepsilon^2\left(
Z_{i0}-\frac1{\mu_i}
\right)
&&\nonumber\\
+W_{\xi'}(1-\eta_{i1})(Z_{i0}-\widehat{Z}_{i0})
+\frac{\varepsilon^2}{\mu_i}\eta_{i1}
-2\nabla\eta_{i1}\nabla(Z_{i0}-\widehat{Z}_{i0})-(Z_{i0}-\widehat{Z}_{i0})\Delta_{a(\varepsilon y)}\eta_{i1}.
\qquad\quad
&&
\end{eqnarray*}
Note that in $\Omega_2$, by (\ref{3.22}) and (\ref{3.47}),
\begin{equation}\label{3.50}
\aligned
|Z_{i0}-\widehat{Z}_{i0}|=O\left(
\frac{1}{p\mu_iR}
\right)
\,\,\quad\quad\,\textrm{and}\,\,\quad\quad\,
|\nabla\big(Z_{i0}-\widehat{Z}_{i0}\big)|=O\left(
\frac{1}{p\mu_i^2R}
\right).
\endaligned
\end{equation}
Moreover $|\nabla\eta_{i1}|=O\left(\mu_i^{-1}\right)$, $|\Delta_{a(\varepsilon y)}\eta_{i1}|=O\left(\mu_i^{-2}\right)$.
Hence in $\Omega_2$, by (\ref{2.52}), (\ref{3.9}), (\ref{3.46}),  (\ref{3.48}) and (\ref{3.50}),
\begin{equation}\label{3.51}
\aligned
\mathcal{L}(\widetilde{Z}_{i0})
=O\left(
\frac{1}{p\mu_i^3R}
\right).
\endaligned
\end{equation}
In $\Omega_3$, by (\ref{3.9}), (\ref{3.20}), (\ref{3.45}) and (\ref{3.46}),
$$
\aligned
\mathcal{L}(\widetilde{Z}_{i0})=&\mathcal{L}(\widehat{Z}_{i0})
=\mathcal{L}(Z_{i0})-\mathcal{L}(Z_{i0}-\widehat{Z}_{i0})
\\[1.6mm]
=&\left[\frac{8\mu_i^2}{(\mu_i^2+|y-\xi'_i|^2)^2}-W_{\xi'}\right]Z_{i0}
+W_{\xi'}\left[\frac1{\mu_i}
-a_{i0}G(\varepsilon y, \xi_i)\right]
+O\left(\frac{\varepsilon\mu_i}{
(\mu_i+|y-\xi'_i|)^3}\right)
+O\left(\frac{\varepsilon^2\mu_i}{(\mu_i+|y-\xi_i'|)^2}\right)
\\[1.3mm]
\equiv&A_1+A_2+O\left(\frac{\varepsilon\mu_i}{
(\mu_i+|y-\xi'_i|)^3}\right)
+O\left(\frac{\varepsilon^2\mu_i}{(\mu_i+|y-\xi_i'|)^2}\right).
\endaligned
$$
To estimate the first two terms, we need to decompose $\Omega_3$ into some subregions:
$$
\aligned
&\Omega_{3,i}=\left\{y\in\Omega_3\left|\,\,\mu_i(R+1)<|z_i|\leq
\mu_i^{1/2}\big/\big(3\varepsilon^{1/2} p^{2\kappa}\big)\,\right.\right\},\\
\Omega_{3,k}=&\left\{y\in\Omega_3\left|\,\,|z_k|\leq\mu_k^{1/2}\big/\big(3\varepsilon^{1/2} p^{2\kappa}\big)\,\right.\right\},
\,\,\,\,k\neq i,\quad\,\,\textrm{and}\quad\,\,
\widetilde{\Omega}_3=\Omega_3\setminus\bigcup_{t=1}^m\Omega_{3,t}.
\endaligned
$$
From (\ref{2.51}), (\ref{2.52}) and (\ref{3.43}) we get
$$
\aligned
A_1=\left\{
\aligned
&\frac{8\mu^2_i}{(\mu_{i}^2+|y-\xi'_i|^2)^2}
\left[\frac1p\left(\omega_1-U_{1,0}-\frac12U_{1,0}^2
\right)\left(\frac{y-\xi'_i}{\mu_i}\right)
+O\left(\frac{\log^4\big(\big|\frac{y-\xi'_i}{\mu_i}\big|+2\big)}{p^2}\right)\right]O\left(\frac1{\mu_i}\right)
\,\,\,\textrm{in}\,\,\,\Omega_{3,i},\\[2mm]
&O\left(\varepsilon^2p^{8\kappa}/\mu_i\right)
\quad\qquad\qquad\qquad\qquad\qquad\,
\qquad\qquad\qquad\qquad\qquad\qquad
\qquad\qquad\qquad\qquad\qquad\,
\textrm{in}\,\,\,\,\widetilde{\Omega}_{3}.
\endaligned
\right.
\endaligned
$$
Moreover, by (\ref{3.22}) and (\ref{3.47}),
$$
\aligned
A_2=
\frac{8\mu^2_i}{(\mu_i^2+|y-\xi'_i|^2)^2} O\left(\frac{\log|y-\xi_i'|-\log\mu_iR+\varepsilon^\alpha|y-\xi_i'|^\alpha}{p\mu_i}\right)
\,\,\,\,\,\,\textrm{in}\,\,\,\Omega_{3,i}\cup \widetilde{\Omega}_{3}.
\endaligned
$$
Then in $\Omega_{3,i}\cup \widetilde{\Omega}_{3}$,
by (\ref{2.1}) and (\ref{2.8}),
\begin{equation}\label{3.52}
\aligned
\mathcal{L}(\widetilde{Z}_{i0})=O
\left(
\frac{1}{(\mu_i^2+|y-\xi'_i|^2)^{(3+\sigma)/2}}\cdot\frac{\mu_i^\sigma}{p}
\right)
+
O
\left(\frac{\log|y-\xi'_i|-
\log\mu_iR}{(\mu_i^2+|y-\xi'_i|^2)^2}\cdot
\frac{\mu_i}{p}
\right)
\endaligned
\end{equation}
with $\sigma>0$  small but fixed,  independent of $p$.
In $\Omega_{3,k}$ with $k\neq i$, by (\ref{2.52}), (\ref{3.9}), (\ref{3.23})
and (\ref{3.46}),
\begin{eqnarray}\label{3.53}
\mathcal{L}(\widetilde{Z}_{i0})=\mathcal{L}(\widehat{Z}_{i0})=
\frac{8\mu_i^2}{(\mu_{i}^2+|y-\xi'_i|^2)^2}Z_{i0}
-\left[
\Delta_{a(\varepsilon y)} Z_{i0}+\frac{8\mu_i^2}{(\mu_{i}^2+|y-\xi'_i|^2)^2}Z_{i0}
\right]
+\varepsilon^2\left(
Z_{i0}-\frac1{\mu_i}
\right)-W_{\xi'}\widehat{Z}_{i0}
&&\nonumber\\[1.5mm]
=O\left(\frac{8\mu^2_k}{(\mu^2_k+|y-\xi'_k|^2)^2}
\cdot\frac{\log p}{p\mu_i}
\right).
\qquad\qquad\qquad\qquad\qquad\qquad\qquad\qquad
\qquad\quad\qquad\qquad\qquad\qquad\,\,\,\,\,
&&
\end{eqnarray}
In $\Omega_4$, by (\ref{3.45}),
$$
\aligned
\mathcal{L}(\widetilde{Z}_{i0})=&\frac{8\mu_i^2}{(\mu_i^2+|y-\xi'_i|^2)^2}\eta_{i2}Z_{i0}
-\eta_{i2}\left[\Delta_{a(\varepsilon y)} Z_{i0}
+\frac{8\mu_i^2}{(\mu_i^2+|y-\xi'_i|^2)^2}Z_{i0}
\right]+\varepsilon^2\eta_{i2}\left(
Z_{i0}-\frac1{\mu_i}
\right)
\\[1.5mm]
&-\eta_{i2}W_{\xi'}\widehat{Z}_{i0}
-2\nabla\eta_{i2}\nabla \widehat{Z}_{i0}-\widehat{Z}_{i0}\Delta_{a(\varepsilon y)}\eta_{i2}.
\endaligned
$$
By (\ref{2.23}) and (\ref{2.51}) we find $W_{\xi'}=O(\varepsilon^{4-\sigma})$ in $\Omega_4$.
In addition,
$|\nabla\eta_{i2}|=O\left(\varepsilon/d\right)$,
$|\Delta_{a(\varepsilon y)}\eta_{i2}|=O\left(\varepsilon^2/d^2\right)$,
\begin{equation}\label{3.54}
\aligned
|\widehat{Z}_{i0}|=O\left(
\frac{|\log d|}{p\mu_i}
\right)\,\,\quad\quad\,\textrm{and}\,\,\quad\quad\,
|\nabla\widehat{Z}_{i0}|=O\left(
\frac{\varepsilon}{pd\mu_i}
\right).
\endaligned
\end{equation}
Hence by (\ref{3.9}) and (\ref{3.46}),
\begin{equation}\label{3.55}
\aligned
\mathcal{L}(\widetilde{Z}_{i0})
=O\left(
\frac{\varepsilon^2|\log d|}{pd^2\mu_i}
\right).
\endaligned
\end{equation}
Combining  (\ref{2.37}), (\ref{3.49}), (\ref{3.51}), (\ref{3.52}), (\ref{3.53}) and (\ref{3.55}), we arrive at
$$
\aligned
\big\|\mathcal{L}(\widetilde{Z}_{i0})\big\|_{*}=O\left(
\frac{\log p}{p\mu_i}\right),
\,\,\quad\,\forall\,\,i=1,\ldots,m.
\endaligned
$$

\vspace{1mm}
\vspace{1mm}

\noindent{\bf Proof of Claim 2.}
Multiplying equation (\ref{3.37}) by
$a(\varepsilon y)\widetilde{Z}_{i0}$, integrating by parts and using the
relations (\ref{3.38})-(\ref{3.39}), we
 get
$$
\aligned
\sum_{k=1}^md_k&
\int_{\Omega_p}a(\varepsilon y)\widetilde{Z}_{k0}\mathcal{L}(\widetilde{Z}_{i0})
\\
=&-\int_{\Omega_p} a(\varepsilon y)h\widetilde{Z}_{i0}
+\int_{\Omega_p}a(\varepsilon y)\widetilde{\phi}\mathcal{L}(\widetilde{Z}_{i0})-\sum_{k=1}^m\sum_{t=1}^{J_k}e_{kt}
\int_{\Omega_p}a(\varepsilon y)\chi_kZ_{kt}\mathcal{L}(\widetilde{Z}_{i0})
\\[1mm]
\leq&C\frac{\|h\|_{*}}{\mu_i}
+C\big\|\mathcal{L}(\widetilde{Z}_{i0})\big\|_{*}
\left(\|\widetilde{\phi}\|_{L^{\infty}(\Omega_p)}
+\sum_{k=1}^m\sum_{t=1}^{J_k}\frac{1}{\mu_k}|e_{kt}|\right)
\\[1mm]
\leq&C\frac{\|h\|_{*}}{\mu_i}
+
C\big\|\mathcal{L}(\widetilde{Z}_{i0})\big\|_{*}
\left[\|h\|_{*}+
\sum\limits_{k=1}^{m}|d_k|\big\|\mathcal{L}(\widetilde{Z}_{k0})\big\|_{*}
+\sum_{k=1}^m\sum_{t=1}^{J_k}|e_{kt}|\left(\frac{1}{\mu_k}
+\big\|\mathcal{L}(\chi_kZ_{kt})\big\|_{*}
\right)
\right],
\endaligned
$$
where we have applied the following  inequality
$$
\aligned
\int_{\Omega_p}\frac{\mu_i^\sigma}{(|y-\xi'_i|+\mu_i)^{2+\sigma}}\leq C,
\,\,\,\quad\,\,\,\forall\,\,\,i=1,\ldots,m.
\endaligned
$$
But  estimates (\ref{3.32}) and  (\ref{3.35})  imply  that for any $i=1,\ldots,m$,
\begin{equation}\label{3.56}
\aligned
|d_i|\left|
\int_{\Omega_p}a(\varepsilon y)\widetilde{Z}_{i0}\mathcal{L}(\widetilde{Z}_{i0})
\right|
\leq
C\frac{\|h\|_{*}}{\mu_i}
+C\sum_{k=1}^m\frac{|d_k|\log^2p}{p^2\mu_i\mu_k}
+\sum_{k\neq i}^m\left|d_k
\int_{\Omega_p}a(\varepsilon y)\widetilde{Z}_{k0}\mathcal{L}(\widetilde{Z}_{i0})
\right|.
\endaligned
\end{equation}
From (\ref{3.20}), (\ref{3.26}) and (\ref{3.45}) we decompose
\begin{equation*}\label{3.57}
\aligned
\int_{\Omega_p}a(\varepsilon y)\widetilde{Z}_{i0}\mathcal{L}(\widetilde{Z}_{i0})
=J+K,
\endaligned
\end{equation*}
where
$$
\aligned
J=&\int_{\Omega_p}a(\varepsilon y)\widetilde{Z}_{i0}\left[
\eta_{i1}\mathcal{L}(Z_{i0}-\widehat{Z}_{i0})
+\eta_{i2}\mathcal{L}(\widehat{Z}_{i0})
\right]\\[0.5mm]
=&\int_{\Omega_p}a(\varepsilon y)\eta_{i2}^2\left\{Z_{i0}-\big(1-\eta_{i1}\big)
\left[\frac1{\mu_i}-a_{i0}G(\varepsilon y,\xi_i)\right]\right\}\times\left\{
W_{\xi'}\big(1-\eta_{i1}\big)\left[\frac1{\mu_i}-a_{i0}G(\varepsilon y,\xi_i)\right]
\right.\\[1mm]
&\left.
-\left[\Delta_{a(\varepsilon y)}Z_{i0}+\frac{8\mu_i^2}{(\mu_i^2+|y-\xi'_i|^2)^2}Z_{i0}
\right]+\left[\frac{8\mu_i^2}{(\mu_i^2+|y-\xi'_i|^2)^2}-W_{\xi'}\right]Z_{i0}
+\varepsilon^2\left(
Z_{i0}-\frac{1}{\mu_i}\right)
+\frac{\varepsilon^2}{\mu_i}\eta_{i1}
\right\},
\endaligned
$$
and
$$
\aligned
K=\int_{\Omega_p}a(\varepsilon y)\widetilde{Z}_{i0}\left[
-(Z_{i0}-\widehat{Z}_{i0})\Delta_{a(\varepsilon y)}\eta_{i1}-2\nabla\eta_{i1}\nabla(Z_{i0}-\widehat{Z}_{i0})
-2\nabla\eta_{i2}\nabla\widehat{Z}_{i0}-\widehat{Z}_{i0}\Delta_{a(\varepsilon y)}\eta_{i2}
\right].
\endaligned
$$
Let us analyze first the behavior of $J$.
By (\ref{2.51}), (\ref{2.52}),
(\ref{3.5}), (\ref{3.8}), (\ref{3.22}) and  (\ref{3.47}),
$$
\aligned
&\int_{|z_i|\leq\mu_i^{1/2}\big/\big(3\varepsilon^{1/2} p^{2\kappa}\big)}
a(\varepsilon y)\eta_{i2}^2\left\{Z_{i0}-(1-\eta_{i1})
\left[\frac1{\mu_i}-a_{i0}G(\varepsilon y,\xi_i)\right]\right\}\times
W_{\xi'}\big(1-\eta_{i1}\big)\left[\frac1{\mu_i}-a_{i0}G(\varepsilon y,\xi_i)\right]
dy\\[0.5mm]
=&O\left(\frac{1}{p\mu_i^2R}\right).
\endaligned
$$
By (\ref{3.5}), (\ref{3.8}), (\ref{3.9}), (\ref{3.22}) and  (\ref{3.47}),
$$
\aligned
&\int_{|z_i|\leq\mu_i^{1/2}\big/\big(3\varepsilon^{1/2} p^{2\kappa}\big)}
\,a(\varepsilon y)\eta_{i2}^2\left\{Z_{i0}-(1-\eta_{i1})
\left[\frac1{\mu_i}-a_{i0}G(\varepsilon y,\xi_i)\right]\right\}\times
\left[\Delta_{a(\varepsilon y)}Z_{i0}+\frac{8\mu_i^2}{(\mu_i^2+|y-\xi'_i|^2)^2}Z_{i0}
\right]dy\\[0.5mm]
=&O\left(\frac{\varepsilon}{\mu_i}\right).
\endaligned
$$
By (\ref{3.5}), (\ref{3.8}),  (\ref{3.22}),  (\ref{3.46}) and (\ref{3.47}),
$$
\aligned
&\int_{|z_i|\leq\mu_i^{1/2}\big/\big(3\varepsilon^{1/2} p^{2\kappa}\big)}
\,a(\varepsilon y)\eta_{i2}^2\left\{Z_{i0}-(1-\eta_{i1})
\left[\frac1{\mu_i}-a_{i0}G(\varepsilon y,\xi_i)\right]\right\}\times
\left[\varepsilon^2\left(
Z_{i0}-\frac{1}{\mu_i}\right)
+\frac{\varepsilon^2}{\mu_i}\eta_{i1}\right]dy\\[1mm]
=&O\left(p\varepsilon^2\right).
\endaligned
$$
By  (\ref{2.52}),  (\ref{3.5}), (\ref{3.8}), (\ref{3.22}) and (\ref{3.47}),
$$
\aligned
&\int_{|z_i|\leq\mu_i^{1/2}\big/\big(3\varepsilon^{1/2} p^{2\kappa}\big)}
\,a(\varepsilon y)\eta_{i2}^2(1-\eta_{i1})
\left[\frac1{\mu_i}-a_{i0}G(\varepsilon y,\xi_i)\right]
\left[\frac{8\mu_i^2}{(\mu_i^2+|y-\xi'_i|^2)^2}-W_{\xi'}\right]Z_{i0}dy\\[1mm]
=&\int_{\mu_iR<|z_i|\leq\mu_i^{1/2}\big/\big(3\varepsilon^{1/2} p^{2\kappa}\big)}
\,\frac{8\mu^2_ia(\varepsilon y)}{(\mu_{i}^2+|y-\xi'_i|^2)^2}
\left[\,\frac1p\left(\omega_1-U_{1,0}-\frac12U_{1,0}^2
\right)\left(\frac{y-\xi'_i}{\mu_i}\right)+
O\left(\frac1{p^2}\log^4\left|\frac{y-\xi'_i}{\mu_i}\right|\right)\right]\\[1.5mm]
&\,\,\,\,
\times
O\left(\frac{\log|y-\xi_i'|-\log(\mu_iR)+
\varepsilon^\alpha|y-\xi_i'|^\alpha}{p\mu_i^2}\right)dy
\\[0.5mm]
=&\,O\left(
\frac{1}{p^2\mu_i^2R}
\right).
\endaligned
$$
While  by  (\ref{2.51}),  (\ref{2.52}), (\ref{3.23}), (\ref{3.26})  and (\ref{3.53}),

$$
\aligned
&\int_{|z_i|>\mu_i^{1/2}\big/\big(3\varepsilon^{1/2} p^{2\kappa}\big)}
\,a(\varepsilon y)\widetilde{Z}_{i0}\left[
\eta_{i1}\mathcal{L}(Z_{i0}-\widehat{Z}_{i0})
+\eta_{i2}\mathcal{L}(\widehat{Z}_{i0})
\right]dy\\[0.5mm]
=&\sum_{k\neq i}^m\int_{\Omega_{3,k}}
a(\varepsilon y)\eta_{i2}^2\widehat{Z}_{i0}\mathcal{L}(\widehat{Z}_{i0})dy
+
\int_{\widetilde{\Omega}_{3}\cup\Omega_4}
a(\varepsilon y)\eta_{i2}^2\widehat{Z}_{i0}\mathcal{L}(\widehat{Z}_{i0})dy
\\
=&\sum_{k\neq i}^mO\left(\int_{0}^{\mu_k^{1/2}\big/\big(3\varepsilon^{1/2} p^{2\kappa}\big)}
\frac{8\mu^2_k}{(\mu_{k}^2+r^2)^2}
\frac{\log^2 p}{p^2\mu_i^2}rdr
\right)+O
\left(\int_{\mu_i^{1/2}\big/\big(3\varepsilon^{1/2} p^{2\kappa}\big)}^{6d/\varepsilon}
\frac{8\mu_i^2}{(\mu_i^2+r^2)^2}
\frac{\log\varepsilon r}{p\mu_i^2}rdr
\right)\\[0.5mm]
=&\,O\left(\frac{\log^2p}{p^2\mu^2_i}\right).
\endaligned
$$
Thus
$$
\aligned
J=\int_{|z_i|\leq\mu_i^{1/2}\big/\big(3\varepsilon^{1/2} p^{2\kappa}\big)}
\,a(\varepsilon y)\eta_{i2}^2Z_{i0}^2
\left[\frac{8\mu_i^2}{(\mu_i^2+|y-\xi'_i|^2)^2}-W_{\xi'}\right]dy
+O\left(\frac{1}{p\mu_i^2R}\right).
\endaligned
$$
Note that in a straightforward but tedious way, by (\ref{2.10}) we can compute
$$
\aligned
\int_{\mathbb{R}^2}e^{U_{1,0}(z)}
Z_0^2(z)\left(
\omega_1-U_{1,0}-\frac12U_{1,0}^2
\right)(z)dz=-8\pi.
\endaligned
$$
Hence by (\ref{2.14}), (\ref{2.52}), (\ref{3.5}), (\ref{3.8}) and (\ref{3.43}),
\begin{equation*}\label{3.58}
\aligned
J=\frac{c_i a(\xi_i)}{p\mu_i^2}\left[1+O\left(\frac1{R}\right)\right],
\quad\forall\,\,i=1,\ldots,m.
\endaligned
\end{equation*}

Next, we estimate $K$.
Integrating by parts the first term and the last term of $K$ respectively, we get
$$
\aligned
K=&-\int_{\Omega_2}a(\varepsilon y)Z_{i0}\nabla\eta_{i1}\nabla(Z_{i0}-\widehat{Z}_{i0})
+\int_{\Omega_2}a(\varepsilon y)\big(Z_{i0}-\widehat{Z}_{i0}\big)\nabla\eta_{i1}\nabla(Z_{i0}-\widehat{Z}_{i0})
\\
&+\int_{\Omega_2}a(\varepsilon y)(Z_{i0}-\widehat{Z}_{i0})^2|\nabla\eta_{i1}|^2+\int_{\Omega_2}a(\varepsilon y)(Z_{i0}-\widehat{Z}_{i0})\nabla\eta_{i1}\nabla\widehat{Z}_{i0}
+\int_{\Omega_4}a(\varepsilon y)
|\widehat{Z}_{i0}|^2|\nabla\eta_{i2}|^2
\\[1.2mm]
\equiv&K_{21}+K_{22}+K_{23}+K_{24}+K_4.
\endaligned
$$
From (\ref{2.14}), (\ref{3.2}), (\ref{3.5}), (\ref{3.8}), (\ref{3.20}), (\ref{3.21}), (\ref{3.42}), (\ref{3.43}) and (\ref{3.47}) we can conclude that
$$
\aligned
K_{21}=&-\frac{a_{i0}}{\mu_i^2}
\int_{\{\mu_iR<|z_i|\leq\mu_i(R+1)\}}
\frac{1}{|y-\xi'_i|}a(\varepsilon y)Z_0\left(\frac{z_i}{\mu_i}\right)
\eta_1'\left(\frac{|z_i|}{\mu_i}\right)
\left[\frac{4}{c_i}+o(1)\right]dy
\\[1.2mm]
=&-\frac{c_i a_{i0}}{4\mu_i}\int_{R}^{R+1}
a(\xi_i)\eta_1'(r)\left[
\frac{4}{c_i}+O\left(\frac1{r^2}\right)
\right]dr\\[1mm]
=&\frac{c_i a(\xi_i)}{p\mu_i^2}\left[1+O\left(\frac1{R^2}\right)\right].
\endaligned
$$
Moreover, by (\ref{3.2}), (\ref{3.5}), (\ref{3.8}), (\ref{3.42}) and (\ref{3.50})  we find
 $|\nabla\widehat{Z}_{i0}|=O\big(\frac1{\mu_i^2 R^3}\big)$ in $\Omega_2$. Furthermore,
$$
\aligned
K_{22}=O\left(\frac{1}{p^2\mu_i^2R}\right),
\,\,\quad\quad\quad\,\,
K_{23}=O\left(\frac{1}{p^2\mu_i^2R}\right)
\,\,\quad\quad\quad\,\,\textrm{and}\,\,\quad\quad\quad\,\,
K_{24}=O\left(\frac{1}{p\mu^2_iR^3}\right).
\endaligned
$$
By (\ref{3.54}),
$$
\aligned
K_4=O\left(\frac{|\log d|^2}{p^2\mu_i^2}\right).
\endaligned
$$
Combining all these estimates, we conclude
\begin{equation}\label{3.59}
\aligned
\int_{\Omega_p}a(\varepsilon y)\widetilde{Z}_{i0}\mathcal{L}(\widetilde{Z}_{i0})
=J+K=\frac{2c_i a(\xi_i)}{p\mu_i^2}\left[1+O\left(\frac1{R}\right)\right],
\quad\forall\,\,i=1,\ldots,m.
\endaligned
\end{equation}

According to (\ref{3.56}), we need just to
 consider
$\int_{\Omega_p}a(\varepsilon y)\widetilde{Z}_{k0}\mathcal{L}(\widetilde{Z}_{i0})$
when $k\neq i$.
Using the previous estimates of
$\mathcal{L}(\widetilde{Z}_{i0})$  and $\widetilde{Z}_{k0}$, we can easily prove that
$$
\aligned
\int_{\Omega_1}a(\varepsilon y)\widetilde{Z}_{k0}\mathcal{L}(\widetilde{Z}_{i0})
=O\left(\frac{R^2
\log p}{p^2\mu_i\mu_k}\right),
\,\qquad\quad\quad\,\,
\int_{\Omega_2}a(\varepsilon y)\widetilde{Z}_{k0}\mathcal{L}(\widetilde{Z}_{i0})
=O\left(\frac{\log p}{p^2\mu_i\mu_k}\right),
\endaligned
$$
$$
\aligned
\int_{\Omega_4}a(\varepsilon y)\widetilde{Z}_{k0}\mathcal{L}(\widetilde{Z}_{i0})=
O\left(\frac{|\log d|^2}{p^2\mu_i\mu_k}\right),
\,\,\,\,\quad\qquad\,\,\,\,
\int_{\Omega_{3,i}\cup\widetilde{\Omega}_{3}}a(\varepsilon y)\widetilde{Z}_{k0}\mathcal{L}(\widetilde{Z}_{i0})
=O\left(\frac{\log p}{p^2\mu_i\mu_k}\right),
\endaligned
$$
and
$$
\aligned
\int_{\Omega_{3,l}}a(\varepsilon y)\widetilde{Z}_{k0}\mathcal{L}(\widetilde{Z}_{i0})
=O\left(\frac{\log^2 p}{p^2\mu_i\mu_k}\right)
\quad\quad\textrm{for all}\,\,\,l\neq i\,\,\,\textrm{and}\,\,\,l\neq k.
\endaligned
$$
It remains to calculate the integral over $\Omega_{3,k}$. From (\ref{3.26}) and  an integration by parts  we get
$$
\aligned
\int_{\Omega_{3,k}}a(\varepsilon y)\widetilde{Z}_{k0}\mathcal{L}(\widetilde{Z}_{i0})
=\int_{\Omega_{3,k}}a(\varepsilon y)\widetilde{Z}_{i0}\mathcal{L}(\widetilde{Z}_{k0})
-\int_{\partial\Omega_{3,k}}a(\varepsilon y)\widehat{Z}_{k0}\frac{\partial\widehat{Z}_{i0}}{\partial\nu}
+\int_{\partial\Omega_{3,k}}a(\varepsilon y)\widehat{Z}_{i0}\frac{\partial\widehat{Z}_{k0}}{\partial\nu}.
\endaligned
$$
As above, we know that
$$
\aligned
\int_{\Omega_{3,k}}a(\varepsilon y)\widetilde{Z}_{i0}\mathcal{L}(\widetilde{Z}_{k0})
=O\left(\frac{\log p}{p^2\mu_i\mu_k}\right).
\endaligned
$$
On $\partial\Omega_{3,k}$, by (\ref{2.3}) and (\ref{3.23})  we have
$$
\aligned
\,\,\,\,\,\,\,\widehat{Z}_{k0}
=O\left(\frac{1}{\mu_k}\right),
\,\,\qquad\,\qquad\qquad\,\qquad\,\,
\widehat{Z}_{i0}
=O\left(\frac{\log p}{p\mu_i}\right),
\quad
\endaligned
$$
and
$$
\aligned
|\nabla\widehat{Z}_{k0}|
=O\left(\frac{\varepsilon^{1/2}p^{2\kappa-1}}{\mu_k^{3/2}}\right),
\,\qquad\,\,\,\qquad\,
|\nabla\widehat{Z}_{i0}|
=O\left(\frac{\,\varepsilon p^{\kappa-1}\,}{\mu_i}\right).
\endaligned
$$
So
$$
\aligned
\int_{\Omega_{3,k}}a(\varepsilon y)\widetilde{Z}_{k0}\mathcal{L}(\widetilde{Z}_{i0})
=O\left(\frac{\log p}{p^2\mu_i\mu_k}\right).
\endaligned
$$
By the above estimates,  we obtain
\begin{equation}\label{3.60}
\aligned
\int_{\Omega_p}a(\varepsilon y)\widetilde{Z}_{k0}\mathcal{L}(\widetilde{Z}_{i0})
=O\left(\frac{\log^2p}{p^2\mu_i\mu_k}
\right),\,\,\quad\,\,\textrm{if}\,\,\,\,i\neq k.
\endaligned
\end{equation}
As a consequence,  replacing (\ref{3.59}) and (\ref{3.60})
to (\ref{3.56}) we get
\begin{equation*}\label{}
\aligned
\frac{|d_i|}{\mu_i}
\leq
Cp\|h\|_{*}
+C\sum_{k=1}^m\frac{|d_k|}{\mu_k}\frac{\log^2p}{p}.
\endaligned
\end{equation*}
Using linear algebra arguments, we then prove Claim 2 for $d_i$ and complete the proof by (\ref{3.32}).
\end{proof}

\vspace{1mm}
\vspace{1mm}
\vspace{1mm}

{\bf Step 4:}
Proof of Proposition 3.3. We first establish the validity of the a priori estimate
\begin{equation}\label{3.62}
\aligned
\|\phi\|_{L^{\infty}(\Omega_p)}\leq Cp\|h\|_{*}
\endaligned
\end{equation}
for any $\phi$, $c_{ij}$ solutions of
problem (\ref{3.1}) and any $h\in C^{0,\alpha}(\overline{\Omega}_p)$.
Step 3 gives
$$
\aligned
\|\phi\|_{L^{\infty}(\Omega_p)}\leq Cp\left(
\|h\|_{*}+\sum_{i=1}^m\sum_{j=1}^{J_i}|c_{ij}|\cdot\|\chi_iZ_{ij}\|_*
\right)\leq Cp\left(
\|h\|_{*}+\sum_{i=1}^m\sum_{j=1}^{J_i}\mu_i|c_{ij}|
\right).
\endaligned
$$
As before, arguing by contradiction to (\ref{3.62}), we
can proceed as in Step 2 and suppose further that
\begin{equation}\label{3.63}
\aligned
\|\phi_n\|_{L^{\infty}(\Omega_{p_n})}=1,
\,\,\ \,\quad\,\,p_n\|h_n\|_{*}\rightarrow0,
\,\,\ \,\quad\,\,
p_n\sum_{i=1}^m\sum_{j=1}^{J_i}\mu_i^n|c_{ij}^n|\geq\delta>0
\,\,\,\quad\,\,\,
\textrm{as}\,\,\,n\rightarrow+\infty.
\endaligned
\end{equation}
We omit the dependence on $n$.
It suffices to estimate the values of the constants $c_{ij}$.
For this aim, let  us consider the cut-off function $\eta_{i2}$
given by (\ref{3.24})-(\ref{3.25}). For any
$i=1,\ldots,m$ and $j=1,J_i$,
multiplying  (\ref{3.1}) by $a(\varepsilon y)\eta_{i2}Z_{ij}$
and integrating by parts we find
\begin{equation}\label{3.64}
\aligned
\int_{\Omega_p}
a(\varepsilon y)\phi \mathcal{L}(\eta_{i2}Z_{ij})=
\int_{\Omega_p}
a(\varepsilon y)h\eta_{i2}Z_{ij}+
\sum_{k=1}^m\sum_{t=1}^{J_k}c_{kt}
\int_{\Omega_p}\chi_kZ_{kt}\eta_{i2}Z_{ij}.
\endaligned
\end{equation}
Notice that
$$
\aligned
\mathcal{L}(\eta_{i2}Z_{ij})=&
\eta_{i2}\mathcal{L}(Z_{ij})
-Z_{ij}\Delta_{a(\varepsilon y)}\eta_{i2}
-2\nabla\eta_{i2}\nabla Z_{ij}
\\[1mm]
=&\eta_{i2}\left[-\Delta_{a(\varepsilon y)}Z_{ij}
-\frac{8\mu_i^2}{(\mu_{i}^2+|y-\xi'_i|^2)^2}Z_{ij}\right]
+\left[\frac{8\mu_i^2}{(\mu_{i}^2+|y-\xi'_i|^2)^2}-W_{\xi'}\right]\eta_{i2}Z_{ij}
-\varepsilon^2\eta_{i2}Z_{ij}
+O\left(\varepsilon^3/d^3\right)\\[1mm]
\equiv&B_1+B_2
-\varepsilon^2\eta_{i2}Z_{ij}
+O\left(\varepsilon^3/d^3\right).
\endaligned
$$
By (\ref{3.2}), (\ref{3.5}), (\ref{3.8}), (\ref{3.42}) and (\ref{3.43}) we can compute
$$
\aligned
B_1=O\left(
\frac{\varepsilon}{(|y-\xi'_i|+\mu_i)^2}
\right).
\endaligned
$$
To estimate $B_2$, we decompose  $\supp(\eta_{i2})$ into several regions:
$$
\aligned
\widehat{\Omega}_{1k}=\left\{y\in\Omega_p\big|\,|z_k|\leq
\frac{\mu_k^{1/2}}{3\varepsilon^{1/2} p^{2\kappa}}\right\},
\,\,\ \,\,\forall\,\,k=1,\ldots,m,
\endaligned
$$
$$
\aligned
\widehat{\Omega}_{2}=\left\{y\in\Omega_p\big|\,
|z_i|\leq\frac{6d}{\varepsilon},
\,\quad\,|z_k|>
\frac{\mu_k^{1/2}}{3\varepsilon^{1/2} p^{2\kappa}},
\,\,\,\,\,\,\,\,\,k=1,\ldots,m\right\},
\endaligned
$$
where $z_k=y-\xi'_k$  for $k=1,\ldots,l$, but
$z_k=F_k^p(y)$  for $k=l+1,\ldots,m$.
Note that, by (\ref{2.3}) and (\ref{3.43}),
\begin{equation}\label{3.65}
\aligned
|y-\xi'_i|\geq |\xi_i'-\xi_k'|-|y-\xi_k'|\geq
|\xi_i'-\xi_k'|-\frac{\mu_k^{1/2}}{\varepsilon^{1/2} p^{2\kappa}}
\geq
\frac{1}{\varepsilon p^\kappa}\left[1-\frac{(\varepsilon\mu_k)^{1/2}}{p^\kappa}
\right],
\endaligned
\end{equation}
uniformly in $\widehat{\Omega}_{1k}$, $k\neq i$.
In $\widehat{\Omega}_{1i}$, by (\ref{2.52}),
(\ref{3.5}) and (\ref{3.8})
we have that for any $i=1,\ldots,l$ and $j=1,2$,
$$
\aligned
B_2=-
\frac{8\mu^2_i(y-\xi'_i)_j}{(\mu_{i}^2+|y-\xi'_i|^2)^{3}}
\left\{\frac1p\left(\omega_1-U_{1,0}-\frac12U_{1,0}^2
\right)\left(\frac{y-\xi'_i}{\mu_i}\right)
+O\left(\frac{\log^4\big(\big|\frac{y-\xi'_i}{\mu_i}\big|+2\big)}{p^2}\right)\right\},
\endaligned
$$
and for any $i=l+1,\ldots,m$ and $j=J_i=1$,
$$
\aligned
B_2=-
\frac{8\mu^2_i}{(\mu_{i}^2+|y-\xi'_i|^2)^{2}}
\frac{(F_i^p(y))_1}{\mu_i^2+|F_i^p(y)|^2}\left\{\frac1p\left(\omega_1-U_{1,0}-\frac12U_{1,0}^2
\right)\left(\frac{y-\xi'_i}{\mu_i}\right)
+O\left(\frac{\log^4\big(\big|\frac{y-\xi'_i}{\mu_i}\big|+2\big)}{p^2}\right)\right\}.
\endaligned
$$
In $\widehat{\Omega}_{1k}$, $k\neq i$, by (\ref{2.52}),  (\ref{3.5}), (\ref{3.8}), (\ref{3.43}) and  (\ref{3.65}),
$$
\aligned
B_2=
O\left(\frac{8\varepsilon p^{\kappa}\mu^2_k}
{(\mu_{k}^2+|y-\xi'_k|^2)^{2}}\right).
\endaligned
$$
In $\widehat{\Omega}_{2}$, by (\ref{2.51}),
$$
\aligned
B_2=
\sum_{k=1}^mO\left(\frac{\mu_k^2}{|y-\xi'_k|^4}\cdot\frac{\varepsilon^{1/2} p^{2\kappa}}{\mu_i^{1/2}}
\right).
\endaligned
$$
Hence by (\ref{2.1}), (\ref{2.8}) and (\ref{3.43}),
\begin{equation}\label{3.66}
\aligned
\int_{\Omega_p}
a(\varepsilon y)\phi \mathcal{L}(\eta_{i2}Z_{ij})=-\frac{1}{p\mu_i}a(\xi_i)E_{j}(\widehat{\phi}_i)
+O\left(\frac{1}{p^2\mu_i}\|\phi\|_{L^{\infty}(\Omega_p)}\right),
\endaligned
\end{equation}
where for any $i=1,\ldots,l$ and $j=1,2$, $\widehat{\phi}_i(z)=\phi\big(\xi_i'+\mu_iz\big)$ and
$$
\aligned
E_{j}(\widehat{\phi}_i)=\int_{B_{\frac1{3\varepsilon^{1/2}\mu_i^{1/2} p^{2\kappa}}}\left(0\right)}
\frac{8z_j}
{(|z|^{2}+1)^3}\left(\omega_1-U_{1,0}-\frac12U_{1,0}^2\right)(|z|)\widehat{\phi}_i(z)
dz,
\endaligned
$$
but for any $i=l+1,\ldots,m$ and $j=1$,
$\widehat{\phi}_i(z)=\phi\big((F_i^p)^{-1}(\mu_iz)\big)$ and
$$
\aligned
E_{j}(\widehat{\phi}_i)=\int_{\mathbb{R}_+^2\bigcap B_{\frac1{3\varepsilon^{1/2}\mu_i^{1/2} p^{2\kappa}}}\left(0\right)}
\frac{8z_j}
{(|z|^{2}+1)^3}\left(\omega_1-U_{1,0}-\frac12U_{1,0}^2\right)(|z|)\widehat{\phi}_i(z)
dz.
\endaligned
$$
On the other hand, since
$\|\eta_{i2}Z_{ij}\|_{L^{\infty}(\Omega_p)}\leq C\mu_i^{-1}$, we obtain
\begin{equation}\label{3.67}
\aligned
\int_{\Omega_p}a(\varepsilon y)h\eta_{i2}Z_{ij}=O\left(\frac{1}{\mu_i}\|h\|_{*}\right).
\endaligned
\end{equation}
Moreover,  if $1\leq k=i\leq l$, by (\ref{3.2}) and (\ref{3.5}),
\begin{equation}\label{3.68}
\aligned
\int_{\Omega_p}\chi_kZ_{kt}\eta_{i2}Z_{ij}=
\int_{\mathbb{R}^2}\chi Z_{t} Z_{j}dz=D_t\delta_{tj},
\endaligned
\end{equation}
and if $l+1\leq k=i\leq m$, by  (\ref{3.2}), (\ref{3.7}) and (\ref{3.8}),
\begin{equation}\label{3.69}
\aligned
\int_{\Omega_p}\chi_kZ_{k1}\eta_{i2}Z_{i1}=
\int_{\mathbb{R}_{+}^2}\chi Z^2_{1} \big[1+O\left(\varepsilon\mu_i|z|\right)\big] dz=\frac12D_1\big[1+O\left(\varepsilon\mu_i\right)\big],
\endaligned
\end{equation}
and  if $k\neq i$, by (\ref{3.65}),
\begin{equation}\label{3.70}
\aligned
\int_{\Omega_p}\chi_kZ_{kt}\eta_{i2}Z_{ij}=O\left(\mu_k\varepsilon p^\kappa\right).
\endaligned
\end{equation}
Inserting estimates (\ref{3.66})-(\ref{3.70}) into  (\ref{3.64}),  we deduce that
for any $i=1,\ldots,m$ and $j=1,J_i$,
$$
\aligned
D_jc_{ij}+O\left(\sum\limits_{k\neq i}^m\sum\limits_{t=1}^{J_k}\mu_k\varepsilon p^\kappa|c_{kt}|
\right)
= O\left(
\frac{1}{\mu_i}\|h\|_{*}
+\frac1{p\mu_i}\|\phi\|_{L^{\infty}(\Omega_p)}
\right).
\endaligned
$$
Furthermore,
\begin{equation}\label{3.71}
\aligned
\sum_{i=1}^m\sum_{j=1}^{J_i}\mu_i|c_{ij}|
= O\left(
\|h\|_{*}+\frac1{p}\|\phi\|_{L^{\infty}(\Omega_p)}
\right).
\endaligned
\end{equation}
Since $\sum\limits_{i=1}^m\sum\limits_{j=1}^{J_i}\mu_i|c_{ij}|=o\left(1\right)$,
as in contradiction arguments of
Step 2, we conclude  that for any $i=1,\ldots,l$,
$$
\aligned
\widehat{\phi}_i\rightarrow
C_i\frac{|z|^2-1}{|z|^2+1}
\,\,\ \,\,\,
\,\,\textrm{uniformly in}\,\,\,C_{loc}^0(\mathbb{R}^2),
\endaligned
$$
but for any $i=l+1,\ldots,m$,
$$
\aligned
\widehat{\phi}_i\rightarrow
C_i\frac{|z|^2-1}{|z|^2+1}
\,\,\ \,\,\,
\,\,\textrm{uniformly in}\,\,\,C_{loc}^0(\mathbb{R}_{+}^2),
\endaligned
$$
with some constant $C_i\in\mathbb{R}$.
Hence in (\ref{3.66}),
we have a better  estimate, since  by
 Lebesgue's theorem we can derive that for
any $i=1,\ldots,l$ and $j=1,2$,
$$
\aligned
E_{j}(\widehat{\phi}_i)
\longrightarrow
C_i\int_{\mathbb{R}^2}\frac{8z_j}{(|z|^2+1)^3}\frac{|z|^2-1}{|z|^2+1}
\left(\omega_1-U_{1,0}-\frac12U_{1,0}^2\right)(|z|)dz=0,
\endaligned
$$
and for
any $i=l+1,\ldots,m$ and $j=1$,
$$
\aligned
E_{j}(\widehat{\phi}_i)
\longrightarrow
C_i\int_{\mathbb{R}^2_{+}}\frac{8z_j}{(|z|^2+1)^3}\frac{|z|^2-1}{|z|^2+1}
\left(\omega_1-U_{1,0}-\frac12U_{1,0}^2\right)(|z|)dz=0.
\endaligned
$$
Therefore,
$$
\aligned
\sum_{i=1}^m\sum_{j=1}^{J_i}\mu_i|c_{ij}|=o(\frac1p)+O(\|h\|_{*}),
\endaligned
$$
which is impossible because of  (\ref{3.63}). So estimate (\ref{3.62})
is established and then by (\ref{3.71}),  we obtain
\begin{equation*}\label{3.72}
\aligned
|c_{ij}|\leq C\frac{1}{\mu_i}\|h\|_{*}.
\endaligned
\end{equation*}
\indent Now consider the Hilbert space
$$
\aligned
H_{\xi}=\left\{\phi\in H^1(\Omega_p)\left|\,\,
\int_{\Omega_p}\chi_iZ_{ij}\phi=0
\quad
\textrm{for any}\,\,\,i=1,\ldots,m,\,\,j=1,J_i;
\quad
\frac{\partial\phi}{\partial\nu}=0\quad\textrm{on}\,\,\,\po_p
\right.
\right\}
\endaligned
$$
with the norm
$\|\phi\|_{H^1(\Omega_p)}^2=\int_{\Omega_p}a(\varepsilon y)\big(|\nabla\phi|^2+\varepsilon^2\phi^2\big)$.
Equation (\ref{3.1}) is equivalent to find
$\phi\in H_\xi$ such that
$$
\aligned
\int_{\Omega_p}a(\varepsilon y)\big(\nabla\phi\nabla\psi+\varepsilon^2\phi\psi\big)
-\int_{\Omega_p}a(\varepsilon y)W_{\xi'}\phi\psi
=\int_{\Omega_p}a(\varepsilon  y)h\psi\,\,\quad\,\,\forall
\,\,\psi\in H_\xi.
\endaligned
$$
By Fredholm's alternative this is equivalent to the uniqueness of solutions to this
problem, which is guaranteed by estimate (\ref{3.62}).
Finally, for  $p\geq p_m$  fixed,
by density of $C^{0,\alpha}(\overline{\Omega}_p)$
in $(C(\overline{\Omega}_p),\,\|\cdot\|_{L^{\infty}(\Omega_p)})$,
we can approximate $h\in C(\overline{\Omega}_p)$
by smooth functions and, by (\ref{3.62}) and
elliptic regularity theory, we find that
for any $h\in C(\overline{\Omega}_p)$,
problem (\ref{3.1}) admits a unique solution which
belongs to $H^2(\Omega_p)$ and satisfies the a priori estimate
(\ref{3.10}).
The proof is complete.
\qquad\qquad\qquad\qquad\qquad\qquad\qquad\qquad\qquad\qquad\qquad\quad
\qquad\qquad\qquad\qquad\qquad\qquad\qquad\qquad\qquad\qquad\qquad\,$\square$

\vspace{1mm}
\vspace{1mm}
\vspace{1mm}
\vspace{1mm}

\noindent{\bf Remark 3.7.}\,\,Given $h\in C(\overline{\Omega}_p)$ with
$\|h\|_*<\infty$, let $\phi$ be the solution to (\ref{3.1}) given by Proposition 3.3.
Testing the first equation of (\ref{3.1}) against $a(\varepsilon y)\phi$, we get
$$
\aligned
\|\phi\|_{H^1(\Omega_p)}^2
=\int_{\Omega_p}a(\varepsilon y)W_{\xi'}\phi^2
+\int_{\Omega_p}a(\varepsilon  y)h\phi.
\endaligned
$$
Furthermore, by (\ref{2.51}) we obtain
$$
\aligned
\|\phi\|_{H^1(\Omega_p)}\leq C\big(\|h\|_*+\|\phi\|_{L^{\infty}(\Omega_p)}\big).
\endaligned
$$

\vspace{1mm}
\vspace{1mm}
\vspace{1mm}

\section{The nonlinear problem}
In order to solve problem (\ref{2.36}) we
 first  consider an auxiliary nonlinear problem: for any points
$\xi=(\xi_1,\ldots,\xi_m)\in\mathcal{O}_p$,
we find a function $\phi$ and scalars $c_{ij}$,
$i=1,\ldots,m$, $j=1,J_i$  such that
\begin{equation}\label{4.1}
\left\{\aligned
&\mathcal{L}(\phi)=-\big[
R_{\xi'}+N(\phi)
\big]
+\frac1{a(\varepsilon y)}\sum\limits_{i=1}^m\sum\limits_{j=1}^{J_i}c_{ij}\chi_i\,Z_{ij}
\,\,\ \,\textrm{in}\,\,\,\,\,\,\Omega_p,\\
&\frac{\partial\phi}{\partial\nu}=0\,\,\,\,\,\,\,
\qquad\qquad\qquad\quad\qquad\qquad\qquad\qquad
\qquad
\textrm{on}\,\,\,\,\partial\Omega_{p},\\[1mm]
&\int_{\Omega_p}\chi_i\,Z_{ij}\phi=0
\qquad\qquad\qquad\quad\,
\forall\,\,i=1,\ldots,m,\,\,\,j=1, J_i,
\endaligned\right.
\end{equation}
where $W_{\xi'}=pV_{\xi'}^{p-1}$ satisfies (\ref{2.51})-(\ref{2.52}),
and $R_{\xi'}$, $N(\phi)$ are given by
(\ref{2.35}).

\vspace{1mm}
\vspace{1mm}
\vspace{1mm}
\vspace{1mm}

\noindent{\bf Proposition 4.1.}\,\,{\it
Let $m$ be a positive integer.
Then there exist constants $C>0$ and  $p_m>1$ such
that for any $p>p_m$ and any points
$\xi=(\xi_1,\ldots,\xi_m)\in\mathcal{O}_p$,
problem {\upshape(\ref{4.1})} admits
a unique solution
$\phi_{\xi'}$
for some coefficients $c_{ij}(\xi')$,
$i=1,\ldots,m$, $j=1,J_i$, such that
\begin{equation}\label{4.2}
\aligned
\|\phi_{\xi'}\|_{L^{\infty}(\Omega_p)}\leq\frac{C}{p^3},\,\,\quad\,\quad\,\,
\sum_{i=1}^m\sum_{j=1}^{J_i}\mu_i|c_{ij}(\xi')|\leq\frac{C}{p^4}\quad\,\quad\textrm{and}\quad\,\quad
\|\phi_{\xi'}\|_{H^1(\Omega_p)}\leq\frac{C}{p^3}.
\endaligned
\end{equation}
Furthermore, the map $\xi'\mapsto\phi_{\xi'}$ is a $C^1$-function in $C(\overline{\Omega}_p)$ and $H^1(\Omega_p)$.
}

\vspace{1mm}

\begin{proof}
Proposition $3.3$, Remarks $3.2$ and $3.7$ allow us to apply the Contraction Mapping Theorem
and the Implicit Function Theorem
to find a solution for problem (\ref{4.1})
satisfying (\ref{4.2}) and the corresponding regularity
of the map $\xi'\mapsto\phi_{\xi'}$.
Since it is a standard procedure, we omit the detailed
proof here, see Lemma 4.1 in \cite{EMP} for a similar proof.
\end{proof}

\vspace{1mm}
\vspace{1mm}

\noindent{\bf Remark 4.2.}\,\,The function $V_{\xi'}+\phi_{\xi'}$, where $\phi_{\xi'}$ is given by Proposition 4.1, is positive in
$\overline{\Omega}_p$. In fact, we  observe first
that $p^2\phi_{\xi'}\rightarrow0$ uniformly in $C(\overline{\Omega}_p)$.
Furthermore, from Remark 2.6 and the definition of $V_{\xi'}$ in
(\ref{2.31}) we have that, in
each region $|y-\xi'_i|<1/(\varepsilon p^{2\kappa})$, $V_{\xi'}+\phi_{\xi'}$ is positive.
Outside these regions, by (\ref{2.20}) and (\ref{2.31}) we may get the same result.

\vspace{1mm}

\section{Variational reduction}
After problem (\ref{4.1}) has been solved, we find a solution of problem (\ref{2.36})
and hence to the original problem  (\ref{1.1}) if we find  $\xi'$ such that
the coefficient $c_{ij}(\xi')$ in (\ref{4.1}) satisfies
\begin{equation}\label{5.1}
\aligned
c_{ij}(\xi')=0\,\quad\,\,\textrm{for all}\,\,\,i=1,\ldots,m,\,\,\,j=1,J_i.
\endaligned
\end{equation}
\indent
Equation (\ref{1.1}) is the Euler-Lagrange equation of the energy functional
$J_p$ given by
\begin{equation*}\label{5.2}
\aligned
J_p(u)=\frac1{2}\int_{\Omega}a(x)(|\nabla
u|^2+u^2)dx-\frac1{p+1}\int_{\Omega}a(x)u_{+}^{p+1}dx,
\,\,\quad\,\,\,\,\,u\in H^1(\Omega).
\endaligned
\end{equation*}
We define
\begin{equation}\label{5.3}
\aligned
F_p(\xi)=J_p(U_{\xi}+\widetilde{\phi}_{\xi})\,\ \,\,\ \ \,\forall\,\,\xi\in\mathcal{O}_p,
\endaligned
\end{equation}
where $U_\xi$ is the function defined in (\ref{2.12}) and
\begin{equation*}\label{5.4}
\aligned
\widetilde{\phi}_{\xi}(x)=\varepsilon^{-2/(p-1)}\phi_{\xi'}(\varepsilon^{-1}x),
\,\quad\,x\in\Omega,
\endaligned
\end{equation*}
with $\phi_{\xi'}$ the unique solution to problem (\ref{4.1}) given by
Proposition 4.1. Critical points of $F_p$ correspond to solutions of (\ref{5.1})
for large $p$, as the following results states.

\vspace{1mm}
\vspace{1mm}
\vspace{1mm}
\vspace{1mm}

\noindent{\bf Proposition 5.1.}\,\,{\it The function $F_p:\mathcal{O}_p\mapsto\mathbb{R}$
is of class $C^1$.
Moreover, for all $p$ sufficiently large,
if  $D_{\xi}F_p(\xi)=0$,  then $\xi'=\xi/\varepsilon$ satisfies {\upshape (\ref{5.1})}.
}

\vspace{1mm}

\begin{proof}
A direct consequence of the results obtained in Proposition 4.1 and the
definition of function $U_{\xi}$ is the fact
that $F_p(\xi)$ is a $C^1$-function of $\xi$ in $\mathcal{O}_p$
since the map
$\xi\mapsto\widetilde{\phi}_{\xi}$ is a $C^1$-map into $H^1(\Omega)$.
Recall that
$$
\aligned
I_p(\upsilon)=\frac12\int_{\Omega_p}a(\varepsilon y)\left(
|\nabla \upsilon|^2+\varepsilon^2\upsilon^2
\right)dy-\frac1{p+1}\int_{\Omega_p}a(\varepsilon y)\upsilon^{p+1}_{+}dy,
\,\quad\,\upsilon\in H^1(\Omega_p),
\endaligned
$$
Then, making a change of variable, we get
\begin{equation}\label{5.5}
\aligned
F_p(\xi)=J_p\big(U_{\xi}+\widetilde{\phi}_{\xi}\big)=\varepsilon^{-4/(p-1)}I_p\big(V_{\xi'}+\phi_{\xi'}\big).
\endaligned
\end{equation}
Since  $\phi_{\xi'}$ solves equation (\ref{4.1})
and $D_{\xi}F_p(\xi)=0$,  we have that for any $k=1,\ldots,m$ and $t=1,J_k$,
\begin{eqnarray}\label{5.6}
&&0
=I'_p\big(V_{\xi'}+\phi_{\xi'}\big)\partial_{(\xi'_k)_t}\big(V_{\xi'}+\phi_{\xi'}\big)\nonumber\\[1mm]
&&\,\,\,\,\,=
\int_{\Omega_p}a(x)\left[
-\Delta_{a(\varepsilon y)}\big(V_{\xi'}+\phi_{\xi'}\big)+\varepsilon^2\big(V_{\xi'}+\phi_{\xi'}\big)-\big(V_{\xi'}+\phi_{\xi'}\big)^p
\right]\partial_{(\xi'_k)_t}\big(V_{\xi'}+\phi_{\xi'}\big)
\nonumber\\
&&\,\,\,\,\,=\sum\limits_{i=1}^m\sum\limits_{j=1}^{J_i}c_{ij}(\xi')\int_{\Omega_p}\chi_i Z_{ij}
\partial_{(\xi'_k)_t}\big(V_{\xi'}+\phi_{\xi'}\big)\nonumber\\
&&\,\,\,\,\,=\sum\limits_{i=1}^m\sum\limits_{j=1}^{J_i}c_{ij}(\xi')\int_{\Omega_p}\chi_i Z_{ij}
\partial_{(\xi'_k)_t}V_{\xi'}
-\sum\limits_{i=1}^m\sum\limits_{j=1}^{J_i}c_{ij}(\xi')\int_{\Omega_p}\phi_{\xi'}
\partial_{(\xi'_k)_t}\big(\chi_i Z_{ij}\big),
\end{eqnarray}
because  $\int_{\Omega_p}\chi_iZ_{ij}\phi_{\xi'}=0$.
Notice first that  by (\ref{3.2}), (\ref{3.5}),  (\ref{3.8})  and (\ref{3.43}),
a direct computation shows
\begin{eqnarray*}\label{5.7}
\big|\partial_{(\xi'_k)_t}\big(\chi_i Z_{ij}\big)\big|
=O\left(
\frac1{\mu_i}\big|\partial_{(\xi'_k)_t}\log\mu_i\big|
+\frac1{\mu_i^2}\delta_{ki}
\right).
\end{eqnarray*}
On the other hand,  since  $D_{\xi'}V_{\xi'}(y)=\varepsilon^{2/(p-1)}D_{\xi'}U_\xi(\varepsilon y)$,
by (\ref{2.4}), (\ref{2.5}), (\ref{2.12}), (\ref{2.20}) and (\ref{2.21})
we obtain
$$
\aligned
\partial_{(\xi'_k)_t}V_{\xi'}(y)
=&\,\sum_{i=1}^{m}
\frac{\varepsilon^{2/(p-1)}}
{\,\,\gamma\mu_i^{2/(p-1)}\,\,}
\left\{\partial_{(\xi'_k)_t}\left[U_{\delta_i,\xi_i}(\varepsilon y)
+\frac1p\omega_1\left(\frac{y-\xi'_i}{\mu_i}\right)
+\frac1{p^2}\omega_2\left(\frac{y-\xi'_i}{\mu_i}\right)
+\gamma\mu_i^{2/(p-1)}H_i(\varepsilon y)\right]\right.\\
&\left.
-
\frac{\,2\,\partial_{(\xi'_k)_t}\log\mu_i\,}{p-1}\big[\,
p+O\left(1\right)\big]
\right\}.
\endaligned
$$
By (\ref{2.1}) and (\ref{3.2}),
$$
\aligned
\partial_{(\xi'_k)_t}U_{\delta_i,\xi_i}(\varepsilon y)=&
\frac4{\mu_i}\delta_{ki}
Z_{t}\left(
\frac{y-\xi'_i}{\mu_i}
\right)
+2Z_{0}\left(
\frac{y-\xi'_i}{\mu_i}
\right)
\partial_{(\xi'_k)_t}\log\mu_i\\
=&\frac4{\mu_i}\delta_{ki}
Z_{t}\left(
\frac{y-\xi'_i}{\mu_i}
\right)
+2\left(
1-\frac{2\mu_i^2}{|y-\xi'_i|^2+\mu_i^2}
\right)\partial_{(\xi'_k)_t}\log\mu_i,
\endaligned
$$
while for $j=1,2$,
by (\ref{2.8}) and (\ref{3.2}),
$$
\aligned
\partial_{(\xi'_k)_t}\omega_{j}\left(\frac{y-\xi'_i}{\mu_i}\right)=&-
\omega_{j}'\left(\frac{|y-\xi'_i|}{\mu_i}\right)\left\{
\frac{(y-\xi'_i)_t}{|y-\xi'_i|}\frac{1}{\mu_i}\delta_{ki}
+\frac{|y-\xi'_i|}{\mu_i}\partial_{(\xi'_k)_t}\log\mu_i
\right\}\\
=&-\frac1{\mu_i}\delta_{ki}\left[C_jZ_{t}\left(
\frac{y-\xi'_i}{\mu_i}
\right)+
O\left(\frac{\mu_i^2}{|y-\xi'_i|^2+\mu_i^2}\right)
\right]
-C_j\left[1+
O\left(\frac{\mu_i^2}{|y-\xi'_i|^2+\mu_i^2}\right)
\right]\partial_{(\xi'_k)_t}\log\mu_i.
\endaligned
$$
Additionally,  as in the proof of Lemma 2.1, we can  prove that
$$
\aligned
\partial_{(\xi'_k)_t}\left[\,\gamma\mu_i^{2/(p-1)}H_i(\varepsilon y)\right]=&\,
\delta_{ki}\left(
1-\frac{C_1}{4p}-\frac{C_2}{4p^2}
\right)c_i\partial_{(\xi'_k)_t} H(\varepsilon y,\xi_i)
+\left(-2+
\frac{C_1}{p}+\frac{C_2}{p^2}
\right)\partial_{(\xi'_k)_t}\log\mu_i\\
&+
O\left(\frac1p\big|\partial_{(\xi'_k)_t}\log\mu_i\big|\right).
\endaligned
$$
Thus by (\ref{2.4}) and (\ref{2.23}),
\begin{eqnarray*}\label{5.8}
\partial_{(\xi'_k)_t}V_{\xi'}(y)=\frac{1}{p^{p/(p-1)}\mu_k^{(p+1)/(p-1)}}
\left\{
\left(
1-\frac{C_1}{4p}-\frac{C_2}{4p^2}
\right)\left[
4Z_{t}\left(
\frac{y-\xi'_k}{\mu_k}
\right)+c_k\mu_k\partial_{(\xi'_k)_t} H(\varepsilon y,\xi_k)
\right]\right.
&&\nonumber\\
\left.+
O\left(
\frac1p\frac{\mu_k^2}{|y-\xi'_k|^2+\mu_k^2}\right)
\right\}+
O\left(
\frac1p\sum_{i=1}^m\big|\partial_{(\xi'_k)_t}\log\mu_i\big|\right).
\qquad\quad\qquad\qquad\qquad\,\,\,\,\,\quad
&&
\end{eqnarray*}
Hence for
each $k=1,\ldots,m$  and  $t=1,J_k$,   (\ref{5.6}) can be written as
$$
\aligned
&\sum_{i,j}\frac{c_{ij}(\xi')}{p^{p/(p-1)}\mu_k^{2/(p-1)}}\left\{
\delta_{ki}\left[
\frac{c_i}{2\pi}\int_{\mathbb{R}^2}\chi(|z|)Z_j(z)Z_t(z)dz+O\left(\frac1p\right)
\right]
+(1-\delta_{ki})O\left(
\frac{\mu_i}{|\xi'_i-\xi'_k|}
\right)
\right\}\\
&+\sum_{ij}|c_{ij}(\xi')|
\left\{
O\left(
\frac{1}{p}\mu_i\sum_{s=1}^m
\big|\partial_{(\xi'_k)_t}\log\mu_s\big|
\right)+\|\phi_{\xi'}\|_{L^{\infty}(\Omega_p)}\,
O\left(\mu_i\big|\partial_{(\xi'_k)_t}\log\mu_i\big|
+\delta_{ki}
\right)
\right\}=0,
\endaligned
$$
and then, by (\ref{2.3}), (\ref{2.23}), (\ref{2.60}) and (\ref{4.2}),
$$
\aligned
\frac{c_kc_{kt}(\xi')}{p^{p/(p-1)}\mu_k^{2/(p-1)}}
\int_0^{R_0+1}\chi(r)\frac{r^{3}}{(\,1+r^{2})^2}
dr
+
\sum_{i=1}^m\sum_{j=1}^{J_i}\big|c_{ij}(\xi')\big|O\left(
\frac{\delta_{ki}}{p^2\mu_k^{2/(p-1)}}
+
\varepsilon\mu_i p^{\kappa-1}
\right)=0,
\endaligned
$$
which  implies $c_{kt}(\xi')=0$
for each $k=1,\ldots,m$  and  $t=1,J_k$.
\end{proof}

\vspace{1mm}

\section{Expansion of the energy}
\noindent{\bf Proposition 6.1.}\,\,{\it
Let $m$ be a positive integer.
With the choice of $\mu_i$'s given by {\upshape(\ref{2.22})},
there exists $p_m>1$ such that for any $p>p_m$ and any
points $\xi=(\xi_1,\ldots,\xi_m)\in\mathcal{O}_p$,
the following expansion  holds
\begin{equation}\label{6.1}
\aligned
F_p(\xi)=\frac{e}{2p}
\sum_{i=1}^mc_ia(\xi_i)\left[
1-\frac{2\log p}{p}
+\frac{\mathcal{K}+2}p-\frac{1}{p}c_i H(\xi_i,\xi_i)
-\frac{1}{p}\sum_{k\neq i}^m
c_k G(\xi_i,\xi_k)
\right]
+O\left(\frac{\log^2p}{p^3}\right),
\endaligned
\end{equation}
where
$$
\aligned
\mathcal{K}=\frac{1}{8\pi}\int_{\mathbb{R}^2}\left[
\frac{8}{(1+|z|^2)^2}U_{1,0}(z)-\Delta\omega_1(z)
\right]dz.
\endaligned
$$
}

\begin{proof}
Multiplying the first equation in (\ref{4.1}) by $a(\varepsilon y)(V_{\xi'}+\phi_{\xi'})$ and integrating by parts, we obtain
$$
\aligned
\int_{\Omega_p}a(\varepsilon y)\left[|\nabla
(V_{\xi'}+\phi_{\xi'})|^2+\varepsilon^2(V_{\xi'}+\phi_{\xi'})^2\right]dy&=\int_{\Omega_p}a(\varepsilon y)(V_{\xi'}+\phi_{\xi'})^{p+1}dy
+\sum_{i=1}^m\sum_{j=1}^{J_i}c_{ij}(\xi')\int_{\Omega_p}\chi_iZ_{ij}V_{\xi'}dy\\
&=\int_{\Omega_p}a(\varepsilon y)(V_{\xi'}+\phi_{\xi'})^{p+1}dy
+O\left(\frac1{p^4}\right)
\endaligned
$$
because $V_{\xi'}=O\left(1\right)$ and $|c_{ij}(\xi')|=O\left(p^{-4}\mu_i^{-1}\right)$.
Hence, by (\ref{2.31}), (\ref{4.2}) and (\ref{5.5})  we can write
\begin{eqnarray*}
F_p(\xi)=\left(\frac1{2}-\frac{1}{p+1}\right)\varepsilon^{-4/(p-1)}\int_{\Omega_p}a(\varepsilon y)\left[|\nabla
(V_{\xi'}+\phi_{\xi'})|^2+\varepsilon^2(V_{\xi'}+\phi_{\xi'})^2\right]dy+O\left(\frac1{p^5}\right)
\,\,\,\,\,\quad
&&\nonumber\\
=\left(\frac1{2}-\frac{1}{p+1}\right)\varepsilon^{-4/(p-1)}\left\{\int_{\Omega_p}a(\varepsilon y)\big(|\nabla
V_{\xi'}|^2+\varepsilon^2V_{\xi'}^2\big)dy
\right.
\qquad\qquad\qquad\quad\,\,
\qquad\qquad\,\,\,\quad\,
&&\nonumber\\
\left.+2\int_{\Omega_p}a(\varepsilon y)\big(\nabla
V_{\xi'}\nabla\phi_{\xi'}+\varepsilon^2V_{\xi'}\phi_{\xi'}\big)dy
+\int_{\Omega_p}a(\varepsilon y)\big(|\nabla
\phi_{\xi'}|^2+\varepsilon^2\phi_{\xi'}^2\big)dy\right\}+O\left(\frac1{p^5}\right)
\quad
&&\nonumber\\
=\left(\frac1{2}-\frac{1}{p+1}\right)\int_{\Omega}a(x)\big(|\nabla
U_{\xi}|^2+U_{\xi}^2\big)dx
+O\left(\frac1{p^3}\left|\int_{\Omega}a(x)\big(|\nabla
U_{\xi}|^2+U_{\xi}^2\big)dx\right|^{1/2}+\frac1{p^5}\right)
&&
\end{eqnarray*}
uniformly for any points
$\xi=(\xi_1,\ldots,\xi_m)\in\mathcal{O}_p$.
Now,  in view of (\ref{2.4}), (\ref{2.5}), (\ref{2.12}),
(\ref{2.13}), (\ref{2.14}), (\ref{2.20}) and (\ref{2.21}) we have
$$
\aligned
\int_{\Omega}a(x)&\big(|\nabla
U_{\xi}|^2+U_{\xi}^2\big)dx=\int_{\Omega}a(x)(-\Delta_{a}
U_\xi+U_{\xi})U_\xi dx\\
&=\sum_{i=1}^m\frac{1}{\gamma\mu_i^{2/(p-1)}}
\int_{\Omega\bigcap B_{\frac1{p^{2\kappa}}}(\xi_i)}
a(x)\left[e^{U_{\delta_i,\xi_i}}
-\frac{1}{p\delta_i^2}\Delta\omega_1\left(
\frac{x-\xi_i}{\delta_i}\right)
-\frac{1}{p^2\delta_i^2}\Delta\omega_2\left(
\frac{x-\xi_i}{\delta_i}\right)
\right]U_\xi dx +o\left(\delta_i\right)\\
&=\sum_{i=1}^m\frac{1}{\gamma^2\mu_i^{4/(p-1)}}
\int_{\big(\frac{1}{\delta_i}(\Omega-\{\xi_i\})\big)\bigcap B_{\frac1{\delta_ip^{2\kappa}}}(0)}
a(\xi_i+\delta_iz)\left[\frac{8}{(1+|z|^2)^2}
-\frac{1}{p}\Delta\omega_1(z)
-\frac{1}{p^2}\Delta\omega_2(z)
\right]\\
&\,\,\,\qquad\,\times\left[
p+U_{1,0}(z)+\frac1p\omega_1(z)+\frac1{p^2}\omega_2(z)+O\left(
\delta_i^\alpha|z|^\alpha
+\sum_{k=1}^m\delta_k^{\alpha/2}
\right)\right]dz+o\left(\delta_i\right)\\
&=\sum_{i=1}^m\frac{c_ia(\xi_i)}{\gamma^2\mu_i^{4/(p-1)}}
\left[\,p
+\mathcal{K}
+O\left(\frac1p\right)
\right]+o\left(\delta_i\right).
\endaligned
$$
Recalling  that $\gamma=p^{p/(p-1)}e^{-p/(2p-2)}$, we get
$$
\aligned
\frac{1}{\gamma^2}=\frac{e}{p^2}\left[
1-\frac{2\log p}{p}+\frac{1}{p}+O\left(
\frac{\log^2 p}{p^2}
\right)
\right],
\,\qquad\qquad\,
\frac{1}{\,\mu_i^{4/(p-1)}}=1-\frac{4\log\mu_i}{p}+O\left(\frac{\log^2\mu_i}{p^2}\right),
\endaligned
$$
and then, by (\ref{2.23}),
$$
\aligned
\frac{1}{\,\gamma^2\mu_i^{4/(p-1)}\,}=
\frac{e}{p^2}\left[
1-\frac{2\log p}{p}-\frac{4\log\mu_i}{p}
+\frac{1}{p}+O\left(
\frac{\log^2p}{p^2}
\right)
\right].
\endaligned
$$
Therefore,
$$
\aligned
F_p(\xi)=\frac{e}{2p}
\sum_{i=1}^mc_ia(\xi_i)\left\{
1-\frac{2\log p}{p}-\frac{4\log\mu_i}{p}
+\frac{\mathcal{K}-1}p
\right\}
+O\left(\frac{\log^2p}{p^3}\right),
\endaligned
$$
which, together with the expansion of $\mu_i$ in
(\ref{2.28}), easily implies that (\ref{6.1}) holds.
\end{proof}

\vspace{1mm}

\section{Proofs of theorems}
\noindent {\bf Proof of Theorem 1.1.}
According to Proposition  5.1,
$U_{\xi}+\widetilde{\phi}_{\xi}$
is a solution of problem (\ref{1.1}) if we adjust
$\xi=(\xi_1,\ldots,\xi_m)\in\mathcal{O}_p$
so that it is
a critical point of $F_p$ defined in (\ref{5.3}).
For this aim, we only choose points $\xi=(\xi_1,\ldots,\xi_m)\in\mathcal{O}_p$ in the
following form
of the parametrization
\begin{equation}\label{7.1}
\aligned
\xi_i\equiv\xi_i(\mathbf{s},\mathbf{t})=s_i-\frac{t_i}{p}\nu(s_i),
\quad i=1,\ldots,l,
\qquad\quad\textrm{but}\quad\qquad
\xi_i\equiv\xi_i(\mathbf{s},\mathbf{t})=s_i,
\quad i=l+1,\ldots,m,
\endaligned
\end{equation}
where $\mathbf{s}=(s_1,\ldots,s_m)\in(\po)^m$ and $\mathbf{t}=(t_1,\ldots,t_l)\in\mathbb{R}^l_{+}$
lie in the configuration space
$$
\aligned
\Lambda_d=\left\{
(\mathbf{s},\,\mathbf{t})\in(\po)^m\times\mathbb{R}^l_{+}\,
\big|\,|s_i-s_k|>2d,
\quad\,d<t_{\tilde{i}}<1/d,
\quad\forall\,\,\,i, k=1,\ldots,m,\,\,\,\tilde{i}=1,\ldots,l,\,\,\,i\neq k
\right\},
\endaligned
$$
for any $d>0$ sufficiently small, fixed and independent of $p$.
Thus we can easily prove that if $(\mathbf{s},\mathbf{t})$
is a critical point of the reduced energy
$\widehat{F}_p\big(\mathbf{s},\mathbf{t}\big):=F_p\big(\xi(\mathbf{s},\mathbf{t})\big)$
in $\Lambda_d$, then the function
$U_{\xi(\mathbf{s},\mathbf{t})}+\widetilde{\phi}_{\xi(\mathbf{s},\mathbf{t})}$
is a solution of problem (\ref{1.1}) with
the qualitative properties  predicted by Theorem 1.1.
Therefore, we  need first to compute the expansion of the reduced energy $\widehat{F}_p$ with the aid of
Lemmas 2.1-2.2, Corollary 2.3 and Proposition 6.1.

Using the smooth property of $a(x)$ over $\oo$,  we can perform
a Taylor expansion around each boundary point $s_i$ along the inner normal vector
$-\nu(s_i)$ to derive that
\begin{equation}\label{7.2}
\aligned
a(\xi_i)=a(s_i)-\frac{t_i}{p}\partial_{\nu}a(s_i)+O\left(\frac{t_i^2}{p^2}\right),
\quad\forall\,\,
i=1,\ldots,l.
\endaligned
\end{equation}
From the expansions of the Robin's function in (\ref{2.1e})-(\ref{2.1f})
and the regularity of the vector function $T(x)$, we obtain
\begin{equation}\label{7.3}
\aligned
H(\xi_i,\xi_i)=-\frac{1}{2\pi}\log\left(\frac{2t_i}{p}\right)
+\tilde{\mathrm{z}}(s_i,s_i)+O\left(\frac{t_i^\alpha}{p^\alpha}\right),
\quad\forall\,\,
i=1,\ldots,l.
\endaligned
\end{equation}
On the other hand,
from the expansions of the regular part of Green's function in (\ref{2.1c})-(\ref{2.1d}),
we find that if $i, k=1,\ldots,l$ with $i\neq k$,
\begin{eqnarray}\label{7.4}
G(\xi_i,\xi_k)=-\frac{1}{2\pi}\log\left|s_i-s_k-\frac{t_i}{p}\nu(s_i)+\frac{t_k}{p}\nu(s_k)
\right|
-\frac{1}{2\pi}\log\left|s_i-s_k-\frac{t_i}{p}\nu(s_i)-\frac{t_k}{p}\nu(s_k)
\right|
\quad\qquad\quad
&&\nonumber\\[0.2mm]
+\frac{t_k}{\pi p}\left\langle \nabla\log a(s_k),\,\nabla T(s_i-s_k)\cdot\nu(s_k)\right\rangle
-\frac{t_k}{\pi p}\left\langle \big(\nabla\times(\nabla\log a)\big)(s_k)\cdot\nu(s_k),\,T(s_i-s_k)\right\rangle
&&\nonumber\\[0.1mm]
+\tilde{\mathrm{z}}(s_i,s_k)
-\frac{1}{p}\left\langle\nabla_{(s_i,s_k)}\tilde{\mathrm{z}}(s_i,s_k),\,\big(t_i\nu(s_i),\,t_k\nu(s_k)\big)\right\rangle
+O\left(\frac{t_i^2+t_k^2}{p^2}\right),
\qquad\qquad\qquad\quad\,\,\,\,
&&
\end{eqnarray}
while if $i=l+1,\ldots,m$  and $k=1,\ldots,l$,
\begin{eqnarray}\label{7.5}
G(\xi_i,\xi_k)=-\frac{1}{2\pi}\log\left|s_i-s_k+\frac{t_k}{p}\nu(s_k)
\right|
-\frac{1}{2\pi}\log\left|s_i-s_k-\frac{t_k}{p}\nu(s_k)
\right|
\qquad\qquad\qquad\qquad\quad\qquad\quad\,\,\,
&&\nonumber\\[0.2mm]
+\frac{t_k}{\pi p}\left\langle \nabla\log a(s_k),\,\nabla T(s_i-s_k)\cdot\nu(s_k)\right\rangle
-\frac{t_k}{\pi p}\left\langle \big(\nabla\times(\nabla\log a)\big)(s_k)\cdot\nu(s_k),\,T(s_i-s_k)\right\rangle
&&\nonumber\\[0.1mm]
+\tilde{\mathrm{z}}(s_i,s_k)
-\frac{t_k}{p}\left\langle\nabla_{s_k} \tilde{\mathrm{z}}(s_i,s_k),\,\nu(s_k)\right\rangle
+O\left(\frac{t_k^2}{p^2}\right).
\qquad\qquad\qquad\qquad\qquad\qquad\qquad\qquad\,
&&
\end{eqnarray}
Substituting (\ref{2.14}), (\ref{7.2})-(\ref{7.5}) into (\ref{6.1}) and using the fact that
$a(\xi_i)G(\xi_i,\xi_k)=a(\xi_k)G(\xi_k,\xi_i)$ for all
$i,k=1,\ldots,m$ with $i\neq k$, we conclude that $\widehat{F}_p$ becomes
\begin{eqnarray}\label{7.6}
\widehat{F}_p\big(\mathbf{s},\mathbf{t}\big)=
\frac{e}{2p}\left\{8\pi\sum_{i=1}^l\left[a(s_i)
+\frac{1}{p}\large\big(
4a(s_i)\log t_i-t_i\partial_{\nu}a(s_i)
\large\big)
\right]+4\pi\sum_{k=l+1}^ma(s_k)
+\Gamma_p(\mathbf{s})
\right\}+O\left(\frac{\,|\mathbf{t}|^\alpha}{p^{2+\alpha}}+\frac{\,\log^2p}{p^3}\right)
\end{eqnarray}
$C^1$-uniformly in $\Lambda_d$,
where $\Gamma_p(\mathbf{s})$ is a smooth function of $\mathbf{s}$  satisfying that
$\Gamma_p(\mathbf{s})$ and $|\nabla\Gamma_p(\mathbf{s})|$ uniformly converge to
zero as $p\rightarrow+\infty$.

We seek a critical point of  $\widehat{F}_p$ in $\Lambda_d$.
Let $\partial_{T(s_i)}$ be the tangential derivative which is defined on $s_i\in\po$.
Set
$$
\aligned
A(s_i,t_i)=a(s_i)
+\frac{1}{p}\large\big(
4a(s_i)\log t_i-t_i\partial_{\nu}a(s_i)
\large\big),\,\,
\quad\,i=1,\ldots,l.
\endaligned
$$
Then
$$
\aligned
\partial_{T(s_i)}A(s_i,t_i)=\left(1+
\frac{4\log t_i}{p}
\right)\partial_{T}a(s_i)
-\frac{t_i}{p}
\partial_{T}\partial_{\nu}a(s_i)
\large,
\,\qquad\,
\partial_{t_i}A(s_i,t_i)=\frac{1}{p}
\left(
\frac{4}{t_i}a(s_i)-\partial_{\nu}a(s_i)
\right).
\endaligned
$$
In view of $\partial_{\nu}a(\xi^*_i)>0$ with
$i=1,\ldots,l$,
 we have that
for any sufficiently small $d>0$ and for
any $s_i\in B_{d}(\xi^*_i)\cap\po$,
there exists a unique positive $t_i=t_i(s_i)=4a(s_i)/\partial_{\nu}a(s_i)$
such that $\partial_{t_i}A(s_i,t_i)=0$ and $\partial^2_{t_it_i}A(s_i,t_i)<0$.
Set $t_i^*=t_i(\xi_i^*)$, $i=1,\ldots,l$.
Since $\xi^*_1,\ldots,\xi^*_m$ are $m$ different strict local maximum or
strict local minimum points of $a(x)$ on $\po$, we find that for
any $i=1,\ldots,l$, any sufficiently small $d>0$
and any sufficiently large $p>1$,  the Brouwer degrees
$\deg\big(\big(\partial_{T(s_i)}A,\,\partial_{t_i}A\big),\,
\big(B_{d}(\xi^*_i)\cap\po\big)\times\big(t_i^*-d,\,t_i^*+d\big),\,0\big)$
and
$\deg\big(\big(\partial_{T(s_{l+1})}a,\ldots,\partial_{T(s_m)}a\big),
\,\prod_{k=l+1}^m\big(B_{d}(\xi^*_k)\cap\po\big),\,0\big)$
are well defined (see \cite{C,L1}). By the definition and homotopy invariance of
topological degree we deduce
$$
\aligned
\deg\big(\big(\partial_{T(s_i)}A,&\,\partial_{t_i}A\big),\,\big(B_{d}(\xi^*_i)\cap\po\big)\times\big(t_i^*-d,\,t_i^*+d\big),\,0\big)\\[1mm]
=&
\deg\big(\big(\partial_{T}a(s_i),\,\partial_{t_i}A\big),\,\big(B_{d}(\xi^*_i)\cap\po\big)\times\big(t_i^*-d,\,t_i^*+d\big),\,0\big)\\[0.5mm]
=&\sign\,\det\left[
\left(\aligned
    \partial^2_{TT}a(\xi^*_i) \quad &\quad  \frac{1}{p}
\left(
\frac{4}{t^*_i}\partial_{T}a(\xi^*_i)-\partial_{T}\partial_{\nu}a(\xi^*_i)
\right) \\
    0 \quad \quad \quad &\quad \quad \quad  -\frac{1}{p}\frac{4}{(t^*_i)^2}a(\xi^*_i) \\
  \endaligned
\right)
\right]=\pm1
\neq0,
\endaligned
$$
and
$$
\aligned
\deg\left(\big(\partial_{T(s_{l+1})}a,\ldots,\partial_{T(s_m)}a\big),
\,\prod_{k=l+1}^m\big(B_{d}(\xi^*_k)\cap\po\big),\,0\right)
=\sign\left(\prod_{k=l+1}^m\partial^2_{TT}a(\xi^*_k)\right)=\pm1
\neq0.
\endaligned
$$
Furthermore, using the properties of Brouwer degree, by (\ref{7.6}) we conclude
$$
\aligned
&\,\deg\left(\nabla_{\left(T(\mathbf{s}),\mathbf{t}\right)}\widehat{F}_p\big(\mathbf{s},\mathbf{t}\big),\,\,\prod_{i=1}^l\big( B_{d}(\xi^*_i)\cap\po\big)\times\prod_{k=l+1}^m\big(B_{d}(\xi^*_k)\cap\po\big)\times\prod_{i=1}^l\big(t_{i}^*-d,\,t_{i}^*+d\big),\,\,0\right)\\[1mm]
=\,&\deg\left(\big(\partial_{T(s_1)}A,\,\ldots,\,\partial_{T(s_l)}A,\,\,
\partial_{T(s_{l+1})}a,\,\ldots,\,\partial_{T(s_m)}a,\,\,\partial_{t_1}A,\,\ldots,\,\partial_{t_l}A\big),
\right.\\[1mm]
&\left.\qquad\quad\,\,\prod_{i=1}^l\big( B_{d}(\xi^*_i)\cap\po\big)\times\prod_{k=l+1}^m\big(B_{d}(\xi^*_k)\cap\po\big)
\times\prod_{i=1}^l\big(t_{i}^*-d,\,t_{i}^*+d\big),\,\,0\right)\\
=\,&
\prod_{i=1}^l\,\deg\large\big(\big(\partial_{T(s_i)}A,\,\partial_{t_i}A\big),\,\,\big(B_{d}(\xi^*_i)\cap\po\big)
\times\big(t_i^*-d,\,t_i^*+d\big),\,\,0\large\big)\\
&\qquad\quad\,\times
\deg\left(\big(\partial_{T(s_{l+1})}a,\ldots,\partial_{T(s_m)}a\big),
\,\,\prod_{k=l+1}^m\big(B_{d}(\xi^*_k)\cap\po\big),\,\,0\right)\\[1mm]
\neq\,&\,0.
\endaligned
$$
This implies that for any $p>1$ large enough, there exists
$(\mathbf{s}^p,\mathbf{t}^p)\in\Lambda_d$ such that
$\nabla_{\left(T(\mathbf{s}),\mathbf{t}\right)}\widehat{F}_p\big(\mathbf{s}^p,\mathbf{t}^p\big)=0$.
In particular, $\mathbf{s}^p=(s^p_1,\ldots,s^p_m)\rightarrow(\xi^*_1,\ldots,\xi^*_m)$ as $p\rightarrow+\infty$,
which completes the proof.
\qquad\qquad\qquad\qquad\quad\qquad\qquad\qquad\qquad\quad$\square$

\vspace{1mm}
\vspace{1mm}
\vspace{1mm}
\vspace{1mm}

\noindent {\bf Proof of Theorem 1.2.}
According to Proposition 5.1, we need to find
a critical point $\xi^p=(\xi^p_1,\ldots,\xi^p_m)\in\Omega^l\times(\po)^{m-l}$
of  $F_p$ such  that points $\xi^p_1,\ldots,\xi^p_m$
accumulate to $\xi_*$.
Using  (\ref{1.3}), (\ref{2.14}), (\ref{6.1}), Lemma 2.1 and the fact that
$a(\xi_i)G(\xi_i,\xi_k)=a(\xi_k)G(\xi_k,\xi_i)$ for all
$i,k=1,\ldots,m$ with $i\neq k$,
we conclude that $F_p$
becomes
\begin{eqnarray}\label{7.7}
F_p(\xi)=\frac{e}{2p^2}
\left\{8\pi\sum_{i=1}^la(\xi_i)\left[
p-2\log p-8\pi H(\xi_i,\xi_i)
-8\pi \sum_{k=1,\,k\neq i}^l
G(\xi_i,\xi_k)
\right]-64\pi^2\sum_{i=1}^l\sum_{k=l+1}^m a(\xi_k)
G(\xi_k,\xi_i)
\right.
&&\nonumber\\
\left.+4\pi\sum_{i=l+1}^ma(\xi_i)\left[
p-2\log p+4\sum_{k=l+1,\,k\neq i}^m
\log|\xi_i-\xi_k|
\right]
+O\left(1\right)
\right\}
\qquad\qquad\qquad\,\,
\qquad\qquad\qquad
&&
\end{eqnarray}
$C^0$-uniformly in
$\mathcal{O}^*_{d,\,p}:=\left\{\,\xi=(\xi_1,\ldots,\xi_m)\in\big(B_d(\xi_*)\cap\Omega\big)^l\times\big(B_d(\xi_*)\cap\po\big)^{m-l}
\,\left|\,\,\,
\min\limits_{i,k=1,\ldots,m,\,i\neq k}\,|\xi_i-\xi_k|>\frac1{p^{\kappa}},
\right.\right.\\
\left.\min\limits_{1\leq i\leq l}\,\dist(\xi_i,\po)>\frac1{p^{\kappa}}
\right\}$
for any $d>0$ sufficiently small, fixed and independent of $p$.
Here we claim that
for any $p>1$ large enough, the
following maximization problem
$$
\aligned
\max\limits_{(\xi_1,\ldots,\xi_m)\in\overline{\mathcal{O}}^*_{d,p}}
F_p(\xi_1,\ldots,\xi_m)
\endaligned
$$
has a solution in the interior of $\mathcal{O}^*_{d,p}$.
Once this claim  is proven, we can easily get the qualitative properties of
solutions of (\ref{1.1}) as predicted in Theorem 1.3.

Let $\xi^p=(\xi^p_1,\ldots,\xi^p_m)$ be the maximizer of $F_p$
over $\overline{\mathcal{O}}^*_{d,p}$. We are led to prove that $\xi^p$
lies in the interior of $\mathcal{O}^*_{d,p}$.
First, we obtain a lower bound for $F_p$ over $\overline{\mathcal{O}}^*_{d,p}$.
Let us consider a smooth change of variables
$$
\aligned
H_{\xi_*}^p(y)=
e^{p/2}H_{\xi_*}(e^{-p/2} y),
\endaligned
$$
where $H_{\xi_*}:B_d(\xi_*)\mapsto\mathcal{M}$
is a  diffeomorphism and
$\mathcal{M}$  is an open neighborhood of the origin such that
$H_{\xi_*}(B_d(\xi_*)\cap\Omega)=\mathcal{M}\cap\mathbb{R}_+^2$
and
$H_{\xi_*}(B_d(\xi_*)\cap\partial\Omega)=\mathcal{M}\cap\partial\mathbb{R}_+^2$.
Let
$$
\aligned
\xi^0_i=\xi_*-\frac{t_i}{\sqrt{p}}\nu(\xi_*),
\quad i=1,\ldots,l,
\qquad\quad\textrm{but}\quad\qquad
\xi^0_i=e^{-p/2}(H_{\xi_*}^p)^{-1}\left(
\frac{e^{p/2}}{\sqrt{p}}\hat{\xi}_i^0
\right),
\quad i=l+1,\ldots,m,
\endaligned
$$
where $t_i>0$  and
$\hat{\xi}_i^0\in \mathcal{M}\cap\partial\mathbb{R}_+^2$  satisfy
$t_{i+1}-t_i=\rho$,
$|\hat{\xi}_i^0-\hat{\xi}^0_{i+1}|=\rho$
for all $\rho>0$ sufficiently
small, fixed and independent of $p$.
From the expansion
$(H_{\xi_*}^p)^{-1}(z)=e^{p/2}\xi_*+z+O(e^{-p/2}|z|)$,
we get
$$
\aligned
\xi^0_i=\xi_*+\frac{1}{\sqrt{p}}\hat{\xi}_i^0
+O\left(\frac{e^{-p/2}}{\sqrt{p}}|\hat{\xi}_i^0|\right),
\quad i=l+1,\ldots,m.
\endaligned
$$
As $\kappa>1$, we find $\xi^0=(\xi^0_1,\ldots,\xi^0_m)\in\mathcal{O}^*_{d,p}$.
Since $\xi_*\in\po$ is a strict local maximum point of $a(x)$ over $\oo$
and satisfies $\partial_{\nu}a(\xi_*)=\langle\nabla a(\xi_*),\,\nu(\xi_*)\rangle=0$, there exists
a constant $C>0$ independent of $p$ such that
$$
\aligned
a(\xi_*)-\frac{C}{p}\leq a(\xi^0_i)<a(\xi_*),
\qquad i=1,\ldots,m.
\endaligned
$$
On the other hand, from definition (\ref{1.3}),
Lemma 2.2 and Corollary 2.3, we  can compute that
for any $i=1,\ldots,l$ and $k=1,\ldots,m$ with $i\neq k$,
$$
\aligned
H(\xi^0_i,\xi^0_i)=\frac1{4\pi}\log p+O\left(1\right),
\qquad\quad
G(\xi^0_k,\xi^0_i)=H(\xi^0_k,\xi^0_i)
-\frac1{2\pi}\log|\xi^0_k-\xi^0_i|=\frac1{2\pi}\log p+O\left(1\right).
\endaligned
$$
Additionally, for any $i,k=l+1,\ldots,m$ with $i\neq k$,
$$
\aligned
\log|\xi^0_i-\xi^0_k|=-\frac1{2}\log p+O\left(1\right).
\endaligned
$$
Hence by  (\ref{7.7}),   we find
\begin{eqnarray}\label{7.8}
\max\limits_{\xi\in\overline{\mathcal{O}}^*_{d,p}}
F_p(\xi)\geq
F_p(\xi^0)\geq
\frac{e}{2p^2}
\big\{
4\pi (m+l)pa(\xi_*)-8\pi(m+l)^2a(\xi_*)\log p
+O(1)
\big\}.
\end{eqnarray}
Next, we suppose $\xi^p=(\xi^p_1,\ldots,\xi^p_m)\in\partial\mathcal{O}^*_{d,p}$. Then
there exist four cases:\\
C1. \,\,There exists an $i_0\in\{1,\ldots,l\}$ such that
$\xi^p_{i_0}\in\partial B_d(\xi_*)\cap\Omega$, in which case,
$a(\xi^p_{i_0})<a(\xi_*)-d_0$ for some\\
\indent\indent $d_0>0$ independent of $p$;\\
C2. \,\,There exists an $i_0\in\{l+1,\ldots,m\}$ such that
$\xi^p_{i_0}\in\partial B_d(\xi_*)\cap\partial\Omega$, in which case,
$a(\xi^p_{i_0})<a(\xi_*)-d_0$ for\\
\indent\indent some $d_0>0$ independent of $p$;\\
C3. \,\,There exists an $i_0\in\{1,\ldots,l\}$ such that
$\dist(\xi_{i_0}^p,\po)=p^{-\kappa}$;\\
C4. \,\,There exist indices $i_0$, $k_0$, $i_0\neq k_0$ such that
$|\xi_{i_0}^p-\xi_{k_0}^p|=p^{-\kappa}$.\\
Observe that for all $i=1,\ldots,l$ and $k=1,\ldots,m$ with $i\neq k$,
by (\ref{1.2}), (\ref{2.1c}), (\ref{2.1e}) and the maximum principle,
\begin{equation}\label{7.9}
\aligned
G(\xi^p_k,\xi^p_i)>0,
\qquad\,\,
H(\xi^p_k,\xi^p_i)>0
\qquad\,\,
\textrm{and}
\qquad\,\,
H(\xi^p_i,\xi^p_i)=-\frac1{2\pi}\log\big[\dist(\xi_i^p,\po)\big]+O\left(1\right).
\endaligned
\end{equation}

In the first and second cases, by (\ref{7.7}) and (\ref{7.9}) we get
\begin{eqnarray}\label{7.10}
\max\limits_{\xi\in\overline{\mathcal{O}}^*_{d,p}}
F_p(\xi)=
F_p(\xi^p)\leq
\frac{e}{2p^2}
\left\{4\pi p\big[(m+l)a(\xi_*)-d_0\big]+O\left(\log p\right)\right\},
&&
\end{eqnarray}
which contradicts to (\ref{7.8}).
This shows that $a(\xi_i^p)\rightarrow a(\xi_*)$. Using the assumption of $a(x)$ over $\oo$,
we deduce $\xi_i^p\rightarrow \xi_*$ for all $i=1,\ldots,m$.

In the third case, by (\ref{7.7}) and (\ref{7.9}) we get
\begin{eqnarray}\label{7.11}
\max\limits_{\xi\in\overline{\mathcal{O}}^*_{d,p}}
F_p(\xi)=
F_p(\xi^p)
\leq
\frac{e}{2p^2}
\left\{
4\pi(m+l)(p-2\log p)a(\xi_*)
-64\pi^2a(\xi^p_{i_0})H(\xi^p_{i_0},\xi^p_{i_0})
+O(1)
\right\}
&&\nonumber\\
\leq
\frac{e}{2p^2}
\left\{
4\pi(m+l)(p-2\log p)a(\xi_*)
-32\pi\kappa a(\xi^p_{i_0})\log p
+O(1)
\right\}.
\quad\,\,\,\,
&&
\end{eqnarray}

In the last case, by (\ref{7.7}) and (\ref{7.9}) we get, if
$i_0\in\{1,\ldots,m\}$ and $k_0\in\{1,\ldots,l\}$,
\begin{eqnarray}\label{7.12}
\max\limits_{\xi\in\overline{\mathcal{O}}^*_{d,p}}
F_p(\xi)=
F_p(\xi^p)
\leq
\frac{e}{2p^2}
\left\{
4\pi(m+l)(p-2\log p)a(\xi_*)
+32\pi a(\xi^p_{i_0})
\log|\xi^p_{i_0}-\xi^p_{k_0}|
+O(1)
\right\}
&&\nonumber\\
\leq
\frac{e}{2p^2}
\left\{
4\pi(m+l)(p-2\log p)a(\xi_*)
-32\pi\kappa a(\xi^p_{i_0})\log p
+O(1)
\right\},
\qquad\,\,\,\,\,
&&
\end{eqnarray}
while if
$i_0\in\{l+1,\ldots,m\}$ and $k_0\in\{l+1,\ldots,m\}$,
\begin{eqnarray}\label{7.13}
\max\limits_{\xi\in\overline{\mathcal{O}}^*_{d,p}}
F_p(\xi)=
F_p(\xi^p)
\leq
\frac{e}{2p^2}
\left\{
4\pi(m+l)(p-2\log p)a(\xi_*)
+16\pi a(\xi^p_{i_0})
\log|\xi^p_{i_0}-\xi^p_{k_0}|
+O(1)
\right\}
&&\nonumber\\
\leq
\frac{e}{2p^2}
\left\{
4\pi(m+l)(p-2\log p)a(\xi_*)
-16\pi\kappa a(\xi^p_{i_0})\log p
+O(1)
\right\}.
\qquad\,\,\,\,\,
&&
\end{eqnarray}

Comparing (\ref{7.11})-(\ref{7.13}) with  (\ref{7.8}), we obtain
\begin{equation}\label{7.14}
\aligned
8\pi(m+l)a(\xi_*)
\log p+32\pi\kappa a(\xi^p_{i_0})
\log p+O(1)\leq
8\pi(m+l)^2a(\xi_*)
\log p
+O(1),
\endaligned
\end{equation}
which  is impossible by the choice of $\kappa$ in (\ref{2.3a}).
\,\qquad\qquad\qquad\qquad\qquad\qquad\qquad\qquad
\qquad\qquad\qquad\qquad\qquad\qquad\quad$\square$

\vspace{1mm}
\vspace{1mm}
\vspace{1mm}

\section*{Acknowledgments}
The author warmly thanks the anonymous referee for his or her nice and valuable  comments on this
manuscript.
This research is supported by
the Fundamental Research Funds for the
Central Universities under Grant No.  KYZ201649,
and the National
Natural Science Foundation of China  under Grant
Nos. 11601232,   11671354 and  11775116.

\vspace{1mm}

\section{Appendix}

\noindent{\bf Proof of lemma 2.4.}\,\,
Observe that for any $\beta\in(0,1)$,
by (\ref{2.1}), (\ref{2.5}) and (\ref{2.8}),
$$
\aligned
\left\{\aligned&
-\Delta_a
H_i+H_i=\frac{1}{\gamma\mu_i^{2/(p-1)}}\left\{
\left(-4+\frac{C_1}{p}+\frac{C_2}{p^2}
\right)\left[
\frac{(x-\xi_i)\cdot\nabla\log a(x)}{\delta_i^2+|x-\xi_i|^2}
-\frac{1}{2}\log\big(\delta_i^2+|x-\xi_i|^2\big)\right]
-\log(8\delta_i^2)
\right.\\[2mm]
&\left.\,\quad\qquad\qquad\qquad
+\left(\frac{C_1}{p}+\frac{C_2}{p^2}\right)\log\delta_i+\frac1p
O_{\large L^{\infty}\big(\Omega\setminus B_{\delta_i^{\beta/2}}(\xi_i)\big)}
\left(\frac{\delta_i}{\delta_i+|x-\xi_i|}
+\frac{\delta_i}{\delta_i^2+|x-\xi_i|^2}\right)
\right.\\[2mm]
&\left.\,\quad\qquad\qquad\qquad
+\frac1p
O_{\large L^{\infty}\big(\Omega\bigcap B_{\delta_i^{\beta/2}}(\xi_i)\big)}\left(\frac{|(x-\xi_i)\cdot\nabla\log a(x)|}{\delta_i^2+|x-\xi_i|^2}
+\log\frac{\delta_i^2+|x-\xi_i|^2}{\delta_i^2}
\right)
\right\}
\quad\quad\qquad\,
\textrm{in}\,\,\,\,\ \,\,\Omega,\\[2mm]
&\frac{\partial H_i}{\partial \nu}=-\frac{1}{\gamma\mu_i^{2/(p-1)}}\left\{
\left(-4+\frac{C_1}{p}+\frac{C_2}{p^2}
\right)
\frac{(x-\xi_i)\cdot\nu(x)}{\delta_i^2+|x-\xi_i|^2}
+\frac1p
O_{\large L^{\infty}\big(\partial\Omega\setminus B_{\delta_i^{\beta/2}}(\xi_i)\big)}\left(\frac{\delta_i}{\delta_i^2+|x-\xi_i|^2}\right)
\right.\\[2mm]
&\left.\,
\quad\qquad
+\frac1p
O_{\large L^{\infty}\big(\partial\Omega\bigcap B_{\delta_i^{\beta/2}}(\xi_i)\big)}\left(\frac{|(x-\xi_i)\cdot\nu(x)|}{\delta_i^2+|x-\xi_i|^2}\right)
\right\}
\,\qquad\qquad\qquad\qquad\qquad\qquad
\qquad\qquad\qquad\,
\,\textrm{on}\,\,\ \ \,\po.
\endaligned\right.
\endaligned
$$
From (\ref{1.2})-(\ref{1.3}) we have that
the regular part of Green's function, $H(x,\xi_i)$,  satisfies
\begin{equation*}\label{2.16}
\left\{\aligned
&-\Delta_a
H(x,\xi_i)+H(x,\xi_i)=
\frac{4}{c_i}\log|x-\xi_i|-
\frac{4}{c_i}\frac{(x-\xi_i)\cdot\nabla\log a(x)}{|x-\xi_i|^2}
\quad
\,\,\textrm{in}\,\,\,\,\,\Omega,\\[1mm]
&\frac{\partial H(x,\xi_i)}{\partial \nu}
=\frac{4}{c_i}\frac{(x-\xi_i)\cdot\nu(x)}{|x-\xi_i|^2}
\qquad\qquad\qquad\qquad
\qquad\qquad\qquad\qquad\quad
\textrm{on}\,\,\,\po.
\endaligned\right.
\end{equation*}
So, if we set
$$
\aligned
Z_i(x)=\gamma\mu_i^{2/(p-1)}H_i(x)
-\left(
4-\frac{C_1}{p}-\frac{C_2}{p^2}
\right)\frac{c_i}{4}H(x,\xi_i)
+\log(8\delta_i^2)-\left(\frac{C_1}{p}+\frac{C_2}{p^2}\right)\log\delta_i,
\endaligned
$$
then $Z_i(x)$ satisfies
$$
\aligned
\left\{\aligned&
-\Delta_a
Z_i+Z_i=
\left(-4+\frac{C_1}{p}+\frac{C_2}{p^2}
\right)\left[\,\frac{1}{2}\log\left(
\frac{|x-\xi_i|^2}{\delta_i^2+|x-\xi_i|^2}\right)
-\frac{(x-\xi_i)\cdot\nabla\log a(x)}{|x-\xi_i|^2}\cdot\frac{\delta_i^2}{\delta_i^2+|x-\xi_i|^2}
\right]\\[2mm]
&\quad\qquad\qquad\qquad
+\frac1p
O_{\large L^{\infty}\big(\Omega\setminus B_{\delta_i^{\beta/2}}(\xi_i)\big)}
\left(\frac{\delta_i}{\delta_i+|x-\xi_i|}
+\frac{\delta_i}{\delta_i^2+|x-\xi_i|^2}\right)\\[2mm]
&\quad\qquad\qquad\qquad
+\frac1p
O_{\large L^{\infty}\big(\Omega\bigcap B_{\delta_i^{\beta/2}}(\xi_i)\big)}\left(\frac{|(x-\xi_i)\cdot\nabla\log a(x)|}{\delta_i^2+|x-\xi_i|^2}
+\log\frac{\delta_i^2+|x-\xi_i|^2}{\delta_i^2}
\right)
\quad\qquad\,
\textrm{in}\,\,\ \,\,\ \,\,\Omega,\\[2mm]
&\frac{\partial Z_i}{\partial \nu}=
\left(-4+\frac{C_1}{p}+\frac{C_2}{p^2}
\right)
\frac{(x-\xi_i)\cdot\nu(x)}{|x-\xi_i|^2}
\cdot
\frac{\delta_i^2}{\delta_i^2+|x-\xi_i|^2}
+\frac1p
O_{\large L^{\infty}\big(\partial\Omega\setminus B_{\delta_i^{\beta/2}}(\xi_i)\big)}\left(\frac{\delta_i}{\delta_i^2+|x-\xi_i|^2}\right)
\\[2mm]
&\qquad\quad
+\frac1p
O_{\large L^{\infty}\big(\partial\Omega\bigcap B_{\delta_i^{\beta/2}}(\xi_i)\big)}\left(\frac{|(x-\xi_i)\cdot\nu(x)|}{\delta_i^2+|x-\xi_i|^2}\right)
\,\qquad\qquad
\quad\qquad\qquad\quad\qquad
\quad\qquad\qquad\,\,\,\textrm{on}\,\ \,\,\,\po.
\endaligned\right.
\endaligned
$$
Using polar coordinates with center $\xi_i$, i.e. $r=|x-\xi_i|$, and changing variables $s=r/\delta_i$, we
estimate that
for any  $q>1$,
$$
\aligned
\int_{\Omega}\left|
\log
\left(\frac{|x-\xi_i|^2}{\delta_i^2+|x-\xi_i|^2}\right)
\right|^qdx
&\leq2\pi\int_{0}^{\diam(\Omega)}\left|\log
\left(\frac{r^2}{\delta_i^2+r^2}\right)
\right|^qrdr\\
&\leq 2\pi\delta_i^2\int_{0}^{+\infty}\left|\log\left(
1+\frac1{s^2}
\right)\right|^qsds\\[1mm]
&\leq C\delta_i^2,
\endaligned
$$
and
$$
\aligned
\int_{\Omega\setminus B_{\delta_i^{\beta/2}}(\xi_i)}\left|
\frac{\delta_i}{\delta_i+|x-\xi_i|}
+\frac{\delta_i}{\delta_i^2+|x-\xi_i|^2}
\right|^qdx
&\leq C\int_{\delta_i^{\beta/2}}^{\diam(\Omega)}\left[\frac{\delta_i^q}{(\delta_i+r)^q}
+\frac{\delta_i^q}{(\delta_i^2+r^2)^q}
\right]rdr\\
&\leq C\int_{\delta_i^{\beta/2-1}}^{+\infty}\left[\frac{\delta_i^2}{(1+s)^q}
+\frac{\delta_i^{2-q}}{(1+s^2)^q}
\right]sds\\[1mm]
&\leq C\delta_i^{\beta+q(1-\beta)},
\endaligned
$$
and for  any $1<q<2$,
$$
\aligned
\int_{\Omega\bigcap B_{\delta_i^{\beta/2}}(\xi_i)}\left|
\frac{|(x-\xi_i)\cdot\nabla\log a(x)|}{\delta_i^2+|x-\xi_i|^2}
+\log\frac{\delta_i^2+|x-\xi_i|^2}{\delta_i^2}
\right|^qdx
\leq&
C\int_{0}^{\delta_i^{\beta/2}}\left[\left(
\frac{r}{\delta_i^2+r^2}
\right)^q+
\left|\log
\left(\frac{\delta_i^2+r^2}{\delta_i^2}\right)
\right|^q\right]rdr\\
=&C\delta_i^2\int_{0}^{\delta_i^{\beta/2-1}}
\left[\delta_i^{-q}\left(\frac{s}{1+s^2}\right)^q
+\log^q\big(
1+s^2
\big)\right]sds\\[1mm]
\leq& C
\delta_i^{\beta(1-q/2)},
\endaligned
$$
and
$$
\aligned
\int_{\Omega}\left|
\frac{(x-\xi_i)\cdot\nabla\log a(x)}{|x-\xi_i|^2}\cdot\frac{\delta_i^2}{\delta_i^2+|x-\xi_i|^2}
\right|^qdx
&\leq C\int_{0}^{\diam(\Omega)}\left|\frac{\delta_i^2}{r(\delta_i^2+r^2)}
\right|^qrdr\\
&\leq C\delta_i^{2-q}\int_{0}^{+\infty}
\frac{s^{1-q}}{(1+s^2)^q}
ds\\[1mm]
&\leq C\delta_i^{2-q}.
\endaligned
$$
Thus for any $\xi_i\in\oo$ and any  $1<q<2$,
\begin{equation*}\label{2.17}
\aligned
\big\|-\Delta_a
Z_i+Z_i\big\|_{L^q(\Omega)}\leq C\delta_i^{\beta(1/q-1/2)}.
\endaligned
\end{equation*}
As for the boundary terms, if $\xi_i\in\Omega$, by (\ref{2.3}) we get, for any $x\in\partial\Omega$,
$$
\aligned
\left|
\frac{(x-\xi_i)\cdot\nu(x)}{|x-\xi_i|^2}
\cdot
\frac{\delta_i^2}{\delta_i^2+|x-\xi_i|^2}
\right|\leq \frac{\delta_i^2}{|x-\xi_i|^3}\leq \delta_i^2p^{3\kappa},
\endaligned
$$
$$
\aligned
\frac{\delta_i}{\delta_i^2+|x-\xi_i|^2}
\leq \frac{\delta_i}{|x-\xi_i|^2}\leq \delta_ip^{2\kappa},
\endaligned
$$
and further,
\begin{equation*}\label{2.18}
\aligned
\left\|\frac{\partial Z_i}{\partial \nu}\right\|_{L^{\infty}(\partial\Omega)}\leq C\delta_ip^{2\kappa-1}.
\endaligned
\end{equation*}
While if $\xi_i\in\po$, using the fact that
$|(x-\xi_i)\cdot\nu(x)|\leq C|x-\xi_i|^2$ for any $x\in\po$ (see \cite{AP}),
we estimate that for any $q>1$,
$$
\aligned
\int_{\partial\Omega}\left|
\frac{(x-\xi_i)\cdot\nu(x)}{|x-\xi_i|^2}
\cdot
\frac{\delta_i^2}{\delta_i^2+|x-\xi_i|^2}
\right|^qdx
&\leq C\int_{0}^{+\infty}\frac{\delta_i^{2q}}{(\delta_i^2+r^2)^q}
dr=C\delta_i\int_{0}^{+\infty}
\frac{1}{(1+s^2)^q}
ds\leq C\delta_i,
\endaligned
$$
$$
\aligned
\int_{\partial\Omega\bigcap B_{\delta_i^{\beta/2}}(\xi_i)}\left|
\frac{(x-\xi_i)\cdot\nu(x)}{\delta_i^2+|x-\xi_i|^2}
\right|^qdx
&\leq
C\int_{\partial\Omega\bigcap B_{\delta_i^{\beta/2}}(\xi_i)}
\frac{|x-\xi_i|^{2q}}{(\delta_i^2+|x-\xi_i|^2)^q}
dx\leq
C\left|\partial\Omega\cap B_{\delta_i^{\beta/2}}(\xi_i)
\right|\leq C\delta_i^{\beta/2},
\endaligned
$$
and
$$
\aligned
\int_{\partial\Omega\setminus B_{\delta_i^{\beta/2}}(\xi_i)}\left|
\frac{\delta_i}{\delta_i^2+|x-\xi_i|^2}
\right|^qdx
&\leq C\int_{\delta_i^{\beta/2}}^{|\po|}\frac{\delta_i^{q}}{(\delta_i^2+r^2)^q}
dr
\leq
C\int_{\delta_i^{\beta/2}}^{|\po|}
\frac{\delta_i^q}{\,\,r^{2q}\,\,}
dr\leq C\big(\delta_i^q+\delta_i^{q(1-\beta)+\beta/2}\big).
\endaligned
$$
Thus for any $\xi_i\in\po$ and any $q>1$,
\begin{equation*}\label{2.19}
\aligned
\left\|\frac{\partial Z_i}{\partial \nu}\right\|_{L^{q}(\partial\Omega)}\leq C\delta_i^{\beta/2q}.
\endaligned
\end{equation*}
Hence by elliptic regularity theory, we obtain that for any $1<q<2$ and any $0<\theta<1/q$,
$$
\aligned
\left\|Z_i\right\|_{W^{1+\theta,q}(\Omega)}\leq
C\left(\big\|-\Delta_a
Z_i+Z_i\big\|_{L^q(\Omega)}+\left\|\frac{\partial Z_i}{\partial \nu}\right\|_{L^{q}(\partial\Omega)}\right)
\leq C\delta_i^{\beta(1/q-1/2)}.
\endaligned
$$
Then by  Morrey's embedding theorem,
$$
\aligned
\left\|Z_i\right\|_{C^\tau(\overline{\Omega})}\leq C\delta_i^{\beta(1/q-1/2)},
\endaligned
$$
where $0<\tau<1/2+1/q$,
which implies that expansion (\ref{2.15}) holds with $\alpha=2\beta(1/q-1/2)$.
\qquad\qquad\qquad\qquad\qquad$\square$

\vspace{1mm}
\vspace{1mm}
\vspace{1mm}
\vspace{1mm}
\vspace{1mm}

\noindent{\bf Proof of lemma 2.5.}\,\,
Notice that if we make the change of variables  $s=1/p$, then system (\ref{2.22}) can be rewritten
in the following vector form
$$
\aligned
S_i(s,\xi,\mu):=
\log\mu_i+\frac{\,\log8+\frac{1}{4}C_1+\frac{1}{4}C_2s\,}{\,4\big(1-\frac{1}{4}C_1s-\frac{1}{4}C_2s^2\big)\,}-\frac14c_iH(\xi_i,\xi_i)
-\frac14\sum_{k=1,\,k\neq i}^m
\left(
\frac{\,\mu_i\,}{\mu_k}
\right)^{2s/(1-s)}
c_k G(\xi_i,\xi_k)=0.
\endaligned
$$
Obviously, from the explicit expression (\ref{2.11})  of the constant $C_1$, we have that for $s=0$,
\begin{equation}\label{2.24}
\aligned
\mu_i(0,\xi)=e\large^{-\frac{3}{4}+\frac14c_iH(\xi_i,\xi_i)
+\frac14\sum_{k=1,\,k\neq i}^m
c_k G(\xi_i,\xi_k)}.
\endaligned
\end{equation}
Using the Taylor expansion of exponential functions, we can conclude that for any
$s>0$ small enough,
$$
\aligned
\left(
\frac{\,C^2\,}{s^{\kappa}}
\right)^{2s/(1-s)}
=e\large^{\frac{2s}{1-s}\big(2\log C+\kappa\log\frac1s\big)}=1+O\left(
s\log\frac{1}{s}
\right),
\endaligned
$$
and then, by (\ref{2.23}),
\begin{equation}\label{2.25}
\aligned
\left(
\frac{\,\mu_i\,}{\mu_k}
\right)^{2s/(1-s)}
=1+O\left(
s\log\frac{1}{s}
\right),
\qquad\forall\,\,i\neq k.
\endaligned
\end{equation}
Moreover, by (\ref{1.3}), (\ref{2.1c}) and the fact that
$a(\xi_i)G(\xi_i,\xi_k)=a(\xi_k)G(\xi_k,\xi_i)$ for all
$i,k=1,\ldots,m$ with $i\neq k$,
we can easily prove that for any points $\xi=(\xi_1,\ldots,\xi_m)\in\mathcal{O}_{1/s}$,
\begin{equation}\label{2.27}
\aligned
G(\xi_i,\xi_k)=O\left(\log\frac1s\right),
\,\qquad\,\forall\,\,i\neq k.
\endaligned
\end{equation}
Consequently, a simple computation shows that
$$
\aligned
\frac{\partial S_i(s,\xi,\mu)}{\partial\mu_i}=\frac1{\mu_i}-\frac{s}{\,2\mu_i(1-s)\,}\sum_{k=1,\,k\neq i}^m
\left(
\frac{\,\mu_i\,}{\mu_k}
\right)^{2s/(1-s)}
c_k G(\xi_i,\xi_k)=\frac1{\mu_i}\left[
1+O\left(s\log\frac1s\right)
\right],
\endaligned
$$
and for any $k\neq i$,
$$
\aligned
\frac{\partial S_i(s,\xi,\mu)}{\partial\mu_k}=\frac{s}{\,2\mu_k(1-s)\,}
\left(
\frac{\,\mu_i\,}{\mu_k}
\right)^{2s/(1-s)}
c_k G(\xi_i,\xi_k)=\frac1{\mu_k}O\left(s\log\frac1s\right).
\endaligned
$$
Then
$$
\aligned
\det\big(\nabla_{\mu}S(s,\xi,\mu)
\big)=\frac1{\,\mu_1\ldots\mu_m\,}\left[
1+O\left(s\log\frac1s\right)
\right]\neq0.
\endaligned
$$
Hence  $\nabla_{\mu}S(s,\xi,\mu)$ is  invertible in the range of points
and variables that we are considering. From the Implicit Function Theorem we find that
$S(s,\xi,\mu)=0$ is solvable in some neighborhood of $\big(0,\xi,\mu(0,\xi)\big)$,
and thus for   any points
$\xi=(\xi_1,\ldots,\xi_m)\in\mathcal{O}_{p}$
and any $p>1$ large enough, system (\ref{2.22}) has a unique solution $\mu=(\mu_1,\ldots,\mu_m)$
satisfying (\ref{2.23}). This, together with (\ref{2.24})-(\ref{2.27}),
implies
$$
\aligned
\mu_i=e\large^{-\frac{3}{4}+\frac14c_iH(\xi_i,\xi_i)
+\frac14\sum\large_{k=1,\,k\neq i}^m
c_k G(\xi_i,\xi_k)}\left[\,1+O\left(\frac{\log^2p}p\right)\right],
\quad\,\,\forall\,\,i=1,\ldots,m.
\endaligned
$$
Moreover,  by (\ref{2.1c}), (\ref{2.1e}), (\ref{2.3}), (\ref{2.22}), (\ref{2.25}) and (\ref{2.27})
we can conclude that estimate (\ref{2.60}) holds.\qquad\qquad\qquad\qquad$ \square$

\end{document}